\documentclass[11pt, leqno]{amsart}
\usepackage{amsmath, amssymb, latexsym}
\usepackage{mathrsfs}
\usepackage{enumerate}
\usepackage{color}
\usepackage[colorlinks=true, pdfstartview=FitV, linkcolor=blue,%
citecolor=blue, urlcolor=blue]{hyperref}
\usepackage[all, knot]{xy}

\numberwithin{equation}{section}

\newtheorem{theorem}{Theorem}[section]
\newtheorem{corollary}[theorem]{Corollary}
\newtheorem{lemma}[theorem]{Lemma}
\newtheorem{proposition}[theorem]{Proposition}

\theoremstyle{definition}
\newtheorem{definition}[theorem]{Definition}
\newtheorem{remark}[theorem]{Remark}

\newcommand{\Hom}{\operatorname{Hom}}

\newcommand{\A}{\mathcal{A}}
\newcommand{\B}{\mathscr{B}_{q}(\g)}
\newcommand{\C}{\mathbf{C}}

\newcommand{\Q}{\mathbf{Q}}
\newcommand{\N}{\mathbf{N}}
\newcommand{\Z}{\mathbf{Z}}
\newcommand{\h}{\mathfrak{h}}
\newcommand{\g}{\mathfrak{g}}

\newcommand{\U}{U_{q}(\mathfrak{g})}

\title[Crystal bases for quantum Borcherds-Bozec superalgebra]
{Crystal bases for quantum Borcherds-Bozec superalgebra}

\author[Zhaobing Fan]{Zhaobing Fan}
\address{Harbin Engineering University,
	Harbin, China}
\email{fanzhaobing@hrbeu.edu.cn}

\author[Jiaqi Huang]{Jiaqi Huang}
\address{Harbin Engineering University,
	Harbin, China}
\email{jiaqihuang@hrbeu.edu.cn}

\author[Ziqi Xin]{Ziqi Xin}
\address{Harbin Engineering University,
	Harbin, China}
\thanks{}



\begin{document}

\begin{abstract}
Let $U_q^{-}(\mathfrak{g})$ denote the negative half of quantum Borcherds--Bozec superalgebra with isotropic vertex, and let $V(\lambda)$ be the integrable irreducible highest weight module. In this paper, we first modify the Kashiwara operators and define crystal bases for $U_q^{-}(\mathfrak{g})$ and $V(\lambda)$. We then give the tensor product rule and establish the existence of crystal bases. Finally, we construct the global bases for $U_q^{-}(\mathfrak{g})$ and $V(\lambda)$.
\end{abstract}

\maketitle

\section*{Introduction}
In \cite{Bor88}, R. E. Borcherds introduced {\it Borcherds algebras} associated with Borcherds-Cartan matrices.
The diagonal entries of these matrices may be non-positive, which means the Borcherds algebras
may have simple roots with the norm $\le 0$, the {\it imaginary} simple roots.
Hence the index set $I$ has a partition $I = I^{\text{re}} \sqcup I^{\text{im}}$,
where  $I^{\text{re}} = \{i \in I \mid a_{ii} = 2\}$,
$I^{\text{im}} = \{i \in I \mid a_{ii} \le 0\}$.
Borcherds algebras
played an important role in Borcherds' proof of the   famous {\it  Moonshine conjecture} \cite{Bor92}.
\vskip 2mm
 T.Bozec  introduced the {\it quantum Borcherds-Bozec algebras} in his geometric investigation of the representation theory of quivers with loops \cite{Bozec2014a,Bozec2014b}. By using simple perverse sheaves on the representation varieties of such quivers, he provided a construction of Lusztig's canonical basis for the positive half of a quantum Borcherds-Bozec algebra.
The Cartan matrices for quantum Borcherds-Bozec algebras are the same for quantum Borcherds algebras, but they have higher degree simple root vectors.
The index set for these simple root vectors is denoted by
$I^{\infty}=(I^{\text{re}}\times \{1\}) \sqcup (I^{\text{im}}\times \Z_{>0})$.
The Weyl-Kac type character formula for the irreducible highest weight modules was proved in \cite{BSV2016}.
In \cite{FKKT21,FKRT21,FHKX}, {\color{black}the first author and his collaborators} have studied the global bases and crystal bases of quantum Borcherds-Bozec algebras.
 Kang, Kim and Tong \cite{KKT} categorified the quantum Borcherds-Bozec algebra  for an arbitrary
 Borcherds-Cartan datum by constructing their associated Khovanov-Lauda-Rouquier algebras and they showed that the cyclotomic Khovanov-Lauda-Rouquier
 algebras provide a categorification of the irreducible highest weight modules.
\vskip 2mm
On the other hand, quantum superalgebras have developed rapidly, particularly with remarkable breakthroughs in structural analysis, representation theory and cross-disciplinary applications.
 In \cite{BKM}, Benkart-Kang-Melville gave the definition of quantum Borcherds superalgebras. In \cite{CHW},  S. Clark, D. Hill, and W. Wang introduced a quantum
 covering and super groups with no isotropic odd simple root and provided substantial research achievements.
In \cite{CHW2}, {\color{black}the first author and his collaborators} showed that the modified
 quantum group and supergroup are isomorphic over the rational function field adjoined
 with $\sqrt{-1}$, by constructing a twistor on the modified covering quantum group. In \cite{BKK}, Benkart, Kang and  Kashiwara developed  the crystal base theory of quantum superalgebra $U_{q}(\mathfrak{g}\mathfrak{l}(m,n))$. In \cite{CHW5}, S. Clark, D. Hill, and W. Wang construct crystal
bases and canonical bases for quantum supergroups of anisotropic type and
for their integrable modules. We proved the
Weyl-Kac type character formula for the irreducible highest
weight modules with dominant integral highest weights of Borcherds-Bozec superalgebras\cite{FHKS}. We define quantum Borcherds-Bozec superalgebras and prove that integrable modules of this quantum superalgebras are completely reducible\cite{FH2}. 
 \vskip 2mm

 The goal of this paper is to construct crystal bases and global bases for $U_q^-(\mathfrak{g})$ and $V(\lambda)$. The main difficulties arise from the countably infinite index set and the imaginary vertex, which necessitate modifications to the Kashiwara operators for both $U_q^-(\mathfrak{g})$ and $V(\lambda)$. Moreover, due to the superalgebra structure, the crystal bases and global bases are not genuine bases in the usual sense.
 
\vskip 2mm
This paper is organized as follows. In Sections 1 and 2, we give the definition of $\U$ and recall that the category $O_{\mathrm{int}}$ is semisimple. In Section 3, we introduce the Boson algebra and construct a bilinear form on $U_q^-(\mathfrak{g})$. In Section 4, we formulate the notions of crystal lattice and crystal basis for integrable modules in the super setting, and also define a polarization on $V(\lambda)$. We establish the tensor product rule for crystal bases and state the main theorems on crystal bases. In Section 5, we prove the existence theorem for crystal bases by using Kashiwara’s grand loop argument, study further properties of the polarization and show that the crystal basis is pseudo-orthonormal but not orthonormal in the usual sense. In Section 6, we define the $\mathcal{A}$-forms $U_{\mathcal{A}}(\mathfrak{g})$ of $\U$ and $V(\lambda)^{\mathcal{A}}$ of $V(\lambda)$. We prove that $U_{\mathcal{A}}(\mathfrak{g})$ and $V(\lambda)^{\mathcal{A}}$ are stable under the Kashiwara operators and establish the existence of global bases.
 \vskip 2mm

{\it Acknowledgements.} Z.Fan was partially supported by the NSF of China grant $12271120$. We would like to express our sincere gratitude to Professor Seok-Jin
 Kang for his helpful discussions. 

\vskip 5mm

\section{Quantum Borcherds-Bozec superalgebras}

Let $I$ be an index set which can be countably infinite.
An integer-valued matrix $A=(a_{ij})_{i,j \in I}$ is called  a {\it
symmetrizable Borcherds-Cartan matrix} if it satisfies the following conditions:
\begin{itemize}
\item[(i)] $a_{ii}=2, 0, -2, -4, ...$,

\item[(ii)] $a_{ij}\le 0$ for $i \neq j$,

\item[(iii)] there exists a diagonal matrix $D=\text{diag} (d_{i} \in \Z_{>0} \mid i \in I)$ such that $DA$ is symmetric.
\end{itemize}

\vskip 2mm

Set $I^{\text{re}}=\{i \in I \mid a_{ii}=2 \}$,
$I^{\text{im}}=\{i \in I \mid a_{ii} \le 0\}$ ,
$I^{\text{iso}}=\{i \in I \mid a_{ii}=0 \}$ and $I^{\infty}:= (I^{\text{re}} \times \{1\}) \sqcup (I^{\text{im}}
\times \Z_{>0})$. For simplicity, we will often write $i$ for $(i,1)$.

\vskip 2mm

{ An {\it odd Borcherds-Cartan datum}} consists of :

\begin{itemize}

\item[(a)] a symmetrizable Borcherds-Cartan matrix $A=(a_{ij})_{i,j \in I}$,

\item[(b)] a free abelian group $P$, the {\it weight lattice},

\item[(c)] $P^{\vee} := \Hom(P, \Z)$, the {\it dual weight lattice},

\item[(d)] $\Pi=\{\alpha_{i} \in P  \mid i \in I \}$, the set of {\it simple roots},

\item[(e)] $\Pi^{\vee}=\{h_i \in P^{\vee} \mid i \in I \}$, the set of {\it simple coroots},

{ \item[(f)] a subset $I_{\bar{1}}$ of $I$, the set of {\it odd indices}},

\end{itemize}

\noindent satisfying the following conditions

\begin{itemize}

\item[(i)] $\langle h_i, \alpha_j \rangle = a_{ij}$ for all $i, j \in I$,

\item[(ii)] $\Pi$ is linearly independent over $\C$,

\item[(iii)] { for all $i \in I_{\bar{1}}$, we have $a_{ij} \in 2 \Z $},

\item[(iv)] for each $i \in I$, there exists an element $\Lambda_{i} \in P$ such that
$$\langle h_j , \Lambda_i \rangle = \delta_{ij} \ \ \ \text{for all} \ i, j \in I,$$

\end{itemize}

\vskip 2mm

{ We denote $I_{\bar{0}} = I \setminus I_{\bar{1}}$, the set of {\it even indices}}.
Given  a symmetrizable Borcherds-Cartan matrix, such
an odd Borcherds-Cartan datum always exists, which
is not necessarily unique. The $\Lambda_i$  $(i \in I)$ are called the {\it fundamental weights}.

An odd Borcherds-Cartan datum is called bar-consistent if it satisfies
$$d_{i}\equiv p(i)\pmod{2}, \ \forall i\in I,$$
where $p$ is a parity function $p:I\rightarrow \{0,1\}$ such that $i\in I_{\overline{p(i)}}$
\vskip 2mm

We define
\begin{equation} \label{eq:P}
P^{+}:=\{\lambda \in P \mid \langle h_i, \lambda \rangle \ge 0 \ \text{for} \ i \in I,
\ \ \langle h_{i}, \lambda \rangle \in 2 \Z \ \text{for} \ i \in I^{\text{re}} \cap I^{\text{odd}}  \}.
\end{equation}
The elements in $P^{+}$ are called the {\it dominant integral weights}.

\vskip 2mm

The free abelian group $Q:= \bigoplus_{i \in I} \Z \, \alpha_i$ is called the {\it root lattice}.
Set $Q_{+}: = \sum_{i \in I} \Z_{\ge 0}\, \alpha_{i}$ and $Q_{-}: = -Q_{+}$.
For $\beta = \sum k_i \alpha_i \in Q\textbf{}_{+}$, we define its {\it height} to be
{ $\text{ht} (\beta):=\sum k_i$}.

\vskip 2mm

Let ${\mathfrak h} := \C \otimes_{\Z} P^{\vee}$ be the {\it Cartan subalgebra}.
We define
a partial order on ${\mathfrak h}^{*}$ by setting $\lambda \ge \mu$
if and only if $\lambda - \mu  \in Q_{+}$ for $\lambda, \mu \in {\mathfrak h}^{*}$.

\vskip 2mm

 We decompose
${\mathfrak h}= \bigoplus_{i \in I} \C \, h_{i}$ .
Since $A$ is symmetrizable, we can define
a symmetric bilinear form $( \ \, , \ \, )$ on ${\mathfrak h}$ by
\begin{equation} \label{eq:h-bilinear}
\begin{aligned}
& (h, h_i) = d_{i}^{-1} \langle h, \alpha_ {i} \rangle \ \ \text{for} \ h \in {\mathfrak h}.
\end{aligned}
\end{equation}

{In particular}, we have
$$(h_i, h_j) = d_{j}^{-1} a_{ij} = d_{i}^{-1} a_{ji} = (h_j, h_i).$$

\vskip 2mm

Moreover, since $\Pi$ is linearly independent,
it is straightforward to verify that $( \ \, , \ \, )$ is non-degenerate on ${\mathfrak h}$.
Hence there is a linear isomorphism $\nu: \mathfrak{h} \rightarrow \mathfrak{h}^{*}$
defined by
\begin{equation} \label{eq:linear_iso}
\nu(h)(h') = (h, h') \ \ \text{for} \  h, h' \in {\mathfrak h},
\end{equation}
{ which implies}
$$\nu(h_i) = d_{i}^{-1} \alpha_{i}, \ \ \nu^{-1}(\alpha_i) = d_i h_i.$$

Then we obtain the induced non-degenerate symmetric bilinear for $(\ \, , \ \, )$
on ${\mathfrak h}^{*}$ defined by
\begin{equation} \label{eq:hstar}
(\lambda, \mu): = (\nu^{-1}(\lambda), \nu^{-1}(\mu)) \ \ \text{for} \ \lambda, \mu \in {\mathfrak h}^{*}.
\end{equation}
In particular, we have
\begin{equation} \label{eq:alpha_i}
(\alpha_i, \lambda) = d_{i} \langle h_{i}, \lambda \rangle \ \ \text{for} \ i \in I.
\end{equation}

\vskip 2mm

For each  $i \in I^{\text{re}}$, we  define the {\it simple reflection}
$r_{i}:{\mathfrak h}^{*} \rightarrow {\mathfrak h}^{*}$ by
\begin{equation} \label{eq:simple_reflection}
r_{i}(\lambda)= \lambda - \langle h_{i}, \lambda \rangle \,  \alpha_{i}
\ \ \text{for} \ \lambda \in {\mathfrak h}^{*}.
\end{equation}

The subgroup $W$ of $GL({\mathfrak h}^{*})$ generated by the simple reflections $r_{i}$ $(i \in I^{\text{re}})$ is called the {\it Weyl group} of the odd Borcherds-Cartan datum given above.
Similar to \cite{Kang2019,FKRT21} , \  $(\ \,  , \ \,  )$ is $W$-invariant.

\vskip 2mm

For any $i\in I$, we set
$$q_{i}=q^{d_{i}},\ q_{(i)}=q^{\frac{(\alpha_{i},\alpha_{i})}{2}}$$
\vskip 3mm

\vskip 3mm
\section{The Quantum Borcherds-Bozec Superalgebras}

 \begin{definition} \label{def:BBsuper}
 {\rm
The {\it quantum Borcherds-Bozec superalgebra} $U_{q}(\g)$ associated with
 an odd Borcherds-Bozec Cartan datum $(A, P, P^{\vee},
 \Pi,\Pi^{\vee}, I_{\bar{1}})$ is associative $\Q(q)$ superalgebra
 generated by the elements $a_{il}$, $b_{il}$ $((i,l)\in I^{\infty})$ and $q^{h}, h\in P^{\vee}$ where the parity is given by $p(a_{il})=p(b_{il})=p(li)$ and $p(q^{h})=0$
 with defining relations
 \begin{equation} \label{eq:BBsuper}
  \begin{aligned}
  & q^{0}=0,\ \ q^{h}q^{h'}=q^{h+h'}, \\
  & q^{h}a_{il}q^{-h}=q^{l\langle h,\alpha_{i}\rangle}a_{il},\ \ q^{h}b_{il}q^{-h}=q^{-l\langle h,\alpha_{i}\rangle}b_{il}, \\
  & \sum_{t+t'=1-ka_{ij}}(-1)^{t'}(-1)^{p(i)(t'p(lj)+\binom{t'}{2})}\begin{bmatrix}1-la_{ij}\\t'\end{bmatrix}_{i}a_{i}^{t}a_{jk}a_{i}^{t'}=0,\ \ i\in I^{re}, i\neq(j,k)\in I^{\infty}, \\
  & \sum_{t+t'=1-ka_{ij}}(-1)^{t'}(-1)^{p(i)(t'p(lj)+\binom{t'}{2})}\begin{bmatrix}1-la_{ij}\\t'\end{bmatrix}_{i}b_{i}^{t}b_{jk}b_{i}^{t'}=0,\ \ i\in I^{re}, i\neq(j,k)\in I^{\infty}, \\
  & a_{il}a_{jk}-(-1)^{p(li)p(kj)}a_{jk}a_{il}=b_{il}b_{jk}-(-1)^{p(li)p(kj)}b_{jk}b_{il}=0,\ \ \text{if}\ a_{ij}=0,  \\
  & a_{il}b_{jk}-(-1)^{p(li)p(kj)}b_{jk}a_{il}=\delta_{ij}\delta_{kl}\frac{K_{i}^{l}-K_{i}^{-l}}{1-(-1)^{p(li)}q_{i}^{2l}},\ K_{i}^{l}=q^{ld_{i}h_{i}}.
  \end{aligned}
 \end{equation}
We shall denote by $U$ (resp. $U^{+},\ U^{0}$ and $U^{-}$) for $U_{q}(\g)$ (resp. $U_{q}^{+}(\g),\ U_{q}^{0}(\g)$ and $U_{q}^{-}(\g)$).

 }
 \end{definition}
The tensor product \(U^{+} \otimes U^{+}\) is made into a \(\Q(q)\)-superalgebra via the multiplication
\[
(x_1 \otimes x_2)(x_3 \otimes x_4) = (-1)^{p(x_2)p(x_3)} q^{(|x_2|,|x_3|)} \, x_1x_3 \otimes x_2x_4.
\]
 We endow $U^{+}$ with a map $\varrho: U^{+}\to U^{+}\otimes U^{+}$ defined by
$$\varrho(a_{il})=a_{il}\otimes 1+1\otimes a_{il}.$$
 
 For each $(i,l)\in I^{\infty}$, let $\varrho_{i,l},\varrho^{\,i,l}: U^{+}\to U^{+}$ be the $\mathbb{Q}$-linear maps defined by
\[
\varrho_{i,l}(1)=\varrho^{\,i,l}(1)=0, \qquad
\varrho_{i,l}(a_{jk})=\varrho^{\,i,l}(a_{jk})=\delta_{ij}\delta_{kl},
\]
and extended by the following twisted Leibniz rules for homogeneous $x,y\in F$:
\begin{align*}
\varrho_{i,l}(xy) &= x\,\varrho_{i,l}(y)
                   + (-1)^{p(li)p(y)} q^{(l\alpha_i,|y|)} \,\varrho_{i,l}(x)\,y, \\[2mm]
\varrho^{\,i,l}(xy) &= \varrho^{\,i,l}(x)\,y
                     + (-1)^{p(li)p(x)} q^{(l\alpha_i,|x|)} \,x\,\varrho^{i,l}(y).
\end{align*}
There is a unique automorphism $\omega:\ U\to U$ such that
$$\omega(a_{il})=(-1)^{p(li)}b_{il},\ \omega(b_{il})=a_{il},\ \omega(q^{h})=q^{-h}$$
for $(i,l)\in I^{\infty},\ h\in P^{\vee}$. We have $\omega(x^{+})=(-1)^{p(x)}x^{-},\ \omega(x^{-})=x^{+}$ for all $x\in F$.

Similarly, there is a unique isomorphism of $\Q(q)$ vector spaces $\sigma: U\to U$ such that
$$\sigma(a_{il})=a_{il},\ \sigma(b_{il})=(-1)^{p(li)}b_{il},\ \sigma(q^{h})=q^{-h}$$
for $(i,l)\in I^{\infty},\ h\in P^{\vee}$ such that $\sigma(uu')=\sigma(u')\sigma(u)$ for $u,u'\in U$. We have
$\sigma(x^{+})=\sigma(x)^{+},\ \sigma(x^{-})=(-1)^{p(x)}\sigma(x)^{-}.$

 \vskip 2mm

Let $U\otimes U$ be a $\Z_{2}$-graded space. We endow $U\otimes U$ with a multiplication given by
$$(a\otimes b)(c\otimes d)=(-1)^{p(b)p(c)}(ac\otimes bd)$$
for homogeneous $a,b,c,d\in U$. This makes $U\otimes U$ into an associative $\Q(q)$-superalgebra with 1.

We have a superalgebra homomorphism $\Delta: U\to U\otimes U$ given by
\begin{align*}
&\Delta(q^{h})=q^{h}\otimes q^{h},\\
&\Delta(a_{il}) =a_{il}\otimes K_{i}^{-l}+1\otimes a_{il},\\
&\Delta(b_{il})=b_{il}\otimes 1+ K_{i}^{l}\otimes b_{il},
\end{align*}
where $i\in I,\ h\in P^{\vee}$. This homomorphism is a coproduct which makes  $U$ a Hopf superalgebra. When $i\in I^{re}$, we have
\begin{align*}
&\Delta(a_{i}^{(n)})=\sum_{t=0}^{n}q_{i}^{t(n-t)}a_{i}^{(t)}\otimes K_{i}^{-t}a_{i}^{(n-t)},\\
&\Delta(b_{i}^{(n)})=\sum_{t=0}^{n} ((-1)^{p(i)}q_{i})^{-t(n-t)}b_{i}^{(n-t)}K_{i}^{t}\otimes b_{i}^{(t)},
\end{align*}
where $a_{i}^{(n)}=\frac{a_i^{n}}{[n]^!_i}$, $b_{i}^{(n)}= \frac{b_i^{n}}{[n]^!_i}$ and ${[n]^!_i}=\prod_{s=1}^{n}\frac{((-1)^{p(i)} q_i)^{s} - q_i^{-s}}{(-1)^{p(i)} q_i - q_i^{-1}}$.
 \vskip 2mm
For  $U$-modules $M$ and $N$, we endow the tensor product $M\otimes N$ with a $U$-module structure induced by $\Delta$.

\begin{proposition}\label{2.2}\cite[Proposition 2.3]{FH2}
For $x\in F$ and $(i,l)\in I^{\infty}$, we have
\begin{itemize}
\item[(a)] $x^{+}b_{il}-(-1)^{p(li)p(x)}b_{il}x^{+}=\frac{\varrho_{i,l}(x)^{+}K_{i}^{l}-(-1)^{p(li)p(x)-p(li)}K_{i}^{-l}\varrho^{i,l}(x)^{+}}{1-(-1)^{p(li)}q_{i}^{2l}},$
\item[(b)]$a_{il}x^{-}-(-1)^{p(li)p(x)}x^{-}a_{il}=\frac{K_{i}^{l}\varrho^{i,l}(x)^{-}-(-1)^{p(li)p(x)-p(li)}\varrho_{i,l}(x)^{-}K_{i}^{-l}}{1-(-1)^{p(li)}q_{i}^{2l}}.$
\end{itemize}

\end{proposition}

\vskip 3mm
Let $M$ be a $\U$-module. $M=M_{\bar{0}}\oplus M_{\bar{1}}$, we say that M has a weight space decomposition if
\[M=\bigoplus_{\mu\in P}M_{\mu}, \ M_{\mu}=\{m\in M|q^{h}m=q^{\langle h,\mu\rangle}m,\ \forall h\in P^{\vee}\},\]
such that $M_{\mu}=M_{\mu,\bar{0}}\oplus M_{\mu,\bar{1}}$, where $M_{\mu,\bar{0}}=M_{\mu}\cap M_{\bar{0}} $ and $M_{\mu,\bar{1}}=M_{\mu}\cap M_{\bar{1}} $.

We denote $wt(M)=\{\mu\in \h^{*}|M_{\mu}\neq 0\}.$ When $\text{dim}\ M_{\mu}<\infty$ for all $\mu\in P$, we define the character of $M$
\[\text{ch}\ M=\sum_{\mu\in P}(\text{dim}\ M_{\mu})e^{\mu},\]
where $e^{\mu}$ denotes the additive basis of $\C[\h^{*}]$.
 \vskip 3mm
\begin{definition}
A $\U$ module V is called a highest weight module with highest weight $\lambda$ if there is a non-zero vector $v_{\lambda}$ in V such that
\begin{itemize}
\item[(a)] $V=\U v_{\lambda},$
\item[(b)] $a_{il}v_{\lambda}=0,\ \forall (i,l)\in I^{\infty},$
\item[(c)]$q^{h}v_{\lambda}=q^{\langle h,\mu\rangle}v_{\lambda},\ \forall h\in P^{\vee}.$
\end{itemize}
The vector $v_{\lambda}$ is called the highest weight vector with highest weight $\lambda$ and $V$ has a weight decomposition $V=\oplus_{\mu \le \lambda}V_{\mu}$.
\end{definition}
\vskip 3mm

 \vskip 3mm
 \begin{proposition}\label{2.4}
 Let $\lambda\in P^{+}$ and let $V(\lambda)=\U v_{\lambda}$ be the irreducible highest weight module with highest weight $\lambda$ and highest weight vector $v_{\lambda}$, $\mu\in \text{wt}(V(\lambda))$, we have
 \begin{itemize}
 \item[(a)] If $i\in I^{re}$,
$b_{i}^{\langle h_{i},\lambda\rangle+1}v_{\lambda}=0$.
 \item[(b)] If $i\in I^{im}$ and $\langle h_{i},\lambda\rangle=0$, then $b_{ik}v_{\lambda}=0$ for all $k>0$.
 \item[(c)] If $i\in I^{im}$, then  $\langle h_{i},\mu\rangle\ge 0$.
 \item[(d)]If $i\in I^{im}$ and $\langle h_{i},\mu\rangle=0$, then $V(\lambda)_{\mu-l\alpha_{i}}=0$ for all $l>0$. In particular, $b_{il}V(\lambda)_{\mu}=0$.
 \item[(e)]If $i\in I^{im}$ and $\langle h_{i},\mu\rangle\le -la_{ii}$, then $a_{il}V(\lambda)_{\mu}=0.$

 \end{itemize}
 \end{proposition}

  \vskip 3mm
 \begin{definition}\label{2.6}
The category $O_{int}$ consists of $\U$-modules $M$ such that
\begin{itemize}
\item[(a)] $M$ has a weight decomposition $M=\oplus_{\mu\in P}M_{\mu}$ and $\text{dim} M_{\mu}<\infty$ for all $\mu\in P$.
\item[(b)]There exist finitely many weights $\lambda_{1},\cdots ,\lambda_{r}\in P$ such that $\text{wt}(M)\subset\cup_{j=1}^{r}(\lambda_{j}-Q_{+})$.
\item[(c)]If $i\in I^{re}$, $b_{i}$ is locally nilpotent on $M$.
\item[(d)] If $i\in I^{im}$, then $\langle h_{i},\mu\rangle\ge0$ for all $\mu\in \text{wt}(M)$.
\item[(e)] If $i\in I^{im}$ and $\langle h_{i},\mu\rangle=0$, then $b_{il}M(\mu)=0$ for all $l>0$.
\item[(f)] If $i\in I^{im}$ and $\langle h_{i},\mu\rangle\le -la_{ii}$, then $a_{il}(M_{\mu})=0$  for all $l>0$.
\end{itemize}
\end{definition}

 The irreducible highest weight $\U$-module $V(\lambda)$ with $\lambda\in P^{+}$ is an object of $O_{int}$ .  Any finite direct sum of modules in $O_{int}$   is also an object of $O_{int}$. By \cite{FH2}, we have the following proposition.

\begin{proposition}
\begin{itemize}
    \item [(a)]Any integrable module in $O_{int}$ is a direct sum of irreducible highest weight modules.
    \item[(b)] Every simple object in the category $O_{int}$ is isomorphic to some irreducible highest weight modulce $V(\lambda)$ with $\lambda \in P^{+}$.
\end{itemize}

\end{proposition}

 \vskip 7mm
\section{The Boson Superalgebra}
We denote the supercommutator $[x,y]:=xy-(-1)^{p(x)p(y)}yx$. Let $(i,l)\in I^{\infty}$ and $y\in U^{-}$, then there exist unique $P,Q\in U^{-}$, such that  $[a_{il},y]=\frac{K_{i}^{l}P-K_{i}^{-l}Q}{1-(-1)^{p(li)}q_{i}^{2l}}$. For $(i,l)\in I^{\infty}$, we define linear operators $e_{il}'$ and $e_{il}''$ in $\text{End}(U^{-})$ by $E_{il}'(y)=Q$ and $E_{il}''(y)=P$.
\vskip 3mm
\begin{lemma}
    For $y\in U^{-}$, we have 
    \begin{align*}
    &e_{il}'(b_{jk}y)=\delta_{ij}\delta_{kl}y+(-1)^{p(li)p(kj)}q_{i}^{-lk\langle h_{i},\alpha_{j}\rangle}b_{jk}e_{il}'(y),\\
        &e_{il}''(b_{jk}y)=\delta_{ij}\delta_{kl}y+(-1)^{p(li)p(kj)}q_{i}^{lk\langle h_{i},\alpha_{j}\rangle}b_{jk}e_{il}''(y).
    \end{align*}
\end{lemma}
\begin{proof}
    By proposition \ref{2.2}, we know that $$e_{il}''(y)=\varrho^{il}(y^{+})^{-}$$ 
    and $$e_{il}'(y)=(-1)^{p(li)p(y)-p(li)}q_{i}^{l\langle h_{i},|y|+l\alpha_{i}\rangle}\varrho_{il}(y^{+})^{-}.$$
    From \cite{FH2}, the lemma follows.
\end{proof}
\vskip 3mm

\begin{corollary}\label{3.2}
    For $x\in U^{-}$, if $e_{il}'(x)=0$ for all $(i,l)\in I^{\infty}$, then $x$ is a constant multiple of $1$
\end{corollary}
\vskip 3mm
\begin{proposition}\label{3.3}
    For $(i,l),(j,k)\in I^{\infty}$, we have $e_{il}'e_{jk}''=(-1)^{p(li)p(kj)}q_{j}^{kl\langle h_{j},\alpha_{i}\rangle}e_{jk}''e_{il}'$
\end{proposition}
\begin{proof}
    Let $y=b_{i_{1}l_{1}}\cdots b_{i_{r}l_{r}}$. Use induction on $r$. If $r=1$, then $$e_{il}'e_{jk}''(y)=(-1)^{p(li)p(kj)}q_{j}^{kl\langle h_{j},\alpha_{i}\rangle}e_{jk}''e_{il}'(y)=0.$$
    If $r\ge 1$, then $y=b_{i_{1}l_{1}}y'$ with $y'\in U^{-}$. We have
    \begin{align*}
e_{il}'e_{jk}''(b_{i_{1}l_{1}}y')&=e_{il}'(\delta_{i_{1}j}\delta_{l_{1}k}y'+(-1)^{p(kj)p(l_{1}i_{1})}q_{j}^{l_{1}ka_{ji_{1}}}b_{i_{1}l_{1}}e_{jk}''(y'))\\
        &=\delta_{i_{1}j}\delta_{l_{1}k}e_{il}'(y')+(-1)^{p(kj)p(l_{1}i_{1})}q_{j}^{l_{1}ka_{ji_{1}}}e_{il}'(b_{i_{1}l_{1}}e_{jk}''(y'))\\
        &=\delta_{i_{1}j}\delta_{l_{1}k}e_{il}'(y')+(-1)^{p(kj)p(l_{1}i_{1})}q_{j}^{l_{1}ka_{ji_{1}}}\delta_{ii_{1}}\delta_{ll_{1}}e_{jk}''(y')\\
      &+(-1)^{p(kj)p(l_{1}i_{1})+p(li)p(l_{1}i_{1})}q_{j}^{l_{1}ka_{ji_{1}}}q_{i}^{-ll_{1}a_{ii_{1}}}b_{i_{1}l_{1}}e_{il}'e_{jk}''(y'),
    \end{align*}
    and 
    \begin{align*}
       e_{jk}''e_{il}' (y')&=e_{jk}''(\delta_{ii_{1}}\delta_{ll_{1}}y'+(-1)^{p(l_{1}i_{1})p(li)}q_{i}^{-ll_{1}a_{ii_{1}}}b_{i_{1}l_{1}}e_{il}'(y'))\\
       &=\delta_{ii_{1}}\delta_{ll_{1}}e_{jk}''(y')+(-1)^{p(l_{1}i_{1})p(li)}q_{i}^{-ll_{1}a_{ii_{1}}}\delta_{ji_{1}}\delta_{kl_{1}}e_{il}'(y')\\
       &+(-1)^{p(l_{1}i_{1})p(li)+p(l_{1}i_{1})p(kj)}q_{j}^{l_{1}ka_{ji_{1}}}q_{i}^{-ll_{1}a_{ii_{1}}}b_{i_{1}l_{1}}e_{jk}''e_{il}'(y').
    \end{align*}
    These equations imply that $$e_{il}'e_{jk}''(b_{i_{1}l_{1}}y')=(-1)^{p(li)p(kj)}q_{j}^{kl\langle h_{j},\alpha_{i}\rangle}e_{jk}''e_{il}'(b_{i_{1}l_{1}}y').$$
    Therefore, we have 
   $e_{il}'e_{jk}''(y)=(-1)^{p(li)p(kj)}q_{j}^{kl\langle h_{j},\alpha_{i}\rangle}e_{jk}''e_{il}'(y)$.
\end{proof}
\vskip 3mm
\begin{definition}
     Let $\mathscr{B}_{q}(\g)$ be the algebra over $\Q(q)$ generated by $a_{il}', b_{il},\  (i,l)\in I^{\infty}$ with defining relations
\begin{align*}
&a_{il}'b_{jk}=\delta_{ij}\delta_{lk}+(-1)^{p(li)p(kj)}q_{i}^{-kla_{ij}}b_{jk}a_{il}',\ \forall (i,l),(j,k)\in I^{\infty},\\
& \sum_{t+t'=1-la_{ij}}(-1)^{t'}(-1)^{p(i)(t'p(lj)+\binom{t'}{2})}\begin{bmatrix}1-la_{ij}\\t'\end{bmatrix}_{i}a_{i}'^{t}a_{jl}'a_{i}'^{t'}=0,\ \ i\in I^{re}, i\neq(j,k)\in I^{\infty}, \\
  & \sum_{t+t'=1-la_{ij}}(-1)^{t'}(-1)^{p(i)(t'p(lj)+\binom{t'}{2})}\begin{bmatrix}1-la_{ij}\\t'\end{bmatrix}_{i}b_{i}^{t}b_{jl}b_{i}^{t'}=0,\ \ i\in I^{re}, i\neq(j,k)\in I^{\infty}, \\
  & a_{il}'a_{jk}'-(-1)^{p(li)p(kj)}a_{jk}'a_{il}'=b_{il}b_{jk}-(-1)^{p(li)p(kj)}b_{jk}b_{il}=0,\ \ \text{if}\ a_{ij}=0. \\
\end{align*}
\end{definition}
\vskip 3mm
\begin{proposition}
    $U^{-}$ is a $\B$-module, where $b_{il}$ acts as
 multiplication by $b_{il}$ and $a_{il}'$ acts by  $e_{il}'$ for all $(i,l)\in I^{\infty}$.
\end{proposition}
\begin{proof}
    Set $N_{ij}=1-la_{ij}$ and
  $$S_{i,(j,l)}=\sum_{t=0}^{N_{ij}}(-1)^{t}(-1)^{p(i)(tp(lj)+\binom{t}{2})}\begin{bmatrix}N_{ij}\\t\end{bmatrix}_{i}a_{i}'^{N_{ij}-t}a_{jl}'a_{i}'^{t}.$$
  $$S_{i,(j,l)}'=\sum_{t=0}^{N_{ij}}(-1)^{t}(-1)^{p(i)(tp(lj)+\binom{t}{2})}\begin{bmatrix}N_{ij}\\t\end{bmatrix}_{i}b_{i}^{N_{ij}-t}b_{jl}b_{i}^{t}.$$
  It suffices to show that $S_{i,(j,l)}$ and $S_{i,(j,l)}'$ act as zero on any $y\in U^{-}$.
  By the definition of the action and serre relations in $U$, $S_{i,(j,l)}'y=0$.

  Regarding $S_{i,(j,l)}$, assume that $y=b\cdot 1$ where $b\in \B$ is the monomial in the generators $b_{il}$ for $(i,l)\in I^{\infty}$. It is easy to show that
  $S_{i,(j,l)}y=cbS_{i,(j,l)}\cdot1$ for some $c\in \Q(q)$. Since $E_{il}'(1)=0$ for all $(i,l)\in I^{\infty}$, $S_{i,(j,l)}y=0$.
  
\end{proof}
\vskip 3mm
\begin{corollary}
    As $\B$-modules, $U^{-}\cong \B\slash \sum_{(i,l)\in I^{\infty}}\B a_{il}'$
\end{corollary}
\begin{proof}
    Since $1$ is annihilated by $a_{il}'$, we have a surjective morphism
    $$\B\slash \sum_{(i,l)\in I^{\infty}}\B a_{il}'\to U^{-}.$$
    Let $B^{-}$ be the subalgebra of $\B$ generated by $b_{il}$, then we have 
    $$B^{-}\xrightarrow{\phi} \B\slash \sum_{(i,l)\in I^{\infty}}\B a_{il}' \xrightarrow{\psi} U^{-},$$
where $\phi$ and $\psi$ are surjective. It is trivial to show that $\psi\circ \phi$ is bijective. Hence $\psi$ is an isomorphism.
\end{proof}
\vskip 3mm
\begin{proposition} \label{3.7}
    There is a unique bilinear form $(\cdot,\cdot)$ on $U^{-}$ satisfying
    $$(1,1)=1,\ (b_{il}u,v)=(u,a_{il}'v),\ \forall u,v\in U^{-}, (i,l)\in I^{\infty}.$$
\end{proposition}
\begin{proof}
 
   There is a unique linear map $f:\B\to\B$ such that $$f(a_{il}')=b_{il},f(b_{il})=a_{il}', f(xy)=f(y)f(x).$$
   $(U^{-})^{*}$ is a $\B$-module via $(x\cdot\phi)(y)=\phi(f(x)\cdot y)$, for $x\in B,\ y\in U^{-}$ and  $\phi\in (U^{-})^{*} $.

   Let $\phi_{0}\in (U^{-})^{*}$ be defined by $\phi_{0}(1)=1$ and $\phi_{0}(\sum_{(i,l)}b_{il}U^{-})=0$. Since $a_{il}'\phi_{0}(x)=\phi_{0}(b_{il}x)=0$ for all $x\in U^{-} $, $a_{il}'\phi_{0}=0$ for any $(i,l)\in I^{\infty}$, we have a homomorphism
   $$\psi:U^{-}\cong\B\slash \sum_{(i,l)\in I^{\infty}}\B a_{il}'\to (U^{-})^{*},$$
   which maps $1$ to $\phi_{0}$. We define a bilinear form $(\cdot,\cdot)$ on $U^{-}$
   $$(u,v)=\psi(u)(v).$$
   Then we have $$(1,1)=1,$$
   $$(b_{il}u,v)=(u,a_{il}'v)\ \ \text{and} \ \ (a_{il}'u,v)=(u,b_{il}v).$$
   It is trivial to show that the bilinear form is unique. Let $(\cdot ,\cdot )'$ be a bilinear form defined by $(u,v)':=(v,u)$. Since $(\cdot ,\cdot )'$ satisfies the same properties, by the uniqueness of such a bilinear form, We have $(u,v)=(v,u)$. 
\end{proof}   

\vskip 3mm
\begin{proposition}
    The bilinear form $(\cdot,\cdot)$ on $U^{-}$ is nondegenerate.
\end{proposition}
\begin{proof}
   When $\text{ht}(\nu)=0$, it is trivial. Assume that $\text{ht}(\nu)>0$, let $y$ be a nonzero element in  $ U_{\nu}^{-}$ such that $(y, U_{\nu}^{-})=0$, then $(a_{il}'y,U_{\nu+l\alpha_{i}}^{-})=(y,b_{il}U_{\nu+l\alpha_{i}}^{-})=0$ for all $(i,l)\in I^{\infty}$. By induction hypothesis, we have $a_{il}'y=0$ for all $(i,l)\in I^{\infty}$. From Corollary \ref{3.2}, $\nu=0$ and we get a contradiction.
\end{proof}
\vskip 3mm




\vskip 3mm
\section{Crystal Bases}
Let $\mathbf{c}=(c_{1},\cdots,c_{r})\in \Z_{\ge0}^{r}$ be a sequence of non-negative integers. We define $|\mathbf{c}|:=\sum_{t=1}^{r}c_{t}$. We denote by $C_{il}$ the set of compositions of $l$. We define 
$$a_{i,\mathbf{c}}=a_{i,c_{1}}\cdots a_{i,c_{r}},\ \ b_{i,\mathbf{c}}=b_{i,c_{1}}\cdots b_{i,c_{r}}$$.

\subsection{Crystal bases for $V(\lambda)$}

\begin{proposition}
     Assume $(i,\lambda)\ge 0$ for some imaginary vertex $i$. We have a decomposition $V(\lambda)=\bigoplus_{\mathbf{c}\in C_{i}}b_{i,\mathbf{c}}\mathcal{K}_{i}$, where $\mathcal{K}_{i}=\cap_{l>0}\text{ker} \ a_{il}$.
\end{proposition}
\begin{proof}
Consider $v\in V(\lambda)$, assume first that $v$ has the form $v=ub_{i,\mathbf{c}}z$ for some $c\in \mathcal{C}_{i}$, $u\in U^{-}$ satisfying $[a_{il},u]=0$ for every $l$, and  $z\in \mathcal{K}_{i}$. Proceed by induction on $|\mathbf{c}|$. If $( \alpha_{i},|u|)=0$, since $i\in I^{im}$, we have $\text{supp}\ |u|\subset\{j\in I| ( \alpha_{i},\alpha_{j})=0\}$. Hence $[u,b_{il}]=0$ for any $l$ and we get the result by induction. If $( \alpha_{i},|u|)>0$, we set 
\begin{align*}
    l&=c_{1},\\
    [u,b_{il}]^{\circ}&=ub_{il}-(-1)^{p(li)p(u)}Rb_{il}u\ \  \text{for\ some}\  R\in \Q(q),\\
    z'&=b_{i,\mathbf{c}\slash c_{1}}z\in V(\lambda)_{\mu}.
\end{align*}
$\forall k>0$, $u\in U^{-}$, since $a_{ik}u-(-1)^{p(u)p(ki)}ua_{ik}=0$, then $\varrho_{ik}(u^{+})=\varrho^{ik}(u^{+})=0$ and 
\begin{align*}
\varrho^{ik}(u^{+}a_{il})&=\delta_{kl}(-1)^{p(ki)p(u)}q^{(k\alpha_{i},-|u|)}u^{+},\ \ \varrho_{ik}(u^{+}a_{il})=\delta_{kl}u^{+},\\
\varrho_{ik}(a_{il}u^{+})&=\delta_{kl}(-1)^{p(ki)p(u)}q^{(k\alpha_{i},-|u|)}u^{+},\ \ \varrho^{ik}(a_{il}u^{+})=\delta_{kl}u^{+}.
\end{align*}
If $$R=\frac{q^{(k\alpha_{i},\mu)}-q^{-(k\alpha_{i},\mu)}}{q^{(k\alpha_{i},\mu+|u|)}-q^{-(k\alpha_{i},\mu+|u|)}},$$
then
\begin{align*}
    [a_{ik},[u,b_{il}]^{\circ}]z'&=\{a_{ik}[u,b_{il}]^{\circ}-(-1)^{p(ki)(p(u)+p(li))}[u,b_{il}]^{\circ}a_{ik}\}z'\\
    &=\delta_{kl}(-1)^{p(li)p(u)}\frac{q^{(k\alpha_{i},\mu)}-q^{-(k\alpha_{i},\mu)}-R(q^{(k\alpha_{i},\mu+|u|)}-q^{-(k\alpha_{i},\mu+|u|)})}{1-(-1)^{p(li)}q_{i}^{2l}}\\
    &=0,
\end{align*}
which is possible since 
$$(\alpha_{i},|u|+\mu)=(\alpha_{i},|u|)+(\alpha_{i},\lambda)+(\alpha_{i},\mu-\lambda)>(\alpha_{i},|u|)+(\alpha_{i},\mu-\lambda)\ge0.$$
We have used that since $i\in I^{im}$ and $\mu-\lambda\in -\N I$, we have 
$$(\alpha_{i},\mu-\lambda)\ge 0.$$
Hence, the equality
$$v=[u,b_{il}]^{\circ}b_{i,\mathbf{c}\slash c_{1}}z+(-1)^{p(li)p(u)}Rb_{il}ub_{i,\mathbf{c}\slash c_{1}}z
$$
along with the induction hypothesis allow us to conclude since $|\mathbf{c}\slash c_{1}|<|c|$ and $\bigoplus_{\mathbf{c}\in C_{i}}b_{i,\mathbf{c}}\mathcal{K}_{i}$ is stable under the multiplication by $b_{il}$.

We will prove the existence of the decomposition for a general $v\in V(\lambda)_{\mu}$ , using induction on $|\lambda-\mu|$. If $v\neq v_{\lambda}$, we can write $v=\sum_{\iota\in I^{\infty}}b_{\iota}v_{\iota}$ for some finitely many nonzero $v_{\iota}\in V(\lambda)$ by induction hypothesis, we have $v=\sum_{\iota\in I^{\infty},\mathbf{c}\in \mathcal{C}_{i}}b_{\iota}b_{i,\mathbf{c}}z_{\iota,\mathbf{c}}$ for some finitely many nonzero $z_{\iota,\mathbf{c}}\in \mathcal{K}_{i}$. And 
$$v=\sum_{l>0,\mathbf{c}\in \mathcal{C}_{i}}b_{i,(l,\mathbf{c})}z_{(i,l),\mathbf{c}}+\sum_{\iota\notin\{\{i\}\times \N_{>0}\},\mathbf{c}\in \mathcal{C}_{i}}b_{\iota}b_{i,\mathbf{c}}z_{\iota,\mathbf{c}}.$$
we get the result Since  $b_{\iota}b_{i,\mathbf{c}}z_{\iota,\mathbf{c}}$ is of the form   $ub_{i,\mathbf{c}}z$.

Next we will prove the uniqueness of the decomposition. Consider the  minimal nontrivial relation of dependence:
$$0=\sum_{\mathbf{c}\in \mathcal{C_{i}}}b_{i,\mathcal{c}}z_{\mathbf{c}},$$
where $z_{\mathbf{c}}\in V(\lambda)_{\mu+|\mathbf{c}|i}\cap \mathcal{K}_{i}.$  When $i\notin I^{iso}$. Consider $r$ maximal such that there exists $\mathbf{c}=(c_{1},\cdots,c_{r})$ such that $z_{\mathbf{c}}\neq 0$. For any $k\in \{1,\cdots,r\}$
$$c_{<k}=(c_{1},\cdots,c_{k-1})$$
$$c_{>k}=(c_{k+1},\cdots,c_{r})$$
$c_{<1}=c_{>r}=\emptyset$. Then, if $l>0$, we get 
\begin{align*}
    [a_{il},b_{i,\mathbf{c}}]&=\frac{K_{i}^{l}\varrho^{il}(a_{i,\mathbf{c}})^{-}-(-1)^{p(li)p(|\mathbf{c}|i)-p(li)}\varrho_{il}(a_{i,\mathbf{c}})^{-}K_{i}^{-l}}{1-(-1)^{p(li)}q_{i}^{2l}}\\
    &=\frac{\sum_{k:c_{k}=l}(-1)^{p(c_{<k}i)p(li)}b_{i,\mathbf{c}\backslash c_{k}}(q^{-(l\alpha_{i},c_{>k}\alpha_{i})})K_{i}^{l}-q^{(l\alpha_{i},c_{>k}\alpha_{i})}K_{i}^{-l}}{1-(-1)^{p(li)}q_{i}^{2l}}.
\end{align*}
Since $z_{\mathbf{c}}\in \mathcal{K}_{i}$, we have 
\begin{align*}
    a_{il}b_{i,\mathbf{c}}z_{\mathbf{c}}&=[a_{il},b_{i,\mathbf{c}}]z_{\mathbf{c}}+(-1)^{p(li)p(|\mathbf{c}i|)}b_{i,\mathbf{c}}a_{il}z_{\mathbf{c}}\\
    &=\frac{\sum_{k:c_{k}=l}(-1)^{p(c_{<k}i)p(li)}b_{i,\mathbf{c}\backslash c_{k}}(q^{(l\alpha_{i},\mu+c_{\leq k}\alpha_i)}-q^{-(l\alpha_{i},\mu+c_{\leq k}\alpha_{i})})}{1-(-1)^{p(li)}q_{i}^{2l}}z_{\mathbf{c}}.
\end{align*}
For any $\mathbf{c}'\in \mathfrak{S}_{r}$, since $r$ is maximal, we have 
\begin{align*}
    0&=a_{i,\mathbf{c}'}\sum_{\mathbf{c}''\in \mathcal{C}_{i}}b_{i,\mathbf{c''}}z_{\mathbf{c}''}\\
    &=a_{i,\mathbf{c'}}\sum_{\mathbf{c}''\in \mathfrak{S}_{r}\mathbf{c}}b_{i,\mathbf{c}''}z_{\mathbf{c}''}\\
    &=\tau_{i,\mathbf{c}}\sum_{\mathbf{c}''\in \mathfrak{S}_{r}\mathbf{c}}P_{\mathbf{c}',\mathbf{c}''}(q)z_{\mathbf{c}''},
\end{align*}
where $P_{\mathbf{c}',\mathbf{c}''}(q)\in \Z[q,q^{-1}]$. Since $i\in I^{im} ,\ i\notin I^{iso}$ and $\lambda-(\mu+|\mathbf{c}|\alpha_{i})\in \N I$,
$$(\alpha_{i},\mu+|\mathbf{c}|\alpha_{i})=(\alpha_{i},\lambda)+(\alpha_{i},\mu+|\mathbf{c}|\alpha_{i})\ge 0.$$
Then
$$\text{max}_{\mathbf{c}''\in \mathfrak{S}_{r}\mathbf{c}}\{\text{deg}(P_{\mathbf{c}',\mathbf{c}''})\}=\text{deg}P_{\mathbf{c}',\overline{\mathbf{c}'}}=\sum_{1\le k\le r}c_{k}(\alpha_{i},\mu+c_{k}\alpha_{i})=t,$$
We can derive that $t\neq 0$ from that $V(\lambda)$ is simple. Hence $\text{det}(P_{\mathbf{c}',\mathbf{c}''})_{\mathbf{c}',\mathbf{c}''\in \mathfrak{S}_{r}\mathbf{c}}=|\mathfrak{S}_{r}\mathbf{c}|t \neq 0$.

Finally, we consider the situation $i\in I^{iso}$. Write a relation of dependence of minimal degree
$$0=\sum_{\nu\in \mathcal{C}_{i}}b_{i,\nu}z_{\nu}$$
where $z_{\nu}\in V(\lambda)_{\mu+|\nu|\alpha_{i}}\cap \mathcal{K}_{i}$. Without loss of generality, We can assume that $(l\alpha_{i},\mu)\neq 0$, We have 
$$ a_{il}\sum_{\nu\in \mathcal{C}_{i}}b_{i,\nu}z_{\nu}=\tau_{i,l}(q^{-(l\alpha_{i},\mu)\neq 0}.-q^{(l\alpha_{i},\mu)})\sum_{\nu'\in \mathcal{C}_{i}}b_{i,\nu'}\{m_{l}(\nu')+1\}z_{(\nu',l)}$$
which contradicts the minimality of the relation. 
\end{proof}

For $i\in I^{re}$, by  \cite{CHW2}, the vector $v\in V(\lambda)$ can be written uniquely as 
$$v=\sum_{k\ge 0}b_{i}^{(k)}v_{k}$$
such that 
\begin{itemize}
    \item [(1)] $a_{i}v_{k}=0$ for all $k\ge 0$,
    \item[(2)] $u_{k}\in V(\lambda)_{\mu+k\alpha_{i}}$,
    \item[(3)] $u_{k}=0$ if $\langle h_{i},\mu+k\alpha_{i}\rangle=0$. 
\end{itemize}

\begin{definition} \label{def:Kashiwara operator} \hfill

\vskip 3mm

(a) For $i \in I^{\text{re}}$, we define
\begin{equation} \label{eq:real Kas}
\begin{aligned}
& \widetilde{e}_{i} u = \sum_{k \ge 1} {\mathtt b}_{i}^{(k-1)} u_{k}, \\
& \widetilde{f}_{i}  u = \sum_{k \ge 0} {\mathtt b}_{i}^{(k+1)} u_{k}.
\end{aligned}
\end{equation}

\vskip 2mm

(b) For $i \in I^{\text{im}}\setminus I^{\text{iso}}$ and $l>0$, we define
\begin{equation} \label{eq: nonisoKas}
\begin{aligned}
& \widetilde{e}_{il} u = \sum_{\mathbf{c} \in {\mathcal C}_{i}:  c_1 = l} \,
{\mathtt b}_{i, \mathbf{c} \setminus c_1} u_{\mathbf{c}},\\
& \widetilde{f}_{il}  u = \sum_{\mathbf{c} \in {\mathcal C}_{i}}
{\mathtt b}_{i,(l, \mathbf{c})} u_{\mathbf{c}}.
\end{aligned}
\end{equation}

\vskip 2mm

(c) For $i \in I^{\text{iso}}$ and $l>0$, we define
\begin{equation} \label{eq:isoKas}
\begin{aligned}
& \widetilde{e}_{il} u = \sum_{\mathbf{c} \in {\mathcal C}_{i}} \,  (-1)^{p(li)p(|c_{<r}|i)}\mathbf{c}_{l} \,
{\mathtt b}_{i, \mathbf{c} \setminus \{l\}} u_{\mathbf{c}},\\
& \widetilde{f}_{il}  u = \sum_{\mathbf{c} \in {\mathcal C}_{i}} \, \frac{1} {\mathbf{c}_{l} + 1} \,
{\mathtt b}_{i,(l,\mathbf{c})} u_{\mathbf{c}},
\end{aligned}
\end{equation}
where ${\mathbf c}_{l}$ denotes the number of $l$ in $\mathbf{c}$, $r=\min\{t|c_{t}=l\}$.

\end{definition}

\vskip 2mm

Let $\A_{0} = \{f \in \Q(q) \mid \text{$f$ is regular at $q=0$} \}$.
 We have an isomorphism
 \begin{equation}\label{tg}
     \A_{0} / q \A_{0} \cong \Q, \quad f + q \A_{0} \longmapsto f(0).
 \end{equation}
\vskip 2mm

\begin{definition}

Let $M$ be a $U_{q}(\g)$-module in the category ${\mathcal O}_{\text{int}}$
and let $L$ be a free $\A_{0}$-submodule of $M$.
The submodule $L$ is called a {\it crystal lattice} of $M$ if the following conditions hold.

\begin{enumerate}

\vskip 2mm

\item[(a)] $\Q \otimes_{\A_{0}} L \cong M$,

\vskip 2mm

\item[(b)] $L =\oplus_{\mu \in P} L_{\mu}$, where $L_{\mu} = L \cap M_{\mu}$,

\vskip 2mm

\item[(c)] $\widetilde{e}_{il} L \subset L$, $\widetilde{f}_{il} L \subset L$
for $(i, l) \in I^{\infty}$.

\end{enumerate}

\end{definition}

\vskip 3mm

Since the operators $\widetilde e_{il}$, $\widetilde f_{il}$ preserve $L$,
they induce the operators
$$\widetilde{e}_{il}, \, \widetilde{f}_{il}: L/qL \longrightarrow L/qL.$$

\vskip 3mm

\begin{definition} \label{def:crystal basis Blambda} \hfill

Let $M$ be a $U_{q}(\g)$-module in the category ${\mathcal O}_{\text{int}}$.
A {\it crystal basis} of $M$ is a pair $(L, B)$ such that
\begin{enumerate}

\vskip 2mm

\item[(a)] $L$ is a crystal lattice of $M$,

\vskip 2mm

\item[(b)] $B$ is a pseudo-basis of $L/qL$, i.e. there exists a $\Q$-basis $B^{\circ}$ of $L/qL$ such that $B^{\circ}\subset B\subset B^{\circ}\cup (-B^{\circ}) $

\vskip 2mm

\item[(c)] $B=\sqcup_{\mu \in P}\ B_\mu$, where
$B_\mu=B\cap{(L/qL)}_\mu$,

\vskip 2mm

\item[(d)] $\widetilde{e}_{il}B\subset B\cup\{0\}$, $\widetilde{f}_{il} B\subset B\cup\{0\}$
for $(i, l) \in I^{\infty}$,

\vskip 2mm

\item[(e)] for any $b,b'\in B$ and $(i,l)\in I^{\infty}$, we have
$\widetilde{f}_{il}b=b'$ if and only if $b=\widetilde{e}_{il}b'$.

\end{enumerate}
\end{definition}

Let $V(\lambda)=U_q(\mathfrak g)v_\lambda$ be the irreducible highest weight $U_q(\g)$-module
with highest weight $\lambda \in P^+$.
Let $L(\lambda)$ be the free $\A_0$-submodule of $V(\lambda)$.
spanned by $\widetilde{f}_{i_1l_1}\cdots\widetilde{f}_{i_rl_r} v_\lambda$ $(r\geq 0, (i_k,l_k)\in I^{\infty})$
and let
$
\mathcal{B}(\lambda):=\{\widetilde{f}_{i_1l_1}\cdots\widetilde{f}_{i_rl_r}v_\lambda+qL(\lambda)\}\setminus\{0\}
$, and $B(\lambda)=\mathcal{B}(\lambda)\cup (-\mathcal{B}(\lambda))$. Clearly, $B(\lambda)=\sqcup_{\alpha\in Q^+} B(\lambda)_{\lambda-\alpha}$ and $B(\lambda)_{\lambda-\alpha}$ is a finite set. There exists a minimal subset $B(\lambda)_{\lambda-\alpha}^{\circ}\subset \mathcal{B}(\lambda)_{\lambda-\alpha}$ such that $B(\lambda)_{\lambda-\alpha}=B(\lambda)_{\lambda-\alpha}^{\circ}\cup(-B(\lambda)_{\lambda-\alpha}^{\circ})$. Define $B(\lambda)^{\circ}=\sqcup_{\alpha\in Q^+}B(\lambda)_{\lambda-\alpha}^{\circ}$.
\vskip 3mm
\begin{proposition}\label{4.5}
Let $\lambda\in P^{+}$, there is a unique symmetric bilinear form $(\cdot,\cdot)_{K}$ on $V(\lambda)$ such that 
\begin{align*}
&(v_{\lambda},v_{\lambda})_{K}=1,\ \ \  (q^{h}u,v)_{K}=(u,q^{h}v)_{K},\\
&(b_{il}u,v)_{K}=-(u,K_{i}^{l}a_{il}v)_{K},\\
&(a_{il}u,v)_{K}=-(u,K_{i}^{-l}b_{il}v)_{K}.
\end{align*}
where $u,v\in V(\lambda)$ and $h\in P^{\vee}$.

\end{proposition}
\begin{proposition}
Assume that $\lambda,\mu\in P^{+}$. Then the symmetric bilinear form on $V(\lambda)\otimes V(\mu)$ given by $(m_{1}\otimes n_{1},m_{2}\otimes n_{2})=(m_{1},m_{2})(n_{1},n_{2})$ satisfies 
$$(\Delta(u)(m_{1}\otimes n_{1}),m_{2}\otimes n_{2})=(m_{1}\otimes n_{1},\Delta(\tau(u))(m_{2}\otimes n_{2})),$$
where $\tau: U\to U $ is a anti-automorphism defined by 
\begin{align*}
    &\tau(a_{il})=-K_{i}^{l}b_{il},\ \tau(b_{il})=-K_{i}^{l}a_{il},\ \tau(q^{h})=q^{h},\ \ (i,l)\in I^{\infty},\ h\in P^{\vee}.
\end{align*}
\end{proposition}
\begin{proof}
Let $m_{1},m_{2}\in V(\lambda)$ and $n_{1},n_2\in V(\mu)$. We have 
\begin{align*}
   (\Delta(a_{il})(m_{1}\otimes n_{1}),m_{2}\otimes n_{2})&=(a_{il}m_{1}\otimes K_{i}^{-l}n_{1},m_{2}\otimes n_{2})+(-1)^{p(li)p(m_{1})}(m_{1}\otimes a_{il}n_{1},m_{2}\otimes n_{2})\\
   &=(a_{il}m_{1},m_{2})(K_{i}^{-l}n_{1},n_{2})+(-1)^{p(li)p(m_{1})}(m_{1},m_{2})(a_{il}n_{1},n_{2})
\end{align*}
and 
\begin{align*}
    (m_{1}\otimes n_{1},\Delta(\tau(a_{il}))(m_{2}\otimes n_{2}))&=-(m_{1}\otimes n_{1},-(K_{i}^{-l}b_{il}\otimes K_{i}^{-l}+1\otimes K_{-l}b_{il})m_{2}\otimes n_{2})\\
    &=(-1)^{p(li)p(m_{2})}(m_{1},m_{2})(n_{1},-K_{i}^{-l}b_{il}n_{2})\\
     &+(m_{1},-K_{i}^{l}b_{il}m_{2})(n_{1},K_{i}^{-l}n_{2}).
\end{align*}
By Proposition \ref{4.5}, we have
$$(\Delta(a_{il})(m_{1}\otimes n_{1}),m_{2}\otimes n_{2})=(m_{1}\otimes n_{1},\Delta(\tau(a_{il}))(m_{2}\otimes n_{2}))$$.
\end{proof}

\vskip 3mm
For $\lambda,\mu\in P^{+}$, we define $U$-linear homomorphisms
\begin{align*}
    &\Phi(\lambda,\mu): V(\lambda+\mu)\to V(\lambda)\otimes V(\mu) ,\ \ \ v_{\lambda+\mu}\to v_{\lambda}\otimes v_{\mu},\\
    &\Psi(\lambda,\mu):V(\lambda)\otimes V(\mu)\to V(\lambda+\mu),\ \ \ v_{\lambda}\otimes v_{\mu}\to v_{\lambda}.
\end{align*}
Clearly, $\Psi\circ\Phi=\text{id}_{V_{\lambda+\mu}}$
For $\lambda\in P^{+}$, let $\pi_{\lambda}: U^{-}\to V(\lambda)$ be the $U^{-}$-module homomorphism given by $P\to Pv_{\lambda}$.
\vskip 3mm
\begin{corollary}\label{4.7}
    Let $\lambda,\mu\in P^{+}$, then we have 
    $$(\Psi(\lambda,\mu)(w),v)=(w,\Phi((\lambda,\mu))(v)),$$
    for $v\in V((\lambda+\mu))$  and  $  w\in V(\lambda)\otimes V(\mu)$
\end{corollary}
\vskip 3mm
\begin{theorem}\label{thm:crystal basis of V}

{\rm
The pair $(L(\lambda),B(\lambda))$ is a crystal basis of $V(\lambda)$.
}
\end{theorem}

\begin{definition}

	An {\it abstract crystal}  is a set $B$ together with the maps ${\rm wt}\colon  B \rightarrow P$,
	$\varphi_i,\varepsilon_i\colon  B\rightarrow \Z\cup \{-\infty\}$ $(i\in I)$ and $\widetilde{e}_{il},\widetilde{f}_{il}\colon  B\rightarrow   B\cup \{0\}$ $((i,l)\in I^\infty)$ satisfying the following conditions:
	\vskip 2mm
	\begin{enumerate}
		\item[{\rm (a)}] $\text{wt}(\widetilde{f}_{il}b)=\text{wt}(b)-l\alpha_i$ if $\widetilde{f}_{il}b\neq 0$,\quad $\text{wt}(\widetilde{e}_{il}b)=\text{wt}(b)+l\alpha_i$ if $\widetilde{e}_{il}b\neq 0$.
		
		\vskip 2mm
		
		\item[{\rm (b)}] $\varphi_i(b)= \langle h_i, \text{wt}(b) \rangle+\varepsilon_i(b)$ for $i\in I$ and $b\in  B$.
		
		\vskip 2mm
		
		\item[{\rm (c)}] $\widetilde{f}_{il}b=b'$ if and only if $b=\widetilde{e}_{il}b'$ for $(i,l)\in I^\infty$ and $b,b'\in  B$.
		
		\vskip 2mm
		
		\item[{\rm (d)}]  For any $i\in I^{\text{re}}$ and $b\in  B$, we have
		\vskip 2mm
		\begin{itemize}
			\item [(1)] $\varepsilon_i(\widetilde{f}_ib)=\varepsilon_i(b)+1$, $\varphi_i(\widetilde{f}_ib)=\varphi_i(b)-1$ if $\widetilde{f}_ib\neq 0$,
			\vskip 2mm
			\item [(2)] $\varepsilon_i(\widetilde{e}_ib)=\varepsilon_i(b)-1$, $\varphi_i(\widetilde{e}_ib)=\varphi_i(b)+1$ if $\widetilde{e}_ib\neq 0$.
		\end{itemize}
		
		\vskip 2mm
		
		\item[{\rm (e)}] For any $i\in I^{\text{im}}$, $l>0$ and $b\in B$, we have
		\vskip 2mm
		\begin{itemize}
			\item [($1'$)] $\varepsilon_i(\widetilde{f}_{il}b)=\varepsilon_i(b)$, $\varphi_i(\widetilde{f}_{il}b)=\varphi_i(b)-la_{ii}$ if $\widetilde{f}_{il}b\neq 0$,
			\vskip 2mm
			\item [($2'$)] $\varepsilon_i(\widetilde{e}_{il}b)=\varepsilon_i(b)$, $\varphi_i(\widetilde{e}_{il}b)=\varphi_i(b)+la_{ii}$ if $\widetilde{e}_{il}b\neq 0$.
		\end{itemize}
		
		\vskip 2mm
		
		\item[{\rm (f)}]  For any $(i,l)\in I^\infty$ and $b\in B$ such that $\varphi_i(b)=-\infty$, we have $\widetilde{e}_{il}b=\widetilde{f}_{il}b=0$.
	\end{enumerate}
\end{definition}
Let $ B_1$ and $ B_2$ be abstract crystals and let
$ B_1 \otimes B_2 = \{ b_1 \otimes b_2 \mid b_1 \in B_1, b_2 \in B_2 \}$.  Define the maps $\text{wt}$, $\varepsilon_i$, $\varphi_i$ $(i \in I)$, $\widetilde{e}_{il}$, $\widetilde{f}_{il}$ $((i,l) \in I^{\infty})$ as follows.
\begin{equation} \label{eq:wt}
	\begin{aligned}
		& \text{wt}(b_1\otimes b_2)=\text{wt}(b_1)+\text{wt}(b_2),\\
		& \varepsilon_i(b_1\otimes b_2)=\text{max}(\varepsilon_i(b_1),\varepsilon_i(b_2)- \langle h_i, \text{wt} (b_1) \rangle),\\
		& \varphi_i(b_1\otimes b_2)=\text{max}(\varphi_i(b_1)+\langle h_i, \text{wt}(b_2) \rangle,\varphi_i(b_2)).
	\end{aligned}
\end{equation}

If $i \in I^{\text{re}}$,
\begin{equation} \label{eq:real}
	\begin{aligned}
		& \widetilde{e}_{i}(b_1 \otimes b_2) = \begin{cases} \widetilde{e}_{i} b_1 \otimes b_2  & \text{if} \  \varphi_{i}(b_1) \ge \varepsilon_{i}(b_2), \\
			(-1)^{p(i)p(b_{1})}b_1 \otimes \widetilde{e}_{i} b_2  & \text{if} \ \varphi_{i}(b_1) < \varepsilon_{i}(b_2),
		\end{cases} \\
		& \widetilde{f}_{i}(b_1 \otimes b_2) = \begin{cases} \widetilde{f}_{i} b_1 \otimes b_2  & \text{if} \ \varphi_{i}(b_1) > \varepsilon_{i}(b_2), \\
		(-1)^{p(i)p(b_{1})}	b_1 \otimes \widetilde{f}_{i} b_2  & \text{if} \ \varphi_{i}(b_1) \le \varepsilon_{i}(b_2).
		\end{cases}
	\end{aligned}
\end{equation}

If $i \in I^{\text{im}}$,
\begin{equation} \label{eq:im}
	\begin{aligned}
		& \widetilde{e}_{il}(b_1 \otimes b_2) =  \begin{cases} \widetilde{e}_{il} b_1 \otimes b_2  & \text{if} \ \varphi_{i}(b_1) > \varepsilon_{i}(b_2) - l a_{ii}, \\
			0  & \text{if} \ \varepsilon_{i}(b_{2}) < \varphi_{i}(b_{1}) \le \varepsilon_{i}(b_{2}) - l a_{ii}, \\
		(-1)^{p(li)p(b_{1})}	b_{1} \otimes \widetilde{e}_{il} b_2 \ \ & \text{if} \ \varphi_{i}(b_1) \le \varepsilon_{i}(b_{2}),
		\end{cases}\\
		& \widetilde{f}_{il}(b_1 \otimes  b_2) =\begin{cases} \widetilde{f}_{il} b_1 \otimes b_2  & \text{if} \ \varphi_{i}(b_1) > \varepsilon_{i}(b_2), \\
		(-1)^{p(li)p(b_{1})}	b_{1} \otimes \widetilde{f}_{il} b_2  & \text{if} \ \varphi_{i}(b_1) \le \varepsilon_{i}(b_2).
		\end{cases}
	\end{aligned}
\end{equation}
\vskip 3mm
As  the same argument of Proposition 3.1 in \cite{FKKT21}, we have the following proposition.
\begin{proposition}\label{tensor product} 
If $B_1$ and $B_2$ are abstract crystals, then $B_1\otimes B_2$ defined in \eqref{eq:wt}--\eqref{eq:im}	is also an abstract crystal.
	\end{proposition}

Let $V(\lambda)$ be an irreducible highest weight module with highest weight $\lambda\in P^{+}$  and let $( L, B)$ be a crystal basis of $M$.
We already have the maps
\begin{align}\label{def: wt and ef}
\text{wt}: B\rightarrow P,\quad \widetilde{e}_{il}, \widetilde{f}_{il}: B\rightarrow B\cup\{0\}.
\end{align}

 Define
 \begin{equation}\label{def:epsilon and phi}
\begin{aligned}
&\epsilon_i(b)=
\begin{cases}
\text{max}\{k\geq0\mid \widetilde{e}_i^kb\neq 0\}  &\text{for} \ i\in I^{\text{re}},\\
 0 & \text{for} \ i\in I^{\text{im}},
\end{cases}\\
&\phi_i(b)=\epsilon_i(b)+ \langle h_i, \text{wt}(b) \rangle \quad \text{for any} \ \  i\in I.
\end{aligned}
\end{equation}
By the same argument of Lemma 3.16 in \cite{FHKX}, we have the following lemma.
\begin{lemma} {
The set $B$ together with the maps defined in \eqref{def: wt and ef}--\eqref{def:epsilon and phi} is an abstract crystal.}
\end{lemma}

\vskip 3mm
\subsection{Crystal bases for $U_{q}(\g)^{-}$}

\begin{lemma} \label{4.121}
\begin{itemize}
 \item [(1)] $\varrho_{il}(a_{i,\mathbf{c}})=\sum_{k:c_{k}=l}(-1)^{p(li)p(\sum_{j>k}c_{j})}q_{(i)}^{2l\sum_{j>k}c_{j}}a_{i,\mathbf{c}\backslash c_{k}},$
 \item [(2)] $a_{il}x^{-}-(-1)^{p(li)p(x)}x^{-}a_{il}=\frac{K_{i}^{l}\varrho^{il}(x)^{-}-(-1)^{p(li)p(x)-p(li)}\varrho_{il}(x)^{-}K_{i}^{-l}}{1-(-1)^{p(li)}q_{i}^{2l}}$ for any $x\in U^{+}$.
 \item[(3)] $\varrho_{il}(xy)=(-1)^{p(li)p(y)}q_{(l\alpha_{i},|y|)}\varrho^{il}(x)y+x\varrho_{il}(y)$
\end{itemize}
\end{lemma}
\begin{lemma}\label{4.12}
    Define $e_{il}':U^{-}\to U^{-}$ by $e_{il}'(z^{-})=(-1)^{p(li)p(y)-p(li)}q_{i}^{l\langle h_{i},|y|+l\alpha_{i}\rangle}\varrho_{il}(y^{+})^{-}$ for any $z\in U^{+}$. Set $$\mathcal{K}_{i}=\cap_{l>0} \text{ker}\ e_{il}'$$
    for any $i\in I^{im}$. We have the following decomposition
    $$U^{-}=\oplus_{\mathbf{c}\in \mathcal{C}_{i}}b_{i,\mathbf{c}}\mathcal{K}_{i}.$$
\end{lemma}
\begin{proof}
    Let $u\in U^{-}$, assume first that $u$ has the form $u=mb_{i,\mathbf{c}}m'$ where $m,m' \in \mathcal{K}_{i}$ and $\mathbf{c}\in \mathcal{C}_{i}$. We use induction on $|\mathbf{c}|$. When $|\mathbf{c}|=0$, we have $mm'\in \mathcal{K}_{i}$. When $|\mathbf{c}|>0$, set $[y,z]^{\circ}=(-1)^{p(y)p(z)}q^{(|y|,|z|)}yz-zy$ for any $y,z\in U^{+}$. Since $\varrho_{i,l}(a_{ik})=\delta_{lk}$, for any $y\in \cap_{l>0} \text{ker} \varrho_{il}$ and any $k>0$, we have 
    \begin{align*}
        \varrho_{il}([y,a_{ik}]^{\circ})&=\varrho^{il}((-1)^{p(y)p(a_{ik})}q^{(|y|,k\alpha_{i})}ya_{ik}-a_{ik}y)\\
        &=0.
    \end{align*}
    Since $|\mathbf{c}\backslash c_{1}|<|\mathbf{c}|$ and $\oplus_{\mathbf{c}\in \mathcal{C}_{i}}b_{i,\mathbf{c}}\mathcal{K}_{i}$ is stable by left multiplication by $b_{i,c_{1}}$, by the induction hypothesis, we can conclude that 
    $$u=(-1)^{p(m)p(c_1i)}q^{-c_{1}(|m|,\alpha_{i})}[m,b_{ic_{1}}]^{\circ}b_{i,\mathbf{c}\backslash c_{1}}m'+q^{c_{1}(|m|,\alpha_{i})}b_{ic_{1}}mb_{i,\mathbf{c}\backslash c_{1}}m'$$
   has the decomposition.  
   
   Next, We consider the general homogeneous $u\in U^{-}$. Using induction on $\text{ht}(u)$. when $\text{ht}(u)>1$ we can write 
   $$u=\sum_{\iota\in I^{\infty}}b_{\iota}u_{\iota}$$
   for some finitely many nonzero $u_{\iota}\in U^{-}$. By induction, We have $$u=\sum_{\iota\in I^{\infty}} b_{\iota}b_{i,\mathbf{c}}z_{\iota,\mathbf{c}}$$
   for some finitely many  nonzero $z_{\iota,\mathbf{c}}\in \mathcal{K}_{i}$. Then
   $$u=\sum_{l>0,\mathbf{c}\in \mathcal{C}_{i}} b_{i,(l,\mathbf{c})}z_{(l,\mathbf{c})}+\sum_{\iota\notin\{(i,\N_{>0})\}}b_{\iota}b_{i,\mathbf{c}}z_{\iota,\mathbf{c}},$$
   and we have the result since $b_{\iota},z_{\iota,\mathbf{c}}\in \mathcal{K}_{i}$.

   Finally, we will prove the uniqueness of the decomposition. Consider a minimal relation of dependence
   $$0=\sum_{\mathbf{c}\in \mathcal{C}_{i}}a_{i,\mathbf{c}}z_{\mathbf{c}}$$
   for some finitely many nonzero $z_{\mathbf{c}}\in \mathcal{K}_{i}^{-}$.
   First we consider the case $i\notin I^{iso}$. Consider $r$ maximal such that there exists $\mathbf{c}=(c_{1},\cdots ,c_{r})$ such that $z_{\mathbf{c}}\neq 0$. Define $\varrho_{i,\mathbf{c}}=\varrho_{i,c_r}\circ\cdots\circ\varrho_{i,c_{1}}$, for any $\mathbf{c}'\in \mathfrak{S}_{r}\mathbf{c}$,
   \begin{align*}
       0&=\varrho_{i,\mathbf{c}'}(\sum_{\mathbf{c}''\in\mathcal{C}_{i}}a_{i,\mathbf{c}''}z_{\mathbf{c}''})\\
       &=\sum_{\mathbf{c}''\in\mathcal{C}_{i}}(-1)^{p(|c|i)p(z_{\mathbf{c}''})}q^{(|\mathbf{c'}|\alpha_{i},|z_{\mathbf{c}''}|)}\varrho_{i,\mathbf{c}'}(a_{i,\mathbf{c}''})z_{\mathbf{c}''}\\
&=(-1)^{p(|c|i)p(z_{\mathbf{c}})}q^{(|\mathbf{c}|\alpha_{i},|z_{\mathbf{c}}|)}\sum_{\mathbf{c}''\in\mathfrak{S}_{r}\mathbf{c}}\varrho_{i,\mathbf{c}'}(a_{i,\mathbf{c}''})z_{\mathbf{c}''}\\
&=(-1)^{p(|c|i)p(z_{\mathbf{c}})}q^{(|\mathbf{c}|\alpha_{i},|z_{\mathbf{c}}|)}\sum_{\mathbf{c}''\in\mathfrak{S}_{r}\mathbf{c}}P_{\mathbf{c}',\mathbf{c}''}(q_{(i)})z_{\mathbf{c}''},
   \end{align*}
which yields $\sum_{\mathbf{c}''\in\mathfrak{S}_{r}\mathbf{c}}P_{\mathbf{c}',\mathbf{c}''}(q_{(i)})z_{\mathbf{c}''}=0$, where $P_{\mathbf{c}',\mathbf{c}''}(q)\in \Z[q_{(i)}]$. Since $(z_{\mathbf{c}''})_{\mathbf{c}''\in \mathfrak{S}_{r}\mathbf{c}}\neq0$, we have to prove that 
$$\Delta(q)=\text{det}(P_{\mathbf{c}',\mathbf{c}''}(q))_{\mathbf{c}',\mathbf{c}''\in \mathfrak{S}_{r}\mathbf{c}}\neq0\in \Z[q_{(i)}].$$
Using Lemma \ref{4.121} (1) repeatly, we know that 
$$\text{max}_{\mathbf{c}''\in \mathfrak{S}_{r}\mathbf{c}}{\text{deg}(P_{\mathbf{c}',\mathbf{c}''})}=\text{deg}(P_{\mathbf{c}',\mathbf{c}'})=\sum_{p<q}c_{p}'c_{q}'=m,$$ which is only reached for $\mathbf{c}''=\mathbf{c}'$. 
If $m=0$, the initial relation of
dependance can be written: $0=\sum_{l>0}a_{il}z_{l}$. Applying $\varrho_{il}$,
we can conclude that $z_{l}=0$ and get a contradiction. Hence $m\neq 0$ and $\text{deg}(\Delta)=|\mathfrak{S}_{r}\mathbf{c}|m\neq0$.

If $i\in I^{iso}$, write a relation of dependence of minimal degree
$$0=\sum_{\mathbf{c}\in \mathcal{C}_{i}}a_{i,\mathbf{c}}z_{\mathbf{c}},$$
where $z_{\mathbf{c}}\in \mathcal{K}_{i}^{-}$. For any $\mathbf{c}$ and $l>0$, set $\mathbf{c}_{l}=\{s:c_{s}=l\}$. If $p(li)=1$, by Definition \ref{def:BBsuper}, $a_{il}^{2}=0$,  we can assume that when $p(li)=1$,  $\mathbf{c}_{l}=1$. Applying $\varrho_{il}$, we have 
\begin{align*}
    \varrho_{il}(\sum_{\mathbf{c}\in \mathcal{C}_{i}}a_{i,\mathbf{c}}z_{\mathbf{c}})&=\sum_{\mathbf{c}\in \mathcal{C}_{i}}(-1)^{t}q^{(|\mathbf{c}|\alpha_{i},|z_{\mathbf{c}}|)}a_{i,\mathbf{c}\backslash l}z_{\mathbf{c}}\\
    &=0
\end{align*}
for some $t\in \{0,1\}$,which contradicts the minimality of the relation. When $p(li)=0$, we have 
\begin{align*}
     \varrho_{il}(\sum_{\mathbf{c}\in \mathcal{C}_{i}}a_{i,\mathbf{c}}z_{\mathbf{c}})&=\sum_{\mathbf{c}\in \mathcal{C}_{i}}(-1)^{t}q^{(|\mathbf{c}|\alpha_{i},|z_{\mathbf{c}}|)}\mathbf{c}_{l}a_{i,\mathbf{c}\backslash l}z_{\mathbf{c}}\\
     &=\sum_{\lambda\in\mathcal{C}_{i}}a_{i,\lambda}(-1)^{t}q^{(|\lambda\cup l|\alpha_{i},|z_{\lambda\cup l}|)}(\lambda_{l}+1)z_{\lambda\cup l}\\
     &=0,
\end{align*}
which also contradicts the minimality of the relation.
\end{proof}
\begin{definition} \label{def:Kashiwara operator} \hfill

\vskip 3mm

(a) For $i \in I^{\text{re}}$, we define
\begin{equation} \label{eq:real Kas}
\begin{aligned}
& \widetilde{e}_{i} u = \sum_{k \ge 1} {\mathtt b}_{i}^{(k-1)} u_{k}, \\
& \widetilde{f}_{i}  u = \sum_{k \ge 0} {\mathtt b}_{i}^{(k+1)} u_{k}.
\end{aligned}
\end{equation}

\vskip 2mm

(b) For $i \in I^{\text{im}}\setminus I^{\text{iso}}$ and $l>0$, we define
\begin{equation} \label{eq: nonisoKas}
\begin{aligned}
& \widetilde{e}_{il} u = \sum_{\mathbf{c} \in {\mathcal C}_{i}:  c_1 = l} \,
{\mathtt b}_{i, \mathbf{c} \setminus c_1} u_{\mathbf{c}},\\
& \widetilde{f}_{il}  u = \sum_{\mathbf{c} \in {\mathcal C}_{i}}
{\mathtt b}_{i,(l, \mathbf{c})} u_{\mathbf{c}}.
\end{aligned}
\end{equation}

\vskip 2mm

(c) For $i \in I^{\text{iso}}$ and $l>0$, we define
\begin{equation} \label{eq:isoKas}
\begin{aligned}
& \widetilde{e}_{il} u = \sum_{\mathbf{c} \in {\mathcal C}_{i}} \,  (-1)^{p(li)p(|\mathbf{c}_{<r}|i)} \mathbf{c}_{l} \,
{\mathtt b}_{i, \mathbf{c} \setminus \{l\}} u_{\mathbf{c}},\\
& \widetilde{f}_{il}  u = \sum_{\mathbf{c} \in {\mathcal C}_{i}} \, \frac{1} {\mathbf{c}_{l} + 1} \,
{\mathtt b}_{i, (l,\mathbf{c})} u_{\mathbf{c}},
\end{aligned}
\end{equation}
where ${\mathbf c}_{l}$ denotes the number of $l$ in $\mathbf{c}$, $r=\min\{t|c_{t}=l\}$.

\end{definition}

\begin{definition} \label{def:crystal lattice U}

\vskip 2mm

 A free $\A_0$-submodule $ L$ of $U_{q}^-(\g)$ is called a {\it crystal lattice} if the following conditions hold.

 \vskip 2mm

\begin{enumerate}

\vskip 2mm

\item[{\rm (a)}] $\Q(q)\otimes_{\mathbf A_0} L\cong U_q^-(\g)$,

\vskip 2mm

\item[{\rm (b)}] $ L=\oplus_{\alpha\in Q_+} L_{-\alpha}$, where $ L_{-\alpha}= L\cap {U_{q}^-(\g)}_{-\alpha}$,

\vskip 2mm

\item[{\rm (c)}] $\widetilde e_{il} L\subset  L$, \ $\widetilde f_{il} L\subset  L$ for all $(i,l)\in I^{\infty}$.

\end{enumerate}

\end{definition}

\vskip 3mm

The condition (c) yields  the $\Q$-linear maps
$${\widetilde e}_{il}, \, {\widetilde f}_{il}: L/qL \longrightarrow L/qL.$$

\vskip 3mm

\begin{definition}\label{4.15}

A {\it crystal basis} of $U_{q}^-(\g)$ is a pair $( L, B)$ such that
	
\vskip 2mm
	
\begin{enumerate}
	
\vskip 2mm
	
\item[{\rm (a)}] $L$ is a crystal lattice of $U_{q}^-(\g)$,
		
\vskip 2mm

\item[{\rm (b)}]  $B$ is a pseudo-basis of $L/qL$, i.e. there exists a $\Q$-basis $B^{\circ}$ of $L/qL$ such that $B^{\circ}\subset B\subset B^{\circ}\cup (-B^{\circ}) $,

\vskip 2mm

\item[{\rm (c)}] $ B=\sqcup_{\alpha \in Q_+} B_{-\alpha}$, where
$ B_{-\alpha}= B\cap{( L/q L)}_{-\alpha}$,

\vskip 2mm

\item[{\rm (d)}] $\widetilde{e}_{il} B\subset B\cup\{0\}$, \ $\widetilde{f}_{il} B \subset B\cup\{0\}$
for $(i,l) \in I^{\infty}$,

\vskip 2mm

\item[{\rm (e)}] for any $b,b'\in B$ and $(i,l) \in I^{\infty}$,
we have $\widetilde{f}_{il}b=b'$ if and only if $b=\widetilde{e}_{il}b'$.

\end{enumerate}

\end{definition}

\vskip 3mm

Let $L(\infty)$ be the $\A_0$-submodule of $U_{q}^-(\g)$ spanned by $\widetilde{f}_{i_1l_1}\cdots\widetilde{f}_{i_rl_r}\mathbf 1$ $(r\geq 0, (i_j,l_j)\in I_{\infty})$,
and $B(\infty)=\{\widetilde{f}_{i_1l_1}\cdots\widetilde{f}_{i_rl_r}\mathbf 1+qL(\infty)\}\backslash \{0\}$.  Set $B(\infty)=\mathcal{B}(\infty)\cup (-\mathcal{B}(\infty))$. Clearly, $B(\infty)=\sqcup_{\alpha\in Q_+} B(\infty)_{-\alpha}$ and $B(\infty)_{-\alpha}$ is a finite set. There exists a minimal subset $B(\infty)_{-\alpha}^{\circ}\subset \mathcal{B}(\infty)_{-\alpha}$ such that $B(\infty)_{-\alpha}=B(\infty)_{-\alpha}^{\circ}\cup(-B(\infty)_{-\alpha}^{\circ})$. Define $B(\infty)^{\circ}=\sqcup_{\alpha\in Q_+}B(\infty)_{-\alpha}^{\circ}$.

\vskip 3mm

\begin{theorem}\label{thm:crystal basis of U}
{\rm
The pair  $(L(\infty),B(\infty))$ is a crystal basis of  $U^-_q(\mathfrak g)$.
}
\end{theorem}

\vskip 3mm

We will prove this theorem in Section. 

\vskip 3mm

\section{Grand-loop argument}

Let $(i,l) \in I^{\infty}$, $\lambda,\mu\in Q_+$ and  $\alpha\in  R_+(r)$, where
$R_+(r)=\{\alpha\in Q_+\mid |\alpha|\leq r \}$.
\vskip 2mm
\begin{itemize}

\item[$\mathbf{A}(r)$:]  $\widetilde{e}_{il}L(\lambda)_{\lambda-\alpha} \subset L(\lambda)$, \
 $\widetilde{e}_{il}B(\lambda)_{\lambda-\alpha} \subset B(\lambda) \cup \{0\}$.

 \vskip 2mm

\item[$\mathbf{B}(r)$:] For $b \in B(\lambda)_{\lambda - \alpha + l \alpha_{i}}$,
$b' \in B(\lambda)_{\lambda - \alpha}$,\,  $\widetilde{f}_{il} b = b'$ if and only if
$\widetilde{e}_{il} b' = b$.

\vskip 2mm

\item[$\mathbf{C}(r)$:] $\Phi_{\lambda, \mu} (L(\lambda + \mu)_{\lambda + \mu - \alpha})
\subset L(\lambda) \otimes L(\mu)$.

\vskip 2mm

\item[$\mathbf{D}(r)$:] $\Psi_{\lambda, \mu}((L(\lambda) \otimes L(\mu))_{\lambda + \mu - \alpha})
\subset L(\lambda + \mu)$, \
 $\Psi_{\lambda, \mu}((B(\lambda) \otimes B(\mu))_{\lambda + \mu - \alpha})\subset B(\lambda + \mu) \cup \{0\}.$

\vskip 2mm

\item[$\mathbf{E}(r)$:] $\widetilde{e}_{il}L(\infty)_{-\alpha} \subset L(\infty)$, \
 $\widetilde{e}_{il}B(\infty)_{-\alpha} \subset B(\infty) \cup \{0\}$.

\vskip 2mm

\item[$\mathbf{F}(r)$:] For $b \in B(\infty)_{ - \alpha + l \alpha_{i}}$,
$b' \in B(\infty)_{- \alpha}$, \, $\widetilde{f}_{il} b = b'$ if and only if
$\widetilde{e}_{il} b' = b$.

\vskip 2mm

\item[$\mathbf{G}(r)$:] $B(\lambda)_{\lambda - \alpha}^{\circ}$ is  a $\Q$-basis of
$(L(\lambda) / q L(\lambda))_{\lambda - \alpha}$, \  $B(\infty)_{- \alpha}^{\circ}$ is a $\Q$-basis of
$(L(\infty) / q L(\infty))_{ - \alpha}$.

\vskip 2mm

\item[$\mathbf{H}(r)$:] $\pi_{\lambda}(L(\infty)_{-\alpha}) = L(\lambda)_{\lambda - \alpha}$.

\vskip 2mm

\item[$\mathbf{I}(r)$:] For $S \in L(\infty)_{-\alpha + l \alpha_{i}}$,\,
$\widetilde{f}_{il}(S\, v_{\lambda}) \equiv (\widetilde{f}_{il} S)\, v_{\lambda}\ \text{mod} \ q L(\lambda)$.

\vskip 2mm

\item[$\mathbf{J}(r)$:] If $B^{\lambda}_{-\alpha}:= \{ b\in B(\infty)_{-\alpha}
\mid \overline{\pi}_{\lambda}(b) \neq 0 \}$,
then $B^{\lambda}_{-\alpha}  \cong B(\lambda)_{\lambda - \alpha}$.

\vskip 2mm

\item[$\mathbf{K}(r)$:] If $b \in B^{\lambda}_{-\alpha}$, then
$\widetilde{e}_{il} \, \overline{\pi}_{\lambda} (b) = \overline{\pi}_{\lambda}\, \widetilde{e}_{il} (b)$.

\end{itemize}

We will prove the statements $\mathbf{A}(r),\cdots \mathbf{K}(r)$ by induction. Note that the above statements are true for $r=0$ and $r=1$, we will assume that $r\ge 2$ and $\mathbf{A}(r-1),\cdots,\mathbf{K}(r-1)$ are true

Let $\mathbf{c}=(c_{1},\cdots,c_k)\in \mathcal{C}_{i}$, $|c|=\sum_{t=1}^{k}c_r$. If $i\in I^{re}$, $b_{i,\mathbf{c}}=b_{i}^{(|\mathbf{c}|)}$. If $i\in I^{im}$, $b_{i,\mathbf{c}}=b_{ic_{1}}\cdots b_{ic_{k}}$. 

\begin{lemma}\label{alpha and b}
{
Let $\alpha\in R_+(r-1)$ and $b\in {B(\lambda)}_{\lambda-\alpha}$.
If  $\widetilde{e}_{il}b=0$ for any $(i,l)\in I^{\infty}$, then we have $\alpha=0$ and $b=v_\lambda$.
}
\end{lemma}

\begin{proof}
If $\alpha\neq 0$, then there exist $(i,l)\in I^{\infty}$ and $b'\in B(\lambda)_{\lambda-\alpha+l\alpha_{i}}$ such that $b=\widetilde{f}_{il}b'$. Then from $\mathbf{B}(r-1)$, we have $\widetilde{e}_{il}b=b'\neq 0$. We get a contradiction, hence $\alpha=0$ and $b=v_{\lambda}$.
	\end{proof}

\vskip 3mm

\begin{lemma}\label{uc in L}

{ 	Let $\alpha \in R_{+}(r-1)$, $i\in I$, and
$u = \sum_{\mathbf{c} \in \mathcal{C}_{i}} {\mathtt b}_{i, \mathbf{c}}
u_{\mathbf{c}} \in V(\lambda)_{\lambda-\alpha}$
be the $i$-string decomposition of $u$.
If $u\in L(\lambda)$, then $u_{\mathbf c}\in L(\lambda)$ for any $\mathbf c\in\mathcal C_i$.	
}
\end{lemma}

\begin{proof}
	Suppose $u=\sum_{\mathbf c\in\mathcal C_i}{\mathtt b_{i,\mathbf c}u_{\mathbf c}\in L(\lambda)}$ for finitely many nozero $u_{\mathbf{c}}$, let $N=\text{max}\{|\mathbf{c}||u_{\mathbf{c}}\neq 0\}$. We shall use the induction on $N$.
	If $N=0$, the assertion follows naturally.
	If $N>0$, by $\mathbf{A}(r-1)$, we have  $\widetilde{e}_{il}u\in L(\lambda)$ for any $l>0$. By Definition \ref{def:Kashiwara operator}, we have
	
	\begin{equation*}
		\widetilde{e}_{il}u=\begin{cases}
        \sum_{m\ge 1} b_{i}^{(m-1)}u_{m}\in L(\lambda),&\text{if} \ i\in I^{\text{re}},\\
			\sum_{\mathbf c:c_1=l}{ b_{i,\mathbf c\backslash c_1}u_{\mathbf c}\in L(\lambda)},& \text{if}\ i\in I^{\text{im}}\setminus I^{\text{iso}},\\
			\sum_{\mathbf c\in\mathcal C_i} (-1)^{p(li)p(\mathbf{|c}_{<k}|i)} {\mathbf c}_{l}  { b_{i,\mathbf c\backslash \{l\}}u_{\mathbf c}\in L(\lambda)},& \text{if}\ i\in I^{\text{iso}}.
		\end{cases}
	\end{equation*}
	 By the induction hypothesis, $u_{\mathbf c}\in L(\lambda)$ for any $\mathbf c\ne\mathbf 0$.
	
	Let $u_{1}:=\sum_{\mathbf c\ne\mathbf 0}{\mathtt b_{i,\mathbf c}u_{\mathbf c}}$. It is clear that
	$u_1\in L(\lambda)$. Hence $u_{ 0}:=u-u_1\in L(\lambda)$, which proves this lemma.
\end{proof}

\vskip 3mm

\begin{lemma}\label{uc=0}
{
	Let $\alpha\in R_{+}(r-1)$, $i\in I$ . Let $u=\sum_{\mathbf c\in\mathcal C_i}{\mathtt b_{i,\mathbf c}u_{\mathbf c}\in {V(\lambda)}_{\lambda-\alpha}}$ for finitely nonzero $u_{\mathbf{c}}$ be the $i$-string decomposition of $u$ and  $N=\text{max}\{|\mathbf{c}||u_{\mathbf{c}}\neq 0\}$. 
	If $ u+qL(\lambda)\in B(\lambda)$, then there exists $\mathbf c\in\mathcal C_i$ such that
	\begin{enumerate}
		\item[{\rm (a)}] $u\equiv	b_{i,\mathbf c}u_{\mathbf c}\!\!\! \mod{qL(\lambda)}$,
		\item[{\rm (b)}] $u_{\mathbf c'}\equiv 0\!\!\! \mod{qL(\lambda)}$ for any $\mathbf c'\ne\mathbf c$.
	\end{enumerate}
	}
\end{lemma}

\begin{proof}
The case for  $N=0$ is trivial. For $N>0$, by $\mathbf{A}(r-1)$, we have $\widetilde{e}_{il}(u+qL(\lambda))\in B(\lambda)\cup \{0\}$ for any $l>0$.

If $\widetilde{e}_{il}u\in qL(\lambda)$ for any $l>0$, by Lemma \ref{uc in L}, we have
	$u_{\mathbf c}\in qL(\lambda)$ for any $\mathbf c\ne\mathbf 0$.
	Then $u\equiv u_{0}\!\!\! \mod{qL(\lambda)}$. 
	
	If $\widetilde{e}_{il}u\notin  qL(\lambda)$ for some $l>0$, then $\widetilde{e}_{il}u+qL(\lambda)\in B(\lambda)$. By induction hypothesis, If $i\notin I^{\text{iso}}$, there exists $\mathbf c\in \mathcal C_i$ with $c_{1}=l$ such that $\widetilde{e}_{il}u\equiv  b_{i,\mathbf{c}\backslash c_1}u_{\mathbf{c}} \mod{qL(\lambda)}$.

If $i\in I^{\text{iso}}$, there exists $\mathbf c\in \mathcal C_i$ with $c_{t}=l$, such that $\widetilde{e}_{il}u\equiv(-1)^{p(li)p(\mathbf{|c}_{<k}|i)}\mathbf{c}_{l}b_{i,\mathbf{c}\backslash \{l\}}u_{\mathbf{c}} \mod{qL(\lambda)}$, where $k=\text{min}\{t|c_{t}=l\}$. Note that $$\widetilde{f}_{il}((-1)^{p(li)p(\mathbf{|c}_{<k}|i)}\mathbf{c}_{l}b_{i,\mathbf{c}\backslash \{l\}}u_{\mathbf{c}})=(-1)^{p(li)p(\mathbf{|c}_{<k}|i)}b_{i,(l,\mathbf{c}\backslash \{l\})}u_{\mathbf{c}}=b_{i,\mathbf{c}}u_{\mathbf{c}}.$$

Hence by $\mathbf{B}(r-1)$, we have 
$$u\equiv \widetilde{f}_{il}\widetilde{e}_{il}u\equiv b_{i,\mathbf{c}}u_{\mathbf{c}}\mod{qL(\lambda)}.$$
By Lemma \ref{uc in L}, $u_{\mathbf{c}'}\in qL(\lambda)$ for $\mathbf{c}'\neq \mathbf{c}$.
\end{proof}

\vskip 2mm

\begin{lemma}
{
	Let $\alpha\in  R_{+}(r-1)$, $i\in I$ and let $u=\sum_{\mathbf c\in\mathcal C_i}{\mathtt b_{i,\mathbf c}u_{\mathbf c}}\in  {U^{-}_q(\g)}_{-\alpha}$ be the $i$-string decomposition of $u$.
	\begin{enumerate}
		\item[{\rm (a)}] If $u\in L(\infty)$, then $u_{\mathbf c}\in L(\infty)$ for any $\mathbf c$.
		\item[{\rm (b)}] If $u+qL(\infty)\in B(\infty)$, then there exists $\mathbf c\in\mathcal C_i$ such that
		\begin{enumerate}
			\item[{\rm (1)}] $u\equiv\widetilde{f}_{i,\mathbf c}u_{\mathbf c}\!\! \mod{qL(\infty)}$,
			\item[{\rm (2)}] $u_{\mathbf c'}\equiv 0\!\! \mod{qL(\infty)}$ for any $\mathbf c'\neq\mathbf c$.
		\end{enumerate}
	\end{enumerate}	

}
\end{lemma}
\begin{proof}
    The proof is similar to that of the preceding lemma.
\end{proof}

\vskip 2mm

\begin{lemma}\label{a-g}

{
	Let $\alpha, \beta\in R_{+}(r-1)$ and $i\in I$.
	\begin{enumerate}
		\item[{\rm (a)}] For all $l>0$, we have
		\begin{align*}
			&	\widetilde{e}_{il}(L(\lambda)_{\lambda-\alpha}\otimes L(\mu)_{\mu-\beta})\subset L(\lambda)\otimes L(\mu),\\	
			&\widetilde{f}_{il}(L(\lambda)_{\lambda-\alpha}\otimes L(\mu)_{\mu-\beta})\subset L(\lambda)\otimes L(\mu).			
		\end{align*}
		\item[{\rm (b)}]
		If $\widetilde{e}_{il}(b\otimes b')\ne 0$, then $b\otimes b'=\widetilde{f}_{il}\widetilde{e}_{il}(b\otimes b')$.
		
		\vskip 2mm
		
		\item[{\rm (c)}]
		If $\widetilde{e}_{il}(b\otimes b')= 0$ for all $l>0$, then $b=v_{\lambda}$.
		
		\vskip 2mm
		
		\item[{\rm (d)}]
        For any $(i,l)\in I^{\infty}$, we have $\widetilde{f}_{il}(b\otimes v_{\mu})=
			\widetilde{f}_{il}b\otimes v_{\mu}$ or $\widetilde{f}_{il}b=0$.	
			
			\vskip 2mm
			
		\item[{\rm (e)}]
		For any $(i_1,l_1),\cdots,(i_r,l_r)\in I^{\infty}$, we have
		\begin{equation*}
			\widetilde{f}_{i_1l_1}\cdots	\widetilde{f}_{i_rl_r}(v_{\lambda}\otimes v_{\mu})\equiv \widetilde{f}_{i_1l_1}\cdots \widetilde{f}_{i_rl_r}v_{\lambda}\otimes v_{\mu}\!\!\!\! \mod{q(L(\lambda)\otimes L(\mu))}
		\end{equation*}
		or $\widetilde{f}_{i_1l_1}\cdots\widetilde{f}_{i_rl_r}v_{\lambda}\in qL(\lambda)$.
	\end{enumerate}
}
\end{lemma}
\begin{proof}
(a)From Lemma \ref{uc=0}, it suffices to prove the following statement: for $u\in \text{ker}(a_{il})\cap L(\lambda)_{\lambda-\alpha+|\mathbf{c}_{1}|\alpha_{i}}$ and 
 $v\in \text{ker}(a_{il})\cap L(\mu)_{\mu-\beta+|\mathbf{c}_{2}|\alpha_{i}}$
, we have 
\begin{align*}
  &\widetilde{f}_{il}(\widetilde{f}_{i,\mathbf{c}_{1}}u\otimes \widetilde{f}_{i,\mathbf{c}_{2}}v)\in L(\lambda)\otimes L(\mu),\\ 
   &\widetilde{e}_{il}(\widetilde{f}_{i,\mathbf{c}_{1}}u\otimes \widetilde{f}_{i,\mathbf{c}_{2}}v)\in L(\lambda)\otimes L(\mu).
\end{align*}
    Let $L$ be the free $\A_{0}$-submodule of $V(\lambda)\otimes V(\mu)$
 generated by $\widetilde{f}_{i,\mathbf{c}}u\otimes \widetilde{f}_{i,\mathbf{c}'}v$. By the tensor product rule, we have 
 $\widetilde{f}_{il}L\subset L$ and $\widetilde{e}_{il}L\subset L$. Since $L\subset L(\lambda)\otimes L(\mu)$, our statement follows.

(b) 
 By tensor product rule, we can directly derived (b)

 (c) Suppose $\widetilde{e}_{il}(b\otimes b')=0$ for any $l>0$. Let $\text{wt}(b)=\lambda-l_{1}\alpha_{i_1}-\cdots-l_{s}\alpha_{i_{s}}$.
Set $m=\langle h_i,\mathrm{wt}(b)\rangle$. If $i\in\{i_1,\cdots,i_s\}$, then  there exists $l>0$ such that	 $-la_{ii}\leq m\leq -(l+1)a_{ii}$.
	By tensor product rule, we have $\widetilde{e}_{ik}(b\otimes b')=\widetilde{e}_{ik}b\otimes b'=0$ for $0<k\leq l$, thus $\widetilde{e}_{ik}b=0$.
  we have  $\widetilde{e}_{ik}b=0$ for any $k>l$, this is because, for any $k\geq l+1$, $\text{wt(b)}+k\alpha_{i}$ is not a weight. Hence, in this situation, $\widetilde{e}_{il}b=0$ for any $l>0$.
	
	If $i\notin\{i_1,\cdots,i_s\}$, then for any $l>0$, $\text{wt(b)}+k\alpha_{i}$ is not a weight, which implies $\widetilde{e}_{il}b=0$ for any $l>0$.
Hence, by Lemma \ref{alpha and b}, we have $b=v_{\lambda}$.

    (d)By tensor product rule, we have $\widetilde{f}_{il}(b\otimes v_{\mu})\equiv \widetilde{f}_{il}b\otimes v_{\mu}$ when $\langle h_i, \text{wt}(b)\rangle>0$. If $\langle h_i, \text{wt}(b)\rangle=0$, by Proposition \ref{2.4}, $\widetilde{f}_{il}b=0$. 
    
    (e) This comes from (d) directly.
\end{proof}

\vskip 2mm

\begin{proposition}  \label{prop:Cr} $({\bf C}(r))$ \
{
	For any $\alpha\in R_{+}(r)$, we have
	\begin{equation*}
		\Phi_{\lambda,\mu}(L(\lambda+\mu)_{\lambda+\mu-\alpha})\subset L(\lambda)\otimes L(\mu).
	\end{equation*}
	
}	
\end{proposition}

\begin{proof}
	We know that
	\begin{equation*}
		{L(\lambda+\mu)}_{\lambda+\mu-\alpha}=\sum_{(i,l)\in I_{\infty}}{\widetilde{f}_{il}({L(\lambda+\mu)}_{\lambda+\mu-\alpha+l\alpha_i})}.
	\end{equation*}
	By ${\bf C}(r-1)$ and Lemma \ref{a-g}, we obtain the result.
\end{proof}

\vskip 2mm

\begin{lemma}\label{ffvvbb}
{
	Let $(i_1,l_1),\cdots,(i_r,l_r)\in I^{\infty}$.
	Suppose that there exists $t$ with $t<r$ satisfying $i_t \ne i_{t+1}=\cdots=i_r$. Then for any $\mu\in P^{+}$ and $\lambda=\Lambda_{i_t}$, we have
	\begin{equation*}
		\widetilde{f}_{i_1l_1}\cdots\widetilde{f}_{i_rl_r}(v_{\lambda}\otimes v_{\mu})\equiv b\otimes b' \!\!\!\!\mod{q(L(\lambda)\otimes L(\mu))}
	\end{equation*}
	for some $b\in {B(\lambda)}_{\lambda-\alpha}\cup\{0\}$, $b'\in {B(\mu)}_{\mu-\beta}\cup\{0\}$ and $\alpha,\beta\in  R_{+}(r-1)$.
}
\end{lemma}

\begin{proof}
Since $\Lambda_{i_t}(h_r)=0$, we have $\mathtt b_{i_rl_r}v_{\lambda}=0$.
Thus for any $v\in V(\mu)$, we get
	\begin{align*}
		\mathtt b_{i_rl_r}(v_{\lambda}\otimes v)=\mathtt b_{i_rl_r}v_{\lambda}\otimes v+K_{i_r}^{l_r}v_{\lambda}\otimes \mathtt b_{i_rl_r}v=v_{\lambda}\otimes \mathtt b_{i_rl_r}v,
	\end{align*}
which gives 
\begin{align*}
    \widetilde{f}_{i_{t+1}l_{t+1}}\cdots\widetilde{f}_{i_rl_{r}}(v_{\lambda}\otimes v_{\mu})&=v_{\lambda}\otimes  \widetilde{f}_{i_{t+1}l_{t+1}}\cdots\widetilde{f}_{i_rl_{r}}v_{\mu}.
\end{align*}
Since $\widetilde{e}_{i_{t}l_t}(v_{\lambda}\otimes \widetilde{f}_{i_{t+1}l_{t+1}}\cdots\widetilde{f}_{i_rl_{r}}v_{\mu})=0$, we have
\begin{align*}
    \widetilde{f}_{i_tl_t}\cdots\widetilde{f}_{i_rl_r}(v_{\lambda}\otimes v_{\mu})&= \widetilde{f}_{i_tl_t}(v_{\lambda}\otimes \widetilde{f}_{i_{t+1}l_{t+1}}\cdots\widetilde{f}_{i_rl_{r}}v_{\mu})\\&=b_{i_tl_t}(v_{\lambda}\otimes \widetilde{f}_{i_{t+1}l_{t+1}}\cdots\widetilde{f}_{i_rl_{r}}v_{\mu})\\
    &\equiv \widetilde{f}_{i_tl_t}v_{\lambda}\otimes\widetilde{f}_{i_{t+1}l_{t+1}}\cdots\widetilde{f}_{i_rl_{r}}v_{\mu} \mod{qL(\lambda)\otimes L(\mu)}
\end{align*}
where $\widetilde{f}_{i_tl_t}v_{\lambda}\in {B(\lambda)}_{\lambda-\alpha}\cup\{0\}$ and $\widetilde{f}_{i_{t+1}l_{t+1}}\cdots	 \widetilde{f}_{i_rl_r}v_{\mu}\in {B(\mu)}_{\mu-\beta}$. Then the lemma follows from the tensor product rule \eqref{eq:im}.
\end{proof}

\vskip 2mm

\begin{lemma}\label{LotimesL}
{
	For any $\alpha\in R_{+}(r)$, we have
	\begin{equation*}
		{(L(\lambda)\otimes L(\mu))}_{\lambda+\mu-\alpha}
		=\sum_{(i,l)\in I_\infty}\mathtt b_{il} {(L(\lambda)\otimes L(\mu))}_{\lambda+\mu-\alpha+l\alpha_i}+v_{\lambda}\otimes {L(\mu)}_{\mu-\alpha}.
	\end{equation*}
}
\end{lemma}
\begin{proof}
    Let $L=(L(\lambda)\otimes L(\mu))_{\lambda+\mu-\alpha}$  and $L'=\sum_{(i,l)\in I_\infty}\mathtt b_{il} {(L(\lambda)\otimes L(\mu))}_{\lambda+\mu-\alpha+l\alpha_i}+v_{\lambda}\otimes {L(\mu)}_{\mu-\alpha}$. By Lemma \ref{a-g}(a), we have $L'\subset L$. By Lemma \ref{a-g}(d), for any $\beta\in Q_{+}(r-1)\backslash \{0\}$ and $b\otimes b'\in B(\lambda)_{\lambda-\beta}\otimes B(\mu)_{\mu-\alpha+\beta}$, there exists $(i,l)\in I^{\infty}$ such that $\widetilde{e}_{il}(b\otimes b')\neq 0$. By Lemma \ref{a-g}(b), we have 
    $$b\otimes b'\equiv \widetilde{f}_{il}\widetilde{e}_{il}(b\otimes b')\mod{qL(\lambda)\otimes L(\mu)}.$$
Hence $L(\lambda)_{\lambda-\beta}\otimes L(\mu)_{\mu-\alpha+\beta}\subset L'+qL$ and we get
\begin{align*}
    L&=L(\lambda)_{\lambda-\alpha}\otimes v_{\mu}+\sum\limits_{\beta\in Q_{+}(r-1)\backslash \{0\}}L(\lambda)_{\lambda-\beta}\otimes L(\mu)_{\mu-\alpha+\beta}+v_{\lambda}\otimes L(\mu)_{\mu-\alpha}\\
    &\subset L(\lambda)_{\lambda-\alpha}\otimes v_{\mu}+L'+qL.
\end{align*}
For any $\widetilde{f}_{i_{1}l_{1}}\cdots \widetilde{f}_{i_{r}l_{r}}v_{\lambda}\in B(\lambda)_{\lambda-\alpha}$, from Lemma \ref{a-g}(d), we have 
$$\widetilde{f}_{i_{1}l_{1}}\cdots \widetilde{f}_{i_{r}l_{r}}v_{\lambda}\otimes v_{\mu}\equiv \widetilde{f}_{i_{1}l_{1}}\cdots \widetilde{f}_{i_{r}l_{r}}(v_{\lambda}\otimes v_{\mu})\mod{qL(\lambda)\otimes L(\mu)}.$$
Hence $L\subset L'+qL$. By Nakayama's Lemma, we can derive that $L=L'$.

\end{proof}

For $\lambda,\mu\in P^+$, define a $U_{q}^-(\g)$-module homomorphism
\begin{align*}
S_{\lambda,\mu}: &V(\lambda)\otimes V(\mu)\rightarrow V(\lambda),\quad u\otimes v_\mu\mapsto u,\\
&V(\lambda)\otimes\sum_{(i,l)\in I_{\infty}}b_{il}V(\mu) \longmapsto 0.
\end{align*}

\vskip 2mm

\begin{lemma}\label{Slm}
{
Let $\lambda,\mu\in P^+$.
\begin{enumerate}
	\item[{\rm (a)}] $S_{\lambda,\mu}(L(\lambda)\otimes L(\mu))=L(\lambda)$.
	\item[{\rm (b)}] For any $\alpha\in R_+(r-1)$ and $w\in{(L(\lambda)\otimes L(\mu))}_{\lambda+\mu-\alpha}$, we have
	\begin{align*}
	S_{\lambda,\mu}\circ\widetilde{f}_{il}(w)\equiv\widetilde{f}_{il}\circ S_{\lambda,\mu}(w)\!\!\!\! \mod{qL(\lambda)}.
	\end{align*}
\end{enumerate}	
}
\end{lemma}

\begin{proof}
(a) is obvious. For (b),  assume that
\begin{align*}
w=u\otimes u'=\mathtt b_{i,\mathbf c}u_{\mathbf c}\otimes \mathtt b_{i,\mathbf c'}u_{\mathbf c'},
\end{align*}	
where $u_{\mathbf c}\in L(\lambda)\cap\mathcal{K}_{i}$, $u_{\mathbf c'}\in L(\mu)\cap \mathcal{K}_{i}$.

Let $L$ be the $\A_0$-submodule of $V(\lambda)\otimes V(\mu)$ generated by $\mathtt b_{i,\mathbf c_{1}}u_{\mathbf{c}}\otimes \mathtt b_{i,\mathbf c_{2}}u_{\mathbf{c}'}$ for all $\mathbf c_{1}$ and  $\mathbf c_{2}$. Thus $L\subset L(\lambda)\otimes L(\mu)$. By the tensor product rule, we have
\begin{align*}
\widetilde{f}_{il}(w)=\widetilde{f}_{il}(u\otimes u')=\begin{cases}
\widetilde{f}_{il}u\otimes u', & \text{if}\ \varphi_i(u)>0,\\
(-1)^{p(li)p(u)}u\otimes \widetilde{f}_{il}u', & \text{if}\ \varphi_i(u)=0.
\end{cases}
\end{align*}

Observe that $S_{\lambda,\mu}\circ\widetilde{f}_{il}(w)\equiv \widetilde{f}_{il}\circ S_{\lambda,\mu}(w)\equiv 0\mod{qL(\lambda)}$ unless $\mathbf{c}'=0$ and $u'\in L(\mu)_{\mu}$. If $\mathbf{c}'=0$ and $u'=v_{\mu}$ then we have 
\begin{align*}
\widetilde{f}_{il}(w)=\widetilde{f}_{il}(u\otimes v_{\mu})\equiv\begin{cases}
\widetilde{f}_{il}u\otimes v_{\mu}\mod{qL}, & \text{if}\ \varphi_i(u)>0,\\
(-1)^{p(li)p(u)}u\otimes \widetilde{f}_{il}v_{\mu}\mod{qL}, & \text{if}\ \varphi_i(u)=0.
\end{cases}
\end{align*}

 If $\varphi_i(u)>0$, then $S_{\lambda,\mu}\circ\widetilde{f}_{il}(u\otimes v_{\mu})\equiv \widetilde{f}_{il}\circ S_{\lambda,\mu}(u\otimes v_{\mu})\mod{qL(\lambda)}$.

 If $\varphi_i(u)=0$, then $\widetilde{f}_{il}(u)=0$. Hence $S_{\lambda,\mu}\circ \widetilde{f}_{il}(w)\equiv 0\equiv \widetilde{f}_{il}\circ S_{\lambda,\mu}(w)\mod{qL(\lambda)}$.

\end{proof}

\vskip 3mm

\begin{lemma}\label{fP and eP}
{
Let $\alpha\in R_+$ and $P \in {U^-_q(\g)}_{-\alpha}$. For any $\lambda\gg 0$, we have
\begin{align*}
&(\widetilde{f}_{il}P)v_\lambda\equiv \widetilde{f}_{il}(Pv_\lambda)\!\!\!\! \mod{qL(\lambda)},\\
&(\widetilde{e}_{il}P)v_\lambda\equiv \widetilde{e}_{il}(Pv_\lambda)\!\!\!\! \mod{qL(\lambda)}.
\end{align*}
}
\end{lemma}

\begin{proof}
We may assume that $P=\mathtt b_{i,\mathbf c}T$ and $e'_{ik}T=0$ for any $k>0$.

 Since $\lambda\gg 0$,  for any $l>0$, we have 
 $$a_{il}(T v_{\lambda})=\frac{q^{ld_{i}\langle h_{i},\lambda+\text{wt}(T)+l\alpha_{i}\rangle}e_{il}''(T)}{1-(-1)^{p(li)}q_{i}^{2l}}v_{\lambda}\in qL(\lambda). $$

 \begin{enumerate}
\item If $i\in I^{\text{re}}$, we have 
\begin{align*}
(\widetilde{f}_{i}P)v_\lambda&=(\widetilde{f}_{i}( b_{i}^{(m)}T))v_\lambda
= b_{i}^{(m+1)}Tv_\lambda
\equiv\widetilde{f}_{i}(Pv_\lambda)\!\!\!\! \mod{qL(\lambda)},\\
 (\widetilde{e}_{i}P)v_\lambda&=(\widetilde{e}_{i} (b_{i}^{(m)}T))v_\lambda
 =( b_{i}^{(m-1)}T)v_\lambda
 \equiv\widetilde{e}_{i}(Pv_\lambda)\!\!\!\! \mod{qL(\lambda)}.
\end{align*}

 	\item If $i\in I^{\text{im}}\backslash I^{\text{iso}}$, we have
 	\begin{align*}
 	(\widetilde{f}_{il}P)v_\lambda=(\widetilde{f}_{il}( b_{i,\mathbf c}T))v_\lambda	=( b_{il}( b_{i,\mathbf c}T))v_\lambda= b_{il}( b_{i,\mathbf c}(Tv_\lambda)).
 	\end{align*}
 Since $a_{il}(Tv_\lambda)\equiv 0 \!\!\mod{qL(\lambda)}$ for any $l>0$, we have
 \begin{align*}
 & b_{il}( b_{i,\mathbf c}(Tv_\lambda))\equiv\widetilde{f}_{il}( b_{i,\mathbf c}Tv_\lambda)\equiv\widetilde{f}_{il}(Pv_\lambda)\!\!\!\!\mod{qL(\lambda)},
\end{align*}
and
\begin{align*}
(\widetilde{e}_{il}P)v_\lambda=&(\widetilde{e}_{il}( b_{i,\mathbf c}T))v_\lambda=( b_{i,\mathbf c\setminus c_1}T)v_\lambda\\
=& b_{i,\mathbf c\setminus c_1}(Tv_\lambda)\equiv\widetilde{e}_{il}( b_{i,\mathbf c}Tv_\lambda)
\equiv \widetilde{e}_{il}(Pv_\lambda)\!\!\!\! \mod{q L(\lambda)}.
\end{align*}
\item If $i\in I^{\text{iso}}$, we have
\begin{align*}
(\widetilde{f}_{il}P)v_\lambda&=(\widetilde{f}_{il}( b_{i,\mathbf c}T))v_\lambda
=\frac{1} {\mathbf c_l+1} \, ( b_{il}b_{i,\mathbf c}T)v_\lambda\\
&=\frac{1}{\mathbf c_l+1}\, b_{il}( b_{i,\mathbf c}Tv_\lambda)\equiv\widetilde{f}_{il}( b_{i,\mathbf c}Tv_\lambda)\equiv\widetilde{f}_{il}(Pv_\lambda)\!\!\!\! \mod{qL(\lambda)},\\
 (\widetilde{e}_{il}P)v_\lambda&=(\widetilde{e}_{il} b_{i,\mathbf c}T)v_\lambda
 =(-1)^{t}{\mathbf c}_l \,( b_{i,\mathbf c\setminus \{l\}}T)v_\lambda
 =(-1)^{t}{\mathbf c}_l ( b_{i,\mathbf c\setminus \{l\}}(Tv_\lambda))\\
&\equiv\widetilde{e}_{il}( b_{i,\mathbf c}Tv_\lambda)\equiv\widetilde{e}_{il}(Pv_\lambda)\!\!\!\! \mod{qL(\lambda)}.
\end{align*}
 \end{enumerate}
\end{proof}

\vskip 3mm

\begin{proposition}  \label{proof of Ir} $({\bf I}(r))$
{
For $\lambda\in P^+$, $\alpha\in R_+(r-1)$ and $P\in {L(\infty)}_{-\alpha}$,
we have
\begin{align*}
(\widetilde{f}_{il}P)v_\lambda\equiv \widetilde{f}_{il}(Pv_\lambda) \!\!\!\!\mod{qL(\lambda)}.
\end{align*}

In particular, we have
\begin{align*}
(\widetilde{f}_{i_1l_1}\cdots\widetilde{f}_{i_rl_r}\mathbf 1)	v_\lambda\equiv \widetilde{f}_{i_1l_1}\cdots\widetilde{f}_{i_rl_r} v_\lambda\!\!\!\! \mod{qL(\lambda)}.
	\end{align*}
}
\end{proposition}

\begin{proof}
Take $\mu \gg 0$ such that $\lambda+\mu\gg 0$. By Lemma \ref{fP and eP}, we have
\begin{align*}
(\widetilde{f}_{il}P)v_{\lambda+\mu}\equiv\widetilde{f}_{il}(Pv_{\lambda+\mu})\!\!\!\! \mod{qL(\lambda+\mu)}.
\end{align*}

\vskip 2mm

By Proposition \ref{prop:Cr},
we obtain
\begin{align}\label{fP}
(\widetilde{f}_{il}P)(v_\lambda\otimes v_\mu)\equiv\widetilde{f}_{il}(P(v_\lambda\otimes v_\mu))\!\!\!\! \mod{q(L(\lambda)\otimes L(\mu))}.
\end{align}

\vskip 2mm

On the other hand, by ${\bf H}(r-1)$ and ${\bf C}(r-1)$, we have
\begin{align*}
P(v_\lambda\otimes v_\mu)=\Phi_{\lambda,\mu}(Pv_{\lambda+\mu})\in L(\lambda)\otimes L(\mu).
\end{align*}

Applying $S_{\lambda,\mu}$ to \eqref{fP}, then Lemma \ref{Slm} yields
\begin{align*}
(\widetilde{f}_{il}P)v_\lambda\equiv \widetilde{f}_{il}(Pv_\lambda) \!\!\!\!\mod{qL(\lambda)}.
\end{align*}
\end{proof}

\vskip 2mm

\begin{proposition} \label{proof of Hr} $({\bf H}(r))$
{
For any $\lambda\in P^+$ and $\alpha\in R_+(r)$, we have
\begin{align*}
\pi_\lambda({L(\infty)}_{-\alpha})={L(\lambda)}_{\lambda-\alpha}.
\end{align*}	
}
\end{proposition}
\begin{proof}
    By proposition \ref{proof of Ir}, we have 
    $$\pi_{\lambda}(\widetilde{f}_{i_1l_1}\cdots\widetilde{f}_{i_rl_r}\mathbf 1	)\equiv \widetilde{f}_{i_1l_1}\cdots\widetilde{f}_{i_rl_r} 	v_\lambda\mod{qL(\lambda)}$$
   Hence, $\pi_{\lambda}(L(\infty)_{-\alpha})\subset L(\lambda)_{\lambda-\alpha}$ and $L(\lambda)_{\lambda-\alpha}\subset \pi_{\lambda}(L(\infty)_{-\alpha})+qL(\lambda)_{\lambda-\alpha}$.
   We get the  result from Nakayama's Lemma.
\end{proof}

\begin{corollary}\label{pi}
{
Consider the $\Q$-linear map
\begin{align*}
\overline{\pi}_\lambda: {L(\infty)}_{-\alpha}/q{L(\infty)}_{-\alpha}\longrightarrow {L(\lambda)}_{\lambda-\alpha}/q{L(\lambda)}_{\lambda-\alpha}.
\end{align*}

\begin{enumerate}
	\item[{\rm (a)}] For any $\beta\in R_+(r-1)$ and $b\in{B(\infty)}_{-\beta}$, we have
	\begin{align*}
	\overline{\pi}_\lambda(\widetilde{f}_{il}b)=\widetilde{f}_{il}(\overline{\pi}_\lambda(b)).
	\end{align*}
\item[{\rm (b)}] For any $\alpha\in  R_+(r)$ and $\lambda\in P^+$, we have
\begin{align*}
\overline{\pi}_\lambda({B(\infty)}_{-\alpha})\backslash \{0\}={B(\lambda)}_{\lambda-\alpha}.
\end{align*}
\item[{\rm (c)}] For any $\alpha\in R_+(r)$ and $\lambda\gg 0$, the map $\pi_\lambda$ induces the isomorphisms
\begin{align*}
{L(\infty)}_{-\alpha}\stackrel{\sim}{\rightarrow}{L(\lambda)}_{\lambda-\alpha},\quad  {B(\infty)}_{-\alpha}
\stackrel{\sim}{\rightarrow}{B(\lambda)}_{\lambda-\alpha}.
\end{align*}
\end{enumerate}
}
\end{corollary}
\begin{proof}
    (a) can be derived from Proposition \ref{proof of Ir}. (b) follows from (a).
    When $\lambda\gg0$, $\pi_{\lambda}$ is an isomorphism, then we obtain (c)
\end{proof}

\vskip 3mm

Fix $\lambda\in  P^+$, $i\in I$, $l_1,\cdots,l_r>0$ and $\alpha=\sum_{j=1}^{r}l_j\alpha_{i_j}$. Take a finite set $T$ containing $\lambda$ and the fundamental weights $\Lambda_{i_1},\cdots,\Lambda_{i_r}$.

\begin{enumerate}
	\item[i)] Since $T$ is a finite set, we can take a sufficient large $N_1\geq 0$ such that
	\begin{align*}
	\widetilde{e}_{il}{L(\tau)}_{\tau-\alpha}\subset q^{-N_1}L(\tau)\ \text{for all}\ \tau\in T\ \text{and}\ l>0.
	\end{align*}
    \item[ii)] Choose $N_2\geq 0$ such that
 $
    \widetilde{e}_{il}{L(\infty)}_{-\alpha}\subset q^{-N_2}L(\infty)\ \text{for all}\ l>0$.

  \vskip 3mm

Then for any $\mu\gg 0$, Lemma \ref{fP and eP} and Proposition \ref{proof of Hr}
yield
\begin{align*}
\widetilde{e}_{il}{L(\mu)}_{\mu-\alpha}&=\widetilde{e}_{il}({L(\infty)}_{-\alpha}v_\mu)\subset (\widetilde{e}_{il}{L(\infty)}_{-\alpha})v_\mu+q{L(\mu)}_{\mu-\alpha}\\
&\subset q^{-N_2}{L(\infty)}_{-\alpha}v_\mu+q{L(\mu)}_{\mu-\alpha}\subset q^{-N_2}L(\mu).
\end{align*}
\end{enumerate}

\vskip 2mm

Therefore, for any $\alpha\in R_+(r)$, there exists $N\geq 0$ such that
\begin{equation}\label{eL}
\begin{aligned}
&\widetilde{e}_{il}{L(\mu)}_{\mu-\alpha}\subset q^{-N}L(\mu)\ \text{for all}\ \mu\gg 0,\\
&\widetilde{e}_{il}{L(\tau)}_{\tau-\alpha}\subset q^{-N}L(\tau)\ \text{for all}\ \tau\in T,\\
&\widetilde{e}_{il}{L(\infty)}_{-\alpha}\subset q^{-N}L(\infty).	
	\end{aligned}
\end{equation}

\vskip 3mm

\begin{lemma}\label{eLL}
{
For any $\alpha\in R_+$, let $N\geq 0$ be a non-negative integer satisfying \eqref{eL}. For any $\mu\gg 0$ and $\tau\in T$, we have
\begin{align*}
\widetilde{e}_{il}{(L(\tau)\otimes L(\mu))}_{\tau+\mu-\alpha}\subset q^{-N}(L(\tau)\otimes L(\mu)).
\end{align*}
}
\end{lemma}

\begin{proof}
Let $u\in {L(\tau)}_{\tau-\beta}$ and $v\in {L(\mu)}_{\mu-\gamma}$	with $\alpha=\beta+\gamma$.
We will show that $\widetilde{e}_{il}(u\otimes v)\in q^{-N}(L(\tau)\otimes L(\mu))$. 
If $\beta\neq 0$ and $\gamma\neq 0$, by Lemma \ref{a-g}(a) and \eqref{eL}, we have $\widetilde{e}_{il}(u\otimes v)\subset q^{-N}L(\tau)\otimes L(\mu)$.

If $\beta=0$, then $\gamma=\alpha$, we may assume that $u=v_\tau$.
Let $v=\sum_{\mathbf c\in \mathcal C_i}\mathtt b_{i,\mathbf c}v_{\mathbf c}$ be the $i$-string decomposition of $v$. By \eqref{eL}, we have
\begin{align*}
\widetilde{e}_{il}v=\begin{cases}
\sum_{m\ge 1}b_{i}^{(m-1)}v_{m}\in q^{-N}L(\mu), &\text{if}\ i\in I^{\text{re}},\\
\sum_{\mathbf c\neq\mathbf 0,c_{1}=l}b_{i,\mathbf c\setminus c_1}v_{\mathbf c}\in q^{-N}L(\mu), & \text{if}\ i\in I^{\text{im}}\backslash I^{\text{iso}},\\
\sum_{\mathbf c\neq\mathbf 0} (-1)^{p(li)p(|\mathbf{c}_{<k}|i)} \mathbf c_l \,
 b_{i,\mathbf c\setminus \{l\}}v_{\mathbf c}\in q^{-N}L(\mu), & \text{if}\ i\in I^{\text{iso}}.
\end{cases}
\end{align*}
Hence by Lemma \ref{uc in L}, we obtain
\begin{align*}
v_{\mathbf c}\in q^{-N}L(\mu)\ \text{for any}\ \mathbf c\neq\mathbf 0.
\end{align*}

\vskip 2mm

Let $L$ be the $\A_0$-submodule of $L(\tau)\otimes L(\mu)$ generated by $\mathtt b_{i,\mathbf c_1}v_\tau\otimes\mathtt b_{i,\mathbf c_2}v_{\mathbf c}$ for $\mathbf c_1,\mathbf c_2,\mathbf c \neq \mathbf 0$. Then $\widetilde{e}_{il}L\subset L$. It follows that
\begin{align*}
\widetilde{e}_{il}(v_\tau\otimes v)=\sum_{\mathbf c\neq\mathbf 0}\widetilde{e}_{il}(v_{\tau}\otimes\mathtt b_{i,\mathbf c}v_{\mathbf c})\in L\subset q^{-N}(L(\tau)\otimes L(\mu)).
\end{align*}
If
 $\beta=\alpha$, $\gamma=0$, we can derive that $\widetilde{e}_{il}(u\otimes v_{\mu})\in q^{-N}(L(\tau)\otimes L(\mu))$ by a similar argument.
	\end{proof}

\begin{lemma}\label{lemma for Ar and Hr}
{
Let $\alpha\in R_+(r)$ and let $N>0$ be the positive integer satisfying \eqref{eL}. Then we have
\begin{enumerate}
	\item[{\rm (a)}] $\widetilde{e}_{il}{L(\mu)}_{\mu-\alpha}\subset q^{1-N}L(\mu)$ for all $\mu\gg 0$,
	\item[{\rm (b)}] $\widetilde{e}_{il}{L(\tau)}_{\tau-\alpha}\subset q^{1-N}L(\tau)$ for all $\tau\in T$,
	\item[{\rm (c)}] $\widetilde{e}_{il}{L(\infty)}_{-\alpha}\subset q^{1-N}L(\infty)$.	
\end{enumerate}
}
\end{lemma}
\begin{proof}
	(a)
Let $u=\widetilde{f}_{i_1l_1}\cdots\widetilde{f}_{i_t l_t}v_\mu\in {L(\mu)}_{\mu-\alpha}$. Suppose $i_1=i_2=\cdots=i_t$. If $i=i_1$, then
\begin{align*}
u=\mathtt b_{i,\mathbf c}v_\mu,\ \mathbf c=(l_1,\cdots,l_t),
\end{align*}
and
\begin{align*}
\widetilde{e}_{il}u=\widetilde{e}_{il}(\mathtt b_{i,\mathbf c}v_\mu)=\begin{cases}
\mathtt b_{i,\mathbf c\setminus c_1}v_\mu, & \text{if}\ i\notin I^{\text{iso}},\ c_1=l,\\
\mathbf c_l \, \mathtt b_{i,\mathbf c\setminus \{l\}}v_\mu, & \text{if}\ i\in I^{\text{iso}},\ l\in\mathbf c,\\
0, & \text{otherwise}.
\end{cases}
\end{align*}
Therefore, we have $\widetilde{e}_{il}u\in L(\mu)$. 

\vskip 2mm

If $i\neq i_1$, then $\widetilde{e}_{il}u=0$.
 We may assume that there exists $s$ with $1\leq s<t$ such that
$i_s\neq i_{s+1}=\cdots=i_t$.
Suppose $\mu \gg 0$ and set $\lambda_0=\Lambda_{i_s}$. Then $\mu':=\mu-\lambda_0\gg 0$.
Set
\begin{align*}
w:=\widetilde{f}_{i_1l_1}\cdots\widetilde{f}_{i_t l_t} (v_{\lambda_0}\otimes v_{\mu'}).
\end{align*}
By Lemma \ref{ffvvbb}, we have
\begin{align*}
w\equiv v\otimes v'\!\!\!\!\! \mod{qL(\lambda_0)\otimes L(\mu')}
\end{align*}
for some $v\in {L(\lambda_0)}_{\lambda_0-\beta}$, $v'\in {L(\mu')}_{\mu'-\gamma}$, $\alpha=\beta+\gamma$ and $\beta,\gamma\in R_+(r-1)$.

\vskip 2mm

Then Lemma \ref{a-g} (a) and Lemma \ref{eLL} imply
\begin{align*}
\widetilde{e}_{il}w&\in L(\lambda_0)\otimes L(\mu')+q \, \widetilde{e}_{il}{(L(\lambda_0)\otimes L(\mu'))}_{\lambda_0+\mu'-\alpha}\\
&\subset L(\lambda_0)\otimes L(\mu')+q^{1-N}L(\lambda_0)\otimes L(\mu')=q^{1-N}L(\lambda_0)\otimes L(\mu').
\end{align*}

Thus we have
\begin{align*}
\widetilde{e}_{il}w\in q^{1-N}{(L(\lambda_0)\otimes L(\mu'))}_{\lambda_0+\mu'-\alpha+l\alpha_i}=q^{1-N}{(L(\lambda_0)\otimes L(\mu'))}_{\mu-\alpha+l\alpha_i}.
\end{align*}
Applying $\Psi_{\lambda_0,\mu'}$ to ${\bf D}(r-1)$, we have
\begin{align*}
\widetilde{e}_{il}u=\widetilde{e}_{il}\widetilde{f}_{i_1l_1}\cdots\widetilde{f}_{i_t l_t}v_\mu\in q^{1-N}L(\mu).
\end{align*}
(b)
Let $\tau\in T$ and set
$
u=\widetilde{f}_{i_1l_1}\cdots\widetilde{f}_{i_tl_t}v_\tau\in {L(\tau)}_{\tau-\alpha}$.
If $u\in qL(\tau)$, our assertion follows from \eqref{eL}.

If $u\notin qL(\tau)$, for any $\mu\in P^+$, Lemma \ref{a-g} (f) gives
\begin{align}\label{ffvv}
\widetilde{f}_{i_1l_1}\cdots\widetilde{f}_{i_tl_t}(v_\tau\otimes v_\mu)\equiv u\otimes v_\mu\!\!\!\! \mod{qL(\tau)\otimes L(\mu)}.
\end{align}

If $\mu\gg 0$, (a) implies
\begin{align*}
\widetilde{e}_{il}\widetilde{f}_{i_1l_1}\cdots\widetilde{f}_{i_tl_t}v_{\tau+\mu}\in q^{1-N}L(\tau+\mu).
\end{align*}

Applying $\Phi_{\tau,\mu}$ and $\mathbf{B}(r-1)$, we obtain
\begin{align*}
\widetilde{e}_{il}\widetilde{f}_{i_1l_1}\cdots\widetilde{f}_{i_tl_t}(v_\tau\otimes v_\mu)\in q^{1-N}L(\tau)\otimes L(\mu).
\end{align*}
By \eqref{ffvv} and Lemma \ref{a-g}, we have
\begin{align}\label{evv}
\widetilde{e}_{il}(v\otimes v_\mu)\in q^{1-N}L(\tau)\otimes L(\mu)+q\widetilde{e}_{il}(L(\tau)\otimes L(\mu))\subset q^{1-N}L(\tau)\otimes L(\mu).
\end{align}

Let $u=\sum_{\mathbf c\in\mathcal C_i}\mathtt b_{i,\mathbf c}u_{\mathbf c}$ be the $i$-string decomposition of $u$. By \eqref{eL}, we have $\widetilde{e}_{il}u\in q^{-N}L(\tau)$.

Recall
\begin{align*}
\widetilde{e}_{il}u=
\begin{cases}
\sum_{\mathbf c\in\mathcal C_i}	\mathtt b_{i,\mathbf c\setminus c_1}u_{\mathbf c}, & \text{if}\ i\notin I^{\text{iso}},\ c_1=l,\\
\sum_{\mathbf c\in\mathcal C_i}	\mathbf c_l \, \mathtt b_{i,\mathbf c\setminus \{l\}}u_{\mathbf c}, & \text{if}\ i\in I^{\text{iso}},\ l\in\mathbf c,\\
	0 & \text{otherwise}.
\end{cases}
\end{align*}
By Lemma \ref{uc in L}, we have $u_{\mathbf c}\in q^{-N}L(\tau)$. Let $L$ be the $\A_0$-submodule of $V(\tau)\otimes V(\mu)$ generated by $\mathtt b_{i,\mathbf c}u_{\mathbf c}\otimes \mathtt b_{i,\mathbf c'}v_{\mu}$ $(c_1=l\ \text{or}\ l\in\mathbf c)$. Then we have $L\subset q^{-N}L(\tau)\otimes L(\mu)$.

\vskip 2mm

The tensor product rule gives
\begin{align*}
\widetilde{e}_{il}(u\otimes v_\mu)\equiv \widetilde{e}_{il}u\otimes v_\mu\!\!\!\! \mod{qL}.
\end{align*}
By \eqref{evv}, we have
\begin{align*}
\widetilde{e}_{il}u\otimes v_\mu\equiv\widetilde{e}_{il}(u\otimes v_\mu)\in q^{1-N}L(\tau)\otimes L(\mu).
\end{align*}
Hence $\widetilde{e}_{il}u\in q^{1-N}L(\tau)$.
\vskip 2mm
(c)
Let $S\in {L(\infty)}_{-\alpha}$ and take $\mu\gg 0$. By Lemma \ref{fP and eP}, we have
$
(\widetilde{e}_{il}S)v_\mu\equiv \widetilde{e}_{il}(Sv_\mu)\!\!\! \mod{qL(\mu)}$.
Thus Proposition \ref{proof of Hr}
implies
\begin{align*}
(\widetilde{e}_{il}S)v_\mu=\widetilde{e}_{il}(Sv_\mu)\in\widetilde{e}_{il}{L(\mu)}_{\mu-\alpha}\subset q^{1-N}L(\mu).
\end{align*}

\noindent
Hence by Corollary \ref{pi} (c), we have
\begin{align*}
\widetilde{e}_{il}S\in q^{1-N}L(\infty).
\end{align*}
	\end{proof}
	
	\vskip 2mm


\begin{proof}
If $b\in {B(\infty)}_{-\alpha}$, then there exist $(i,l)\in I^{\infty}$ and $b'\in{B(\infty)}_{-\alpha+l\alpha_i}$ such that $b=\widetilde{f}_{il}b'$. By $\mathbf{G}(r-1)$, the set ${B(\infty)}_{-\alpha+l\alpha_i}$ forms a $\Q$-basis of ${L(\infty)}_{-\alpha+l\alpha_i}/q{L(\infty)}_{-\alpha+l\alpha_i}$. Then we have $b'\neq 0$. Hence $b\neq 0$.
 	\end{proof}
	
\vskip 2mm

\begin{lemma}\label{pie=epi}
{
Let $\alpha\in R_+(r)$, $(i,l)\in I_{\infty}$, $\lambda\gg 0$ and $b\in{B(\infty)}_{-\alpha}$.
Then we have
\begin{equation*}
\overline{\pi}_\lambda(\widetilde{e}_{il}b)=\widetilde{e}_{il}\overline{\pi}_\lambda(b).
\end{equation*}
}
\end{lemma}

\begin{proof}

The assertion follows directly from Lemma \ref{fP and eP}.
	\end{proof}
	
\vskip 2mm

\begin{corollary}\label{fLLeLL} 
{
Let $\lambda,\mu\in P^+$ and $\alpha,\beta\in R_+(r)$.
\begin{enumerate}
	\item[{\rm (a)}] For the $i$-string decomposition $u=\sum_{\mathbf c\in\mathcal C_i}\mathtt b_{i,\mathbf c}u_{\mathbf c}\in {L(\lambda)}_{\lambda-\alpha}$, we have $u_{\mathbf c}\in L(\lambda)$
	for all $\mathbf{c} \in \mathcal{C}_{i}$.
	\item[{\rm (b)}] For any $(i,l)\in I^{\infty}$, we have
	\begin{align*}
	&\widetilde{f}_{il}({L(\lambda)}_{\lambda-\alpha}\otimes {L(\mu)}_{\mu-\beta})\subset L(\lambda)\otimes L(\mu),\\
	&\widetilde{e}_{il}({L(\lambda)}_{\lambda-\alpha}\otimes {L(\mu)}_{\mu-\beta})\subset L(\lambda)\otimes L(\mu).
	\end{align*}
\end{enumerate}
}
\end{corollary}

\begin{proof}
By applying Lemma \ref{lemma for Ar and Hr} repeatedly, we conclude that $N=0$ and we have 
$$\widetilde{e}_{il}L(\infty)_{-\alpha}\subset L(\infty),\ \ \ \widetilde{e}_{il}L(\lambda)_{\lambda-\alpha}\subset L(\lambda)$$
for all $\lambda\in P^{+}$ and $\alpha\in \Q_{+}(r)$.
Since Lemma \ref{a-g} depends only on ${\bf A}(r-1)$, the
corollary follows from the proof of Lemma \ref{uc in L}.
	\end{proof}

\vskip 2mm

\begin{lemma}\label{euv}
{
Let $\lambda,\mu\in P^+$ and $\alpha\in R_+(r)$. For any $u\in{L(\lambda)}_{\lambda-\alpha}$, we have
\begin{align*}
\widetilde{e}_{il}(u\otimes v_\mu)\equiv \widetilde{e}_{il}u\otimes v_\mu\!\!\!\! \mod{q(L(\lambda)\otimes L(\mu))}.	
	\end{align*}
}
\end{lemma}

\begin{proof}
Assume that $u=b_{i,\mathbf{c}}u_{\mathbf{c}}\neq 0$ with $a_{il}u_{\mathbf{c}}=0$ and $u_{\mathbf{c}}\in L(\lambda)_{\lambda-\alpha+|\mathbf{c}|\alpha_{i}}$. Let $L$ be the $\A_{0}$-submodule generated by $b_{i,\mathbf{c}_{1}}u_{\mathbf{c}}\otimes b_{i,\mathbf{c}_{2}}v_{\mu}$. Then the tensor product rule geives
$$\widetilde{e}_{il}(b_{i,\mathbf{c}}u_{\mathbf{c}}\otimes v_{\mu})\equiv b_{i,\mathbf{c}\backslash l}u_{\mathbf{c}}\otimes v_{\mu}\mod{qL}.$$
Our assertion follows from $L\subset L(\lambda)\otimes L(\mu)$
	\end{proof}
	
\vskip 2mm

\begin{proposition}  \label{proof of Kr} $({\bf K}(r))$
{
	Let $\lambda\in P^+$ and $\alpha\in R_+(r)$. If $b\in{B(\infty)}_{-\alpha}$ and $\overline{\pi}_\lambda(b)\neq 0$, then we have
	\begin{align*}
		\widetilde{e}_{il}\overline{\pi}_\lambda(b)=\overline{\pi}_\lambda(\widetilde{e}_{il}b).
		\end{align*}
		
}
\end{proposition}

\begin{proof}
Let
\begin{align*}
P&=\widetilde{f}_{i_1l_1}\cdots\widetilde{f}_{i_tl_t}\mathbf 1\in {L(\infty)}_{-\alpha},\\
b&=P +q{L(\infty)}_{-\alpha}\in {B(\infty)}_{-\alpha},\\
u&=\widetilde{f}_{i_1l_1}\cdots\widetilde{f}_{i_tl_t}v_\lambda.
\end{align*}
By Proposition  \ref{proof of Ir},
we have
$
u\equiv Pv_{\lambda}\mod{qL(\lambda)}.
$

\noindent
Since $\overline{\pi}_\lambda(b)\neq 0$, we have  $u\notin qL(\lambda)$.  For any $\mu\in P^+$, by Lemma \ref{a-g} (e), we have
\begin{align*}
\widetilde{f}_{i_1l_1}\cdots\widetilde{f}_{i_tl_t}(v_\lambda\otimes v_\mu)\equiv \widetilde{f}_{i_1l_1}\cdots\widetilde{f}_{i_tl_t}v_\lambda\otimes v_\mu\equiv u\otimes v_\mu\!\!\!\! \mod{q(L(\lambda)\otimes L(\mu))}.
\end{align*}

\noindent
Hence by Lemma \ref{euv}, we have
\begin{align}\label{effvv1}
\widetilde{e}_{il}(\widetilde{f}_{i_1l_1}\cdots\widetilde{f}_{i_tl_t}(v_\lambda\otimes v_\mu))\equiv\widetilde{e}_{il}(u\otimes v_\mu)\equiv\widetilde{e}_{il}u\otimes v_\mu\!\!\!\! \mod{q(L(\lambda)\otimes L(\mu))}.	
	\end{align}
	\vskip 3mm
	
On the other hand, for $\mu\gg 0$, by Lemma \ref{pie=epi}, we have
\begin{align*}
\widetilde{e}_{il}(\widetilde{f}_{i_1l_1}\cdots\widetilde{f}_{i_tl_t}v_{\lambda+\mu})\equiv \widetilde{e}_{il}(Pv_{\lambda+\mu})\equiv (\widetilde{e}_{il}P)v_{\lambda+\mu}\!\!\!\! \mod{qL(\lambda+\mu)}.
\end{align*}
Applying $\Phi_{\lambda,\mu}$, from  Proposition \ref{prop:Cr}, we obtain 
\begin{align}\label{effvv2}
\widetilde{e}_{il}(\widetilde{f}_{i_1l_1}\cdots\widetilde{f}_{i_tl_t}(v_\lambda\otimes v_\mu))\equiv(\widetilde{e}_{il}P)(v_\lambda\otimes v_\mu)\!\!\!\! \mod{q(L(\lambda)\otimes L(\mu))}.
\end{align}
Then \eqref{effvv1} and \eqref{effvv2} yield
\begin{align*}
\widetilde{e}_{il}u\otimes v_\mu\equiv (\widetilde{e}_{il}P)(v_\lambda\otimes v_\mu)\!\!\!\! \mod{q(L(\lambda)\otimes L(\mu))}.	
	\end{align*}
Applying $S_{\lambda,\mu}$, we have
\begin{align*}
\widetilde{e}_{il}u\equiv (\widetilde{e}_{il}P)v_\lambda\!\!\!\! \mod{qL(\lambda)}.
\end{align*}
Hence $\widetilde{e}_{il}\overline{\pi}_\lambda(b)= \overline{\pi}_{\lambda}(\widetilde{e}_{il}b).$
	\end{proof}
	
\vskip 2mm

\begin{proposition}  \label{proof of Er} $({\bf E}(r))$ \
{	For every $\alpha\in R_+(r)$, we have
	\begin{align*}
	\widetilde{e}_{il}{L(\infty)}_{-\alpha}\subset L(\infty),\quad	\widetilde{e}_{il}{B(\infty)}_{-\alpha}\subset B(\infty)\cup \{0\}.
	\end{align*}
}
\end{proposition}

\begin{proof}
Applying Lemma \ref{lemma for Ar and Hr} (c) repeatedly, the first assertion holds.
For the second assertion, let
	\begin{align*}
		S=\widetilde{f}_{i_1l_1}\cdots\widetilde{f}_{i_tl_t}\mathbf 1\in{L(\infty)}_{-\alpha},\quad b=S+q{L(\infty)}_{-\alpha}\in {B(\infty)}_{-\alpha}.
	\end{align*}
	
	If $i_1=i_2=\cdots=i_t$, our assertion is true as we have seen in the proof of Lemma \ref{lemma for Ar and Hr} (a). Here,  we may assume that there exists $s$ with  $1\leq s<t$ such that $i_s\neq i_{s+1}=\cdots= i_t$. Take $\mu\gg 0$ and set $\lambda_0=\Lambda_{i_s}$, $\lambda=\lambda_0+\mu\gg 0$. Then Lemma \ref{ffvvbb} yields
	\begin{align*}
		S(v_{\lambda_0}\otimes v_\mu)=\widetilde{f}_{i_1l_1}\cdots\widetilde{f}_{i_tl_t}(v_{\lambda_0}\otimes v_\mu)\equiv v\otimes v'\!\!\!\! \mod{q(L(\lambda_0)\otimes L(\mu))}
	\end{align*}
	for some $v\in{L(\lambda_0)}_{\lambda_0-\beta}$, $v'\in{L(\mu)}_{\mu-\gamma}$, $\beta,\gamma\in R_+(r-1)\setminus\{0\}$ and $\alpha=\beta+\gamma$ such that
	\begin{align*}
		v+qL(\lambda_0)\in B(\lambda_0)\cup \{0\},\quad v'+qL(\mu)\in B(\mu)\cup \{0\}.
	\end{align*}
   Therefore we have
	\begin{align*}
		\widetilde{e}_{il}(\widetilde{f}_{i_1l_1}\cdots\widetilde{f}_{i_tl_t}(v_{\lambda_0}\otimes v_\mu))\equiv\widetilde{e}_{il}v\otimes v' \text{or} (-1)^{p(li)p(u)}v\otimes \widetilde{e}_{il}v'\!\!\!\! \mod{q(L(\lambda_0)\otimes L(\mu))}.
	\end{align*}
	By ${\bf A}(r-1)$, we have
	\begin{align*}
		\widetilde{e}_{il}(\widetilde{f}_{i_1l_1}\cdots\widetilde{f}_{i_tl_t}(v_{\lambda_0}\otimes v_\mu))+q(L(\lambda_0)\otimes L(\mu))\in (B(\lambda_0)\otimes B(\mu))\cup \{0\}.
	\end{align*}
	The map $\Psi_{\lambda_0,\mu}$ and $\mathbf{D}(r-1)$ yield
	\begin{align*}
		\widetilde{e}_{il}\overline{\pi}_\lambda(b)=\widetilde{e}_{il}(\widetilde{f}_{i_1l_1}\cdots\widetilde{f}_{i_tl_t}v_\lambda+qL(\lambda))\in B(\lambda)\cup \{0\}.
	\end{align*}
	
	Since $\lambda\gg 0$, Lemma \ref{pie=epi} and Corollary \ref{pi} (c) yield
	\begin{align*}
		\widetilde{e}_{il}b=\widetilde{e}_{il}(\widetilde{f}_{i_1l_1}\cdots\widetilde{f}_{i_tl_t}\mathbf 1+qL(\infty))\in B(\infty)\cup\{0\}.
	\end{align*}
\end{proof}

\vskip 2mm

\begin{proposition}  \label{proof of Ar} $({\bf A}(r))$
{
For any $\lambda\in P^+$ and $\alpha\in R_+(r)$, we have
\begin{align*}
\widetilde{e}_{il}{L(\lambda)}_{\lambda-\alpha}\subset L(\lambda),\quad
\widetilde{e}_{il}{B(\lambda)}_{\lambda-\alpha}\subset B(\lambda)\cup \{0\}.
\end{align*}
}
\end{proposition}

\begin{proof}
	Proposition follows from Lemma \ref{lemma for Ar and Hr}, Proposition  \ref{proof of Kr}, Proposition \ref{proof of Er},
	Corollary \ref{pi} (b) and Proposition \ref{proof of Hr}.	
	\end{proof}

\vskip 2mm

For $(i,l) \in I^{\infty}$, let $u = {\mathtt b}_{i,\mathbf{c}} u_{\mathbf{c}}$ such that $a_{ik} u_{\mathbf{c}} = 0$ for all $k>0$.
Define an operator $Q_{il}: V(\lambda) \rightarrow V(\lambda)$ by

\begin{equation} \label{eq:Qill}
Q_{il}(u) = \begin{cases}
(\mathbf{c}_{l}+1) u & \text{if} \ i\in I^{\text{iso}}, \\
u & \text{otherwise}
\end{cases}
\end{equation}

\vskip 2mm

\begin{lemma}\label{lem:LLA0}
{
	Let $\lambda\in P^{+}$ and $\alpha\in R_{+}(r)$.
	\begin{enumerate}	
		\item[{\rm (a)}] For any $u\in {L(\lambda)}_{\lambda-\alpha+l\alpha_i}$ and $v\in {L(\lambda)}_{\lambda-\alpha}$, we have
		\begin{equation*}
			(\widetilde{f}_{il} Q_{il}(u),v)_K\equiv (u,\widetilde{e}_{il}v)_K\!\!\!\! \mod{q \A_0}.
		\end{equation*}
	
\vskip 2mm
		\item[{\rm (b)}]
		$({L(\lambda)}_{\lambda-\alpha},{L(\lambda)}_{\lambda-\alpha})_K\subset \A_0$.
	\end{enumerate}
}
\end{lemma}

\begin{proof}
	(a)We will use induction on $\text{ht}(\alpha)$. When $\text{ht}(\alpha)=0$, our assertion is trivial.  For $u, v\in L(\lambda)_{\lambda-\alpha}$, Let $u=b_{i,\mathbf{c}}u_{\mathbf{c}}$ and $v=b_{i,\mathbf{c}'}v_{\mathbf{c'}}$, where $u_{\mathbf{c}},v_{\mathbf{c'}}\in \mathcal{K}_{i}, |\mathbf{c}|=n$ and $|\mathbf{c}'|=m$. From  Corollary \ref{fLLeLL}, we konw that $u_{\mathbf{c}},v_{\mathbf{c'}}\in L(\lambda)$.
    When $i\in I^{im}$, by Proposition \ref{2.2} Lemma \label{4.6} and Corollary \ref{fLLeLL} , we have $ -K_{i}^{l}a_{il}v \in(1+q\mathcal{A}_{0})\widetilde{e}_{il}v$. When $i\in I^{re}$, Since the bilinear form is symmetric, we can assume that $n\ge m$. Then by Lemma 1.7 in \cite{FH2} and Proposition \ref{2.2},  we have 
      $
         -[n+1]_{ i}^{-1} K_{i}a_{i}v\in(1+q\mathcal{A}_{0})\widetilde{e}_{i}(v).
      $
 Hence,  $(\widetilde{f}_{il} Q_{il}(u),v)_K\in (1+q\mathcal{A}_{0}) (u,\widetilde{e}_{il}v)_K$. By induction hypothesis, $(\widetilde{f}_{il} Q_{il}(u),v)_K\equiv (u,\widetilde{e}_{il}v)_K\!\!\!\! \mod{q \A_0}.$

 \vskip 2mm
	
	(b) follows directly from (a).
\end{proof}

\vskip 2mm

\begin{lemma}\label{(P,Q)}
{
	Let $\alpha=l_1\alpha_{i_1}+\cdots+l_t\alpha_{i_t}\in R_{+}$, $S,T\in {U^-_q(\g)}_{-\alpha}$ and $m\in \Z$. For any $\lambda\gg 0$, we have
	\begin{equation*}
		(Sv_{\lambda},Tv_{\lambda})_K=\prod_{k=1}^{t}{(1-(-1)^{p(l_{k}i_{k})}q_{i_{k}}^{2l_k})}^{-1} (S,T)_K\!\!\!\!\mod{q^{m}\A_0}.
	\end{equation*}
}
\end{lemma}
\begin{proof}
	We will prove the assertion by induction on $\mathrm{ht}(\alpha)$. If  $\alpha=0$, Our assertion is obvious.
	 Assume that $S=\mathtt b_{il}W$
	for some $W\in {U^{-}_q(\g)}_{-\alpha+l\alpha_i}$.
	Then we have
	\begin{align*}
		(Sv_{\lambda},Tv_{\lambda})_K=&-(Wv_{\lambda},K_{i}^{l}a_{il}(Tv_{\lambda}))_K\\
		=&(Wv_{\lambda},\frac{e'_{il}(T)-K_i^{2l}e''_{il}(T)}{1-(-1)^{p(li)}q^{2l}_i}v_{\lambda})_K\\
		=&\frac{\{(Wv_{\lambda},e'_{il}(T)v_{\lambda})_K-q_i^{2l\langle h_i,\lambda-\alpha\rangle+2 l^2 a_{ii}}(Wv_{\lambda},e_{il}''(T)v_{\lambda})_K\}}{1-(-1)^{p(li)}q_i^{2l}}.
	\end{align*}

	Since $\lambda\gg 0$, then $q_i^{2l\langle h_i,\lambda-\alpha\rangle+2l^2a_{ii}}\in q^{m}\A_0$, which implies $-\frac{q_i^{2l\langle h_i,\lambda-\alpha\rangle+2 l^2 a_{ii}}(Wv_{\lambda},e_{il}''(T)v_{\lambda})_K}{1-(-1)^{p(li)}q_i^{2l}}\in q^{m}\A_{0}$.
By induction, we obtain
	\begin{align*}
		(Sv_{\lambda},Tv_{\lambda})_K&\equiv {(1-(-1)^{p(li)}q_i^{2l})}^{-1}(Wv_{\lambda},e'_{il}(T)v_{\lambda})_K\\
		&\equiv{(1-(-1)^{p(li)}q_i^{2l})}^{-1}(W,e'_{il}T)_K\\
        &\equiv{(1-(-1)^{p(li)}q_i^{2l})}^{-1}\prod_{k=2}^{t}{(1-(-1)^{p(l_{k}i_{k})}q_{i_{k}}^{2l_k})}^{-1} (S,T)_K\!\!\!\! \mod{q^m \A_0}.
	\end{align*}
\end{proof}

\vskip 2mm

Let $L$ be a finitely generated $\A_0$-submodule of ${V(\lambda)}_{\lambda-\alpha}$ and set
\begin{align*}
L^{\vee}:=\{u\in {V(\lambda)}_{\lambda-\alpha}\ |\ (u,L)_K\subset \A_0\}.
\end{align*}

Similarly, let $L$ be a finitely generated $\mathbf A_0$-submodule of ${L(\infty)}_{-\alpha}$ and set
\begin{align*}
	L^{\vee}:=\{u\in {U^-_q(\g)}_{-\alpha}\ |\ (u,L)_K\subset \A_0\}.
\end{align*}

Then $(L^{\vee})^{\vee}=L$ and we obtain.

\begin{lemma}\label{piL=L}
{
	If $\lambda\gg 0$ and $\alpha\in R_{+}(r)$, we have $\pi_\lambda({L(\infty)}^\vee_{-\alpha})={L(\lambda)}_{\lambda-\alpha}^{\vee}$.
}
\end{lemma}

\begin{proof}
	Let ${\{S_k\}}$ be an $\A_0$-basis of ${L(\infty)}_{-\alpha}$, and let ${\{T_k\}}$ be its dual basis: $(S_i,T_j)_K=\delta_{ij}$. Then ${L(\infty)}_{-\alpha}^{\vee}=\sum_{j}{\A_0 T_j}$. Since $\lambda\gg 0$, from the previous lemma, we know that $(S_{i}v_{\lambda}, T_jv_{\lambda})_K\equiv (S_i,T_j)_K \mod{q\A_{0}}$.

	By Proposition \ref{proof of Hr},
	we have $L(\lambda)=\sum_{k}{\A_0 (S_kv_{\lambda})}$. 
	Hence we conclude that $\text{for}\ \lambda\gg 0,$
	\begin{equation*}
		{L(\lambda)}_{\lambda-\alpha}^{\vee}=\sum_{j\in I}{\mathbf A_0T_jv_{\lambda}}=\pi_{\lambda}({L(\infty)}_{-\alpha}^{\vee}).\ 
	\end{equation*}
\end{proof}

\vskip 3mm
\begin{lemma}\label{PsiLL}
{
	Let $\lambda\in P^{+}$, $\mu\gg 0$ and $\alpha\in R_{+}(r)$, then we have
	\begin{equation*}
		\Psi_{\lambda,\mu}({(L(\lambda)\otimes L(\mu))}_{\lambda+\mu-\alpha})\subset {L(\lambda+\mu)}_{\lambda+\mu-\alpha}.
	\end{equation*}
}
\end{lemma}

\begin{proof}
	Recall Lemma \ref{LotimesL}, we have
	\begin{equation*}
		{(L(\lambda)\otimes L(\mu))}_{\lambda+\mu-\alpha}=\sum_{(i,l)\in I^\infty}{\widetilde{f}_{il}({(L(\lambda)\otimes L(\mu))}_{\lambda+\mu-\alpha+l\alpha_i})+v_{\lambda}\otimes {L(\mu)}_{\mu-\alpha}}.
	\end{equation*}
	By induction hypothesis ${\bf D}(r-1)$, we get
	\begin{align*}
		&\Psi_{\lambda,\mu}(\sum_{(i,l)\in I^\infty}{\widetilde{f}_{il}({(L(\lambda)\otimes L({\mu}))}_{\lambda+\mu-\alpha+l\alpha_i})}\\
		=&\sum_{(i,l)\in I^\infty}{\widetilde{f}_{il}\Psi_{\lambda,\mu}({(L(\lambda)\otimes L({\mu}))}_{\lambda+\mu-\alpha+l\alpha_i})}\\
		\subset&\sum_{(i,l)\in I^\infty}{\widetilde{f}_{il}{L(\lambda+\mu)}_{\lambda+\mu-\alpha+l\alpha_i}}={L(\lambda+\mu)}_{\lambda+\mu-\alpha}.
	\end{align*}
	
	Thus it remains to show
	\begin{equation*}
		\Psi_{\lambda,\mu}(v_{\lambda}\otimes {L(\mu)}_{\mu-\alpha})\subset {L(\lambda+\mu)}_{\lambda+\mu-\alpha}.
	\end{equation*}
	
	Let $u\in {L(\lambda+\mu)}_{\lambda+\mu-\alpha}^{\vee}$. By Lemma \ref{piL=L}, we have $u=Sv_{\lambda+\mu}$ for some $S\in {L(\infty)}_{-\alpha}^{\vee}$.
	Note that
	\begin{equation*}
		\Delta(S)=S\otimes\mathbf 1+(\text{intermediate terms})+K_{\alpha}\otimes S,
	\end{equation*}
	
	which yields
	\begin{align*}
		&(\Phi_{\lambda,\mu}(u),v_{\lambda}\otimes {L(\mu)}_{\mu-\alpha})
		=(\Delta(S)(v_{\lambda}\otimes v_{\mu}),v_{\lambda}\otimes {L(\mu)}_{\mu-\alpha})\\
		=&(Sv_{\lambda}\otimes v_{\mu}+(\text{intermediate terms})+K_{\alpha}v_{\lambda}\otimes Sv_{\mu}, v_{\lambda}\otimes {L(\mu)}_{\mu-\alpha})\\
		=&q^{(\alpha,\lambda)}(Sv_{\mu},{L(\mu)}_{\mu-\alpha}).
	\end{align*}
	
	Since $\mu\gg 0$, Lemma \ref{piL=L} implies that $Sv_{\mu}\in {L(\mu)}^{\vee}$. By Corollary \ref{4.7}, we have 
	\begin{align*}
		(u,\Psi_{\lambda,\mu}(v_{\lambda}\otimes {L(\mu)}_{\mu-\alpha}))=(\Phi_{\lambda,\mu}(u),v_{\lambda}\otimes {L(\mu)}_{\mu-\alpha})=q^{(\alpha,\lambda)}(Sv_{\mu},{L(\mu)}_{\mu-\alpha})\subset \A_0.
	\end{align*}
	
	Hence $\Psi_{\lambda,\mu}(v_{\lambda}\otimes {L(\mu)}_{\mu-\alpha})\subset {({L(\lambda+\mu)}_{\lambda+\mu-\alpha}^{\vee})}^\vee={L(\lambda+\mu)}_{\lambda+\mu-\alpha}$.
\end{proof}

\vskip 2mm

\begin{proposition} \label{prop: proof of Fr} 	{ $({\bf F}(r))$ \
	Let $\alpha\in R_{+}(r)$ and $b\in {B(\infty)}_{-\alpha}$.
	If $\widetilde{e}_{il}b\ne 0$, then $b=\widetilde{f}_{il}\widetilde{e}_{il}b$.	
}
\end{proposition}

\begin{proof}
	Let $b=\widetilde{f}_{i_1l_1}\cdots\widetilde{f}_{i_tl_t}\mathbf 1\in {B(\infty)}_{-\alpha}$. We assume $\widetilde{e}_{il}b\ne 0$. If $i_1=\cdots=i_t$ and $i\ne i_{1}$, then
	\begin{equation*}
		\widetilde{e}_{il}b=\widetilde{e}_{il}\widetilde{f}_{i_1l_1}\cdots\widetilde{f}_{i_tl_t}\mathbf 1=(-1)^{y}\widetilde{f}_{i_1l_1}\cdots\widetilde{f}_{i_tl_t}\widetilde{e}_{il}\mathbf 1=0,\quad\text{for some} \ y\in \{0,1\}.
	\end{equation*}
	Hence we must have $i=i_1=\cdots=i_t$.
	In this case, our assertion follows easily.
	
	Assume that there exists $s$ with $1\leq s<t$ such that $i_s\ne i_{s+1}=\cdots=i_t$. Take $\mu\gg 0$ and set $\lambda_0=\Lambda_{i_s},~\lambda=\lambda_0+\mu$.
	
	Then Lemma \ref{ffvvbb} yields
	\begin{equation*}
		\widetilde{f}_{i_1l_1}\cdots\widetilde{f}_{i_tl_t}	(v_{\lambda_0}\otimes v_{\mu})\equiv  v\otimes v'\!\!\!\! \mod{q(L(\lambda_0)\otimes L(\mu))}
	\end{equation*}
	for some $v\in {L(\lambda_0)}_{\lambda_0-\beta}$, $v'\in {L(\mu)}_{\mu-\gamma}$, $\beta,\gamma\in \mathrm Q_{+}(r-1)\backslash\{0\}$ and $\alpha=\beta+\gamma$ such that
	$
	v+qL(\lambda_0)\in B(\lambda_0)\cup\{0\},\ v'+qL(\mu)\in B(\mu)\cup\{0\}.
$
	
	By Corollary \ref{fLLeLL} (b), we have
	\begin{equation*}
		\widetilde{e}_{il}(\widetilde{f}_{i_1l_1}\cdots\widetilde{f}_{i_tl_t}(v_{\lambda_0}\otimes v_{\mu}))\equiv \widetilde{e}_{il}(v\otimes v')\!\!\!\! \mod{q (L(\lambda_0)\otimes L(\mu))}.
	\end{equation*}
	
	Applying $\Psi_{\lambda_0,\mu}$, $\mathbf{H}(r-1)$ yields
	\begin{align*}
		\pi_{\lambda}(\widetilde{e}_{il}\widetilde{f}_{i_1l_1}\cdots \widetilde{f}_{i_tl_t} \mathtt 1)&=\widetilde{e}_{il}\widetilde{f}_{i_1l_1}\cdots\widetilde{f}_{i_tl_t}v_{\lambda_0+\mu}\equiv \Psi_{\lambda_0,\mu}(\widetilde{e}_{il}(v\otimes v'))\!\!\!\! \mod{qL(\lambda)}.
	\end{align*}
Thus $\widetilde{e}_{il}(v\otimes v')\notin q(L(\lambda_0)\otimes L(\mu))$.
By Lemma \ref{a-g} (b), we have
	\begin{align*}
		\widetilde{f}_{i_1l_1}\cdots\widetilde{f}_{i_tl_t}(v_{\lambda_0}\otimes v_{\mu})&\equiv v\otimes v'\equiv	 \widetilde{f}_{il}\widetilde{e}_{il}(v\otimes v')\\
		&\equiv \widetilde{f}_{il}\widetilde{e}_{il}(\widetilde{f}_{i_1l_1}\cdots\widetilde{f}_{i_tl_t}(v_{\lambda_0}\otimes v_{\mu})) \!\!\!\!\mod{q (L(\lambda_0)\otimes L(\mu))}.
	\end{align*}
	
	Therefore, applying $\Psi_{\lambda_0,\mu}$, Lemma \ref{PsiLL} implies
	\begin{align*}
		 \widetilde{f}_{i_1l_1}\cdots\widetilde{f}_{i_tl_t}v_{\lambda}&=\widetilde{f}_{il}\widetilde{e}_{il}(\widetilde{f}_{i_1l_1}\cdots\widetilde{f}_{i_tl_t}v_{\lambda})\!\!\!\! \mod{qL(\lambda)}.
	\end{align*}
	Since $\lambda\gg 0$, by Lemma \ref{fP and eP} and Corollary \ref{pi}, we have $b=\widetilde{f}_{il}\widetilde{e}_{il}b\!\! \mod{qL(\infty)}$.
\end{proof}

\vskip 2mm

\begin{proposition} \label{proof of Br} $({\bf B}(r))$ \
{
	Let $\lambda\in P^{+}$ and $\alpha\in R_{+}(r)$.
	For $b\in {B(\lambda)}_{\lambda-\alpha+l\alpha_i}$ and $b'\in {B(\lambda)}_{\lambda-\alpha}$, we have $\widetilde{f}_{il}b=b'$ if and only if $b=\widetilde{e}_{il}b'$.
	}
\end{proposition}
\begin{proof}
	Suppose $\widetilde{f}_{il}b=b'\in B(\lambda)_{\lambda-\alpha}$. Since $b\in B(\lambda)_{\lambda-\alpha+l\alpha_i}$, by Lemma \ref{uc=0}, there exists $\mathbf c\in \mathcal C_i$ with $|\mathbf c|\geq l$, such that
	$
	b\equiv \mathtt b_{i,\mathbf c}u_{\mathbf{c}}\mod{qL(\lambda)},\ \text{where}\ u_{\mathbf{c}}\in \mathcal{K}_i\cap L(\lambda).$
    
	If $i\notin I^{\text{iso}}$, then
	\begin{align*}
	&\widetilde{f}_{il}b= b_{i,(l,\mathbf c)}u_{\mathbf{c}}=b',\\
	&\widetilde{e}_{il}b'=\widetilde{e}_{il} b_{i,(l,\mathbf c)}u_{\mathbf{c}}= b_{i,\mathbf c}u_{\mathbf{c}}=b.
	\end{align*}
If $i\in I^{\text{iso}}$, then
	$
		\widetilde{f}_{il}b= \frac{1}{\mathbf c_l+1} \, \mathtt b_{i,(l,\mathbf c)}u_{\mathbf{c}}=b'.
$
	Hence
	\begin{align*}
	\widetilde{e}_{il}b'=\frac{\mathbf c_l+1}{\mathbf c_l+1}\mathtt b_{i,\mathbf c} u_{\mathbf{c}}=b.
	\end{align*}
	
	Conversely, suppose $b'\in {B(\lambda)}_{\lambda-\alpha}$ and $b=\widetilde{e}_{il}b'\in {B(\lambda)}_{\lambda-\alpha+l\alpha_i}$.
	By Corollary \ref{pi} (b), we have $b'=\overline{\pi}_{\lambda}(b_0')$ for some $b_0'\in {B(\infty)}_{-\alpha}$. By Proposition    \ref{proof of Kr}, we have	
	\begin{equation*}
		\overline{\pi}_{\lambda}(\widetilde{e}_{il}b_0')=\widetilde{e}_{il}(\overline{\pi}_{\lambda}(b_0'))=\widetilde{e}_{il}b'\ne 0.
	\end{equation*}
	Hence $\widetilde{e}_{il}b_0'\ne 0$ and Proposition \ref{prop: proof of Fr} gives $b_0'=\widetilde{f}_{il}\widetilde{e}_{il}b_0'$.
	Applying $\overline{\pi}_{\lambda}$, we get
	\begin{equation*}
		\widetilde{f}_{il}b=\widetilde{f}_{il}(\widetilde{e}_{il}b')=\widetilde{f}_{il}	 \overline{\pi}_{\lambda}(\widetilde{e}_{il}b_0')=\overline{\pi}_{\lambda}(\widetilde{f}_{il}\widetilde{e}_{il}b_0')=\overline{\pi}_{\lambda}(b_0')=b'.
	\end{equation*}
\end{proof}

\vskip 2mm

\begin{proposition}  \label{proof of Gr} $({\bf G}(r))$ \
{	Let $\lambda\in P^{+}$ and $\alpha\in R_{+}(r)$. We have the following facts.
	\begin{enumerate}	
		\item[{\rm (a)}] ${B(\lambda)}_{\lambda-\alpha}^{\circ}$ is a $\Q$-basis of ${L(\lambda)}_{\lambda-\alpha}/q{L(\lambda)}_{\lambda-\alpha}$.
		\item[{\rm (b)}] ${B(\infty)}_{-\alpha}^{\circ}$ is a $\Q$-basis of ${L(\infty)}_{-\alpha}/q{L(\infty)}_{-\alpha}$.
	\end{enumerate}
}
\end{proposition}

\begin{proof}
We will prove (a) only.	Suppose that we have  $\Q$-linear dependence relation : $$\sum_{b\in {B(\lambda)}_{\lambda-\alpha}^{\circ}}{a_bb=0},$$ where $a_b\in \Q.$
	
Since $\widetilde{e}_{il}{B(\lambda)}_{\lambda-\alpha}\subset B(\lambda)\cup \{0\}$ for any $(i,l)\in I^{\infty}$, we have
	\begin{align*}
	\widetilde{e}_{il}(\sum_{b \in {B(\lambda)}_{\lambda-\alpha}^{\circ}}{a_bb})=\sum_{ b\in {B(\lambda)}_{\lambda-\alpha}^{\circ},\atop   \widetilde{e}_{il}b\ne 0}(-1)^{t_{il}(b)}
	{a_b((-1)^{t_{il}(b)}\widetilde{e}_{il}b)}=0.
	\end{align*}
	when $\widetilde{e}_{il}b\in B(\lambda)_{\lambda-\alpha+l\alpha_{i}}$, $t_{il}(b)=0$, otherwise, $t_{il}(b)=1$

    By ${\bf G}(r-1)$ and Proposition \ref{proof of Br},
	we have $a_b=0$ whenever  $\widetilde{e}_{il}b\ne 0$.
 For each $b\in {B(\lambda)}_{\lambda-\alpha}^{\circ}$, there exists $(i,l)\in I^{\infty}$ such that $\widetilde{e}_{il}b\ne 0$.
Thus $a_b=0$ for any $b\in {B(\lambda)}_{\lambda-\alpha}^{\circ}$. Hence, the proposition holds.
\end{proof}

\vskip 2mm

\begin{lemma}\label{uLVLU}
{
	Let $\lambda\in P^{+}$ and $\alpha\in \Q_{+}(r)\backslash \{0\}$.
	\begin{enumerate}	
		\item[{\rm (a)}] If $u\in {L(\lambda)}_{\lambda-\alpha}/ q{L(\lambda)}_{\lambda-\alpha}$ and $\widetilde{e}_{il}u=0$ for any $(i,l)\in I^{\infty}$, then $u=0$.	
		\item[{\rm (b)}]
		If $u\in {V(\lambda)}_{\lambda-\alpha}$ and $\widetilde{e}_{il}u\in L(\lambda)$ for any $(i,l)\in I^{\infty}$, then $u\in {L(\lambda)}_{\lambda-\alpha}$.
		\item[{\rm (c)}]
		If $u\in {L(\infty)}_{-\alpha}/q{L(\infty)}_{-\alpha}$ and $\widetilde{e}_{il}u=0$ for any $(i,l)\in I^{\infty}$, then $u=0$.
		\item[{\rm (d)}]
		If $u\in {U^-_q(\g)}_{-\alpha}$ and $\widetilde{e}_{il}u\in L(\infty)$ for any $(i,l)\in I^{\infty}$, then $u\in {L(\infty)}_{-\alpha}$.
	\end{enumerate}
}
\end{lemma}
\begin{proof}
	(a) Let $u=\sum_{b\in {B(\lambda)}_{\lambda-\alpha}^{\circ}}{a_bb=0}$ ($a_b\in \Q$).
	For any $(i,l)\in I^{\infty}$, we have
	\begin{equation*}
		\widetilde{e}_{il}u=\sum_{ b\in {B(\lambda)}_{\lambda-\alpha}^{\circ},\atop   \widetilde{e}_{il}b\ne 0}(-1)^{t_{il}(b)}
	{a_b((-1)^{t_{il}(b)}\widetilde{e}_{il}b)=0}.
	\end{equation*}
	when $\widetilde{e}_{il}b\in B(\lambda)_{\lambda-\alpha+l\alpha_{i}}$, $t_{il}(b)=0$, otherwise, $t_{il}(b)=1$.
By the same argument as in the proof of Proposition  \ref{proof of Gr},
	all $a_b=0$, which yields $u=0$.
	
\vskip 2mm
	
	(b) Choose the smallest $N\geq 0$ such that $q^{N}u\in L(\lambda)$.
	If $N>0$, we have
	\begin{align*}
	\widetilde{e}_{il}(q^{N}u)=q^{N}(\widetilde{e}_{il}u)\in qL(\lambda)
	\end{align*}
for all $(i,l)\in I^{\infty}$.
		By (a), we have
	$q^{N}u\in qL(\lambda)$, i.e., $q^{N-1}u\in L(\lambda)$ which contradicts to the minimality of $N$.
		Hence $N=0$ and $u\in L(\lambda)$.
	The proofs of (c) and (d) are similar to that of (a) and (b).
\end{proof}

\vskip 2mm

\begin{proposition} \label{proof of Jr} $({\bf J}(r))$ \
{
	Let $\lambda\in P^{+}$ and $\alpha\in  R_{+}(r)$, then we have
	\begin{equation*}
		B_{-\alpha}^{\lambda}:=\{b\in {B(\infty)}_{-\alpha}\mid \overline{\pi}_{\lambda}(b)\ne 0\}\xrightarrow{\sim} {B(\lambda)}_{\lambda-\alpha}.
	\end{equation*}
}
\end{proposition}
\begin{proof}
If $\alpha=0$, it is trivial. We assume $\alpha\neq0$ and we already know that $$\overline{\pi}_{\lambda}B(\infty)_{-\alpha}\backslash \{0\}=B(\lambda)_{\lambda-\alpha}.$$
For $b,b'\in B(\infty)_{-\alpha}$, $\overline{\pi}_{\lambda}(b)=\overline{\pi}_{\lambda}(b')\neq0$. Choose $(i,l)\in I^{\infty}$ such that $\widetilde{e}_{il}\overline{\pi}_{\lambda}b\neq 0$. From Proposition\ref{proof of Kr}, we have 
$$\overline{\pi}_{\lambda}(\widetilde{e}_{il}b)=\overline{\pi}_{\lambda}(\widetilde{e}_{il}b')\neq0$$
By induction hypothesis  ${\bf J} (r-1)$ and Proposition \ref{proof of Ar}, we have $\widetilde{e}_{il}b=\widetilde{e}_{il}b'$. Then $b=\widetilde{f}_{il}\widetilde{e}_{il}b=\widetilde{f}_{il}\widetilde{e}_{il}b'=b'$. 
\end{proof}
\vskip 2mm

Using all the statements we have proved so far, we can show that Lemma \ref{a-g} holds for all
 $\alpha\in R_+(r)$.

 \vskip 2mm

 In particular, we have

\begin{lemma}\label{eBB}
{
Let $\lambda,\mu\in  P^+$ and $\alpha\in R_+(r)$.

\vskip 2mm

\begin{enumerate}
\item[{\rm (a)}] For all $(i,l)\in I^{\infty}$, we have
\begin{align*}
\widetilde{e}_{il}{(B(\lambda)\otimes B(\mu))}_{\lambda+\mu-\alpha}\subset (B(\lambda)\otimes B(\mu))\cup \{0\}.
\end{align*}
\item[{\rm (b)}] If $b\otimes b'\in {(B(\lambda)\otimes B(\mu))}_{\lambda+\mu-\alpha}$ and $\widetilde{e}_{il}(b\otimes b')\neq 0$, then we have
\begin{align*}
b\otimes b'=\widetilde{f}_{il}\widetilde{e}_{il}(b\otimes b').
\end{align*}
	\end{enumerate}
}
\end{lemma}

\vskip 2mm

\begin{proposition} \label{proof of Dr} $({\bf D}(r))$
{
For every $\lambda,\mu\in P^+$ and $\alpha\in R_+(r)$, we have
\begin{enumerate}
\item[\rm (a)] $\Psi_{\lambda,\mu}({(L(\lambda)\otimes L(\mu))}_{\lambda+\mu-\alpha})\subset L(\lambda+\mu)$,
\item[\rm (b)] $\Psi_{\lambda,\mu}({(B(\lambda)\otimes B(\mu))}_{\lambda+\mu-\alpha})\subset B(\lambda+\mu)\cup \{0\}$.
\end{enumerate}
}
\end{proposition}

\begin{proof}
(a)Assume that $\alpha\neq0$. By Corollary \ref{fLLeLL}, we have 
$$\widetilde{e}_{il}{((L(\lambda)\otimes L(\mu))}_{\lambda+\mu-\alpha})\subset (L(\lambda)\otimes L(\mu))_{\lambda+\mu-\alpha+l\alpha_{i}}$$
for every $(i,l)\in I^{\infty}$. Apply $\Psi_{\lambda,\mu}$, ${\bf D}(r-1)$ yields
$$\widetilde{e}_{il}\Psi((L(\lambda)\otimes L(\mu))_{\lambda+\mu-\alpha})\subset \Psi_{\lambda,\mu}((L(\lambda)\otimes L(\mu))_{\lambda+\mu-\alpha+l\alpha_{i}})\subset L(\lambda+\mu)$$
for every $(i,l)\in I^{\infty}$. Then, (a) follows from Lemma \ref{uLVLU}(b). 

(b) We assume $\alpha\neq 0$, let $b\otimes b'\in (B(\lambda)\otimes B(\mu))_{\lambda+\mu-\alpha}$. If $\widetilde{e}_{il}(b\otimes b')=0$ for every $(i,l)\in I^{\infty}$, then $\widetilde{e}_{il}\Psi_{\lambda,\mu}(b\otimes b')=0$ for every $(i,l)\in I^{\infty}$. By Lemma \ref{uLVLU}, we have $\Psi_{\lambda,\mu}(b\otimes b')=0$.

If $\widetilde{e}_{il}(b\otimes b')\neq 0$ for some $(i,l)\in I^{\infty}$, By Lemma \ref{eBB} and ${\bf D}(r-1)$, we have 
\begin{align*}
   \Psi_{\lambda,\mu}(b\otimes b')&\equiv \Psi_{\lambda,\mu}(\widetilde{f}_{il}\widetilde{e}_{il}(b\otimes b'))\equiv \widetilde{f}_{il} \Psi_{\lambda,\mu}(\widetilde{e}_{il}(b\otimes b'))\\
   &\in\widetilde{f}_{il}(B(\lambda+\mu))\cup \{0\}\subset B(\lambda+\mu)\cup \{0\}
\end{align*}
\end{proof}
\vskip 2mm
Thus we have completed the proofs of all the statements in Kashiwara's grand-loop argument,
which proves
Theorem \ref{thm:crystal basis of V} and Theorem \ref{thm:crystal basis of U}.
In the rest of this section, we will show the properties  of polarization.
\vskip 2mm
Let $(\ ,\ )_K^0$ denote the $\Q$-valued inner product on $L(\lambda)/qL(\lambda)$ (resp. $L(\infty)/qL(\infty)$) by taking crystal limit of $(\ ,\ )_K$ on $L(\lambda)$ (resp. $L(\infty)$).

\begin{proposition}\label{lem:orthogonal}
{ The crystal $B(\lambda)$ (resp. $B(\infty)$) forms an pseudo-orthogonal basis of $L(\lambda)/qL(\lambda)$ (resp. $L(\infty)/qL(\infty)$) with respect to $(\ ,\ )_K^0$.}\end{proposition}

\begin{proof}
We first consider  $B(\lambda)^{\circ}$. For all $b,b'\in {B(\lambda)}_{\lambda-\alpha}^{\circ}$, we shall prove $(b,b')_K^0\in \delta_{b,b'} \Z$ by using induction on $\mathrm{ht}(\alpha)$, where $\alpha\in R_+(r)$.

\vskip 2mm

If $\mathrm{ht}(\alpha)=0$, then our conclusion is trivial.

\vskip 2mm

If $\mathrm{ht}(\alpha)>0$, we choose $(i,l)\in I^{\infty}$ such that $\widetilde{e}_{il}b\neq 0$.
By $\mathbf B(r)$ and Lemma \ref{lem:LLA0}, we have
\begin{equation*}
(b,b')_K^0=(\widetilde{f}_{il} \widetilde{e}_{il}b,b')_K^0 = (\widetilde{e}_{il}b,\widetilde{e}_{il}b')_K^0
\in \delta_{\widetilde{e}_{il}b,\widetilde{e}_{il}b'} \Z=\delta_{b,b'} \Z.
\end{equation*}

\vskip 2mm

By Lemma \ref{(P,Q)} and a similar approach above, it is easy to show that the crystal $B(\infty)$ is an pseudo-orthogonal basis of  $L(\infty)/qL(\infty)$ with respect to $(\ ,\ )_K^0$.
\end{proof}

\begin{lemma}
For $(i,l)\in I^{\infty}$ and $P,\ Q\in \mathbf U^{-}$, we have 
$$(Pb_{il},Q)=(-1)^{p(Q)p(li)+p(li)}(P,q^{ld_{i}\langle h_{i}, |Q|+l\alpha_{i}\rangle}e_{il}''(Q)).$$
\end{lemma}
\begin{proof}
    When $|P|=0$ and  $Q=b_{il}$, our claim follows directly from Proposition\ \ref{3.7}. Assume that $|P|>1$ and $P$ satisfies $(Pb_{il},Q)=(-1)^{p(Q)p(li)+p(li)}(P,q^{ld_{i}\langle h_{i}, |Q|+l\alpha_{i}\rangle},e_{il}''(Q))$ for any $Q$, it suffices to prove that for any $(j,k)\in I^{\infty}$,
    $$(b_{jk}Pb_{il},Q)=(-1)^{p(Q)p(li)+p(li)}(b_{jk}P,q^{ld_{i}\langle h_{i}, |Q|+l\alpha_{i}\rangle}e_{il}''(Q)).$$

By using Proposition \ref{3.7} we have 
\begin{equation}
    \begin{aligned}
        (b_{jk}Pb_{il},Q) &=(Pb_{il},e_{jk}'(Q))\\
        &=(-1)^{(p(Q)+p(kj)+1)p(li)}q^{ld_{i}\langle h_{i}, |Q|+k\alpha_{j}+l\alpha_{i}\rangle}(P,e_{il}''\circ e_{jk}'(Q))\\
        &=(-1)^{(p(Q)+1)p(li)}q^{ld_{i}\langle h_{i}, |Q|+l\alpha_{i}\rangle}(P,e_{jk}'\circ e_{il}''(Q))\\
        &=(-1)^{(p(Q)+1)p(li)}q^{ld_{i}\langle h_{i}, |Q|+l\alpha_{i}\rangle}(b_{jk}P, e_{il}''(Q)).
    \end{aligned}
\end{equation}
The third identity  upside comes from Proposition \ref{3.3} and we complete the proof.
\end{proof}

\vskip 2mm
We define an anti-involution $\varrho$ on $\mathbf U$ i.e. $\rho(xy)=\rho(y)\rho(x)$, such that 
$$\rho(a_{il})=(-1)^{p(li)}a_{il},\ \rho(b_{il})=b_{il},\ \rho(q^{h})=q^{-h}.$$
\vskip 2mm
\begin{lemma}
    We have $\rho(a_{il}'\rho(Q))=(-1)^{(p(Q)+1)p(li)}q^{ld_{i}\langle h_{i}, |Q|+l\alpha_{i}\rangle}a_{il}''(Q)$
\end{lemma}
\begin{proof}
    First of all, recall that  
    $$[a_{il},Q]=\frac{K_{i}^{l}e_{il}''(Q)-K_{i}^{-l}e_{il}'(Q)}{1-(-1)^{p(li)}q_{i}^{2l}},$$
where $e_{il}''(Q)=\varrho^{il}(Q^{+})^{-}$ and $e_{il}'(Q)=(-1)^{p(li)p(Q)-p(li)}q_{i}^{l\langle h_{i},|Q|+l\alpha_{i}\rangle}\varrho_{il}(Q^{+})^{-}$.

Applying $\rho$ to $[a_{il},Q]$, we have
\begin{equation}
    \begin{aligned}
      \rho([a_{il},Q])&= [\rho(Q),(-1)^{p(li)}a_{il}]= (-1)^{p(li)+p(Q)p(li)+1}[a_{il},\rho(Q)]\\
      &=\frac{\rho(e_{il}''(Q))K_{i}^{-l}-\rho(e_{il}'(Q))K_{i}^{l}}{1-(-1)^{p(li)}q_{i}^{2l}}\\
      &=\frac{q^{ld_{i}\langle h_{i},|Q|+l\alpha_{i}\rangle  }K_{i}^{-l}\rho(e_{il}''(Q))-q^{-ld_{i}\langle h_{i},|Q|+l\alpha_{i}\rangle  }K_{i}^{l}\rho(e_{il}'(Q))}{1-(-1)^{p(li)}q_{i}^{2l}},
    \end{aligned}
\end{equation}
which implies 
$$[a_{il},\rho(Q)]=(-1)^{p(li)+p(Q)p(li)}\frac{q^{-ld_{i}\langle h_{i},|Q|+l\alpha_{i}\rangle  }K_{i}^{l}\rho(e_{il}'(Q))-q^{ld_{i}\langle h_{i},|Q|+l\alpha_{i}\rangle  }K_{i}^{-l}\rho(e_{il}''(Q))}{1-(-1)^{p(li)}q_{i}^{2l}}.$$
Hence
    $$    e_{il}'(\rho(Q))=(-1)^{p(li)+p(Q)p(li)}q^{ld_{i}\langle h_{i},|Q|+l\alpha_{i}\rangle  }\rho(e_{il}''(Q)),$$
then
$$\rho(e_{il}'\rho(Q))=(-1)^{(p(Q)+1)p(li)}q^{ld_{i}\langle h_{i}, |Q|+l\alpha_{i}\rangle}e_{il}''(Q).$$

\end{proof}

\vskip 2mm
\begin{proposition}
For $P,\ Q \in U^{-}$, we have 
   $(\rho(P),\rho(Q))=(P,Q)$ 
\end{proposition}
\begin{proof}
    The claim is trivial when $P=1$, so it suffices to prove that 
    $$(\rho(Pb_{il}), \rho(Q))=(Pb_{il},Q).$$
  We have
\begin{equation}
    \begin{aligned}
        (\rho(Pb_{il}), \rho(Q)) &=(b_{il}\rho(P),\rho(Q))=(\rho(P),e_{il}'\rho(Q))\\
        &=(P,\rho(e_{il}'\rho(Q)))=(-1)^{(p(Q)+1)p(li)}q^{ld_{i}\langle h_{i}, |Q|+l\alpha_{i}\rangle}(P,e_{il}''(Q))\\
        &=(Pb_{il},Q).
    \end{aligned}
\end{equation}
\end{proof}

As in [FHL], we have the following proposition 
\begin{proposition}\label{5.38}
    $\rho(L(\infty))=L(\infty)$.
\end{proposition}

\section{Global Bases}

\subsection{$\A$-form of $U_{q}^{-}(\g)$}
For each $(i,l)\in I^{\infty}$, set
\begin{equation}\label{new}
A_i=(-q_{i}^{l})\mathtt a_{il},
\end{equation}
which yields the following commutation relation
\begin{equation} \label{eq:comm2}
A_{il}b_{il}-b_{il}A_{il}=\frac{K_{i}^{l}-K_i^{-l}}{(-1)^{p(li)}q_i-q_i^{-l}}.
\end{equation}

Let $\A=\Q[q,q^{-1}]$. Let $U_{\A}(\mathfrak{g})$ be the $\A$-subalgebra of $\U$  generated by $A_{i}^{(n)},\ b_{i}^{(n)}\ (i\in I^{re},\ n\ge 0)$, $A_{il},\ b_{il}\ (i\in I^{im},\ l\ge 1)$, $q^{h}\ (h\in P^{\vee})$ and $\begin{bmatrix}
    K_{i}q^{n}\\
    m
\end{bmatrix}_{i}\ (i\in I^{re},\ m\in \Z_{\ge 0 },\ n\in \Z ) $, where
$$
\begin{bmatrix}
    K_{i}q^{n}\\
    m
\end{bmatrix}_{i}=\frac{1}{[m]_{i}!}\prod_{s=1}^{m}\frac{(-1)^{p(i)}K_{i}q_{i}^{n+1-s}-K_{i}^{-1}q_{i}^{-n-1+s}}{(-1)^{p(i)}q_{i}-q_{i}^{-1}}.$$
We denote by $U_{\A}^{-}(\mathfrak{g})$ the $\A$-subalgebra of $\U$ generated by $b_{i}^{(n)}\ (i\in I^{re},\ n\ge 0)$ and $b_{il} \ (i\in I^{im},l\ge 1)$.  Then $U_{\A}(\mathfrak{g})$ and $U_{\A}^{-}(\mathfrak{g})$ is stable under the $\Q$-linear map $\Psi:\Q(q)\to \Q(q),\ q\to -q^{-}$, we also use the notation $\overline{f(q)}=\Psi(f(q))$ and call this map bar involution.

As in \cite{CHW2}, if $u\in U_{\A}^{-}(\mathfrak{g})$ and $u=\sum_{n\ge0}b_{i}^{(n)}u_{n}$
with $i\in I^{re}$ and $a_{i}'u_{n}=0$ for all $n\ge0$, then $u_{n}\in U_{\A}^{-}(\mathfrak{g}) $. If $i\in I^{im}$, we have the following lemma.
\begin{lemma}
    For each $i\in I^{im}$ and $u=\sum_{}b_{i,\mathbf{c}}u_{\mathbf{c}}\in U_{\A}^{-}(\mathfrak{g})$, we have $u_{\mathbf{c}}\in  U_{\A}^{-}(\mathfrak{g})$.
\end{lemma}
\begin{proof}
    We use a similar argument in Lemma \ref{4.12}. Assume that $u\in  U_{\A}^{-}(\mathfrak{g})$ has the form $u=mb_{i,\mathbb{c}}m'$ for some homogeneous elements $m,m'\in \mathcal{K}_{i}\cap  U_{\A}^{-}(\mathfrak{g})$, where $\mathcal{K}_{i}=\cap_{l\ge 1} ker\ e_{il}'$. We shall show that $u$ can be written into the form 
    $$u=\sum_{}b_{i,\mathbf{c}}u_{\mathbf{c}}\ \text{with}\ u_{\mathbf{c}}\in \mathcal{K}_{i}\cap U_{\A}^{-}(\mathfrak{g}).  $$

    If $|c|=0$, then $u=mm'\in \mathcal{K}_{i}\cap U_{\A}^{-}(\mathfrak{g})$. Otherwise, we have
    $$u=q^{-c_{1}(|m|,\alpha_{i})}[m,b_{ic_{1}}]^{\circ}b_{i,\mathbf{c}\backslash c_{1}}m'+q^{-c_{1}(|m|,\alpha_{i})}b_{ic_{1}}mb_{i,\mathbf{c}\backslash c_{1}}m',$$
    where $[m,b_{ic_{1}}]^{\circ}=(-1)^{p(m)p(c_{1}i)}q^{-(c_{1}\alpha_{i},|m|)}mb_{ic_{1}}-b_{ic_{1}}m\in \mathcal{K}_{i}\cap U_{\A}^{-}(\mathfrak{g}) $, our claim follows by induction on $\text{ht}(u)$.

    Next we will show that if $u\in  U_{\A}^{-}(\mathfrak{g})$, then $u$ can be written into the form  $$u=\sum_{}b_{i,\mathbf{c}}u_{\mathbf{c}}\ \text{with}\ u_{\mathbf{c}}\in \mathcal{K}_{i}\cap U_{\A}^{-}(\mathfrak{g}).  $$

    Assume that $u$ is a monomial in $U_{\A}^{-}(\mathfrak{g})$. When $\text{ht}(u)>1$, $u=b_{jk}u'$ for some $j\in I^{im}$ or $u=b_{j}^{(n)}u'$ for some $j\in I^{re}$.  By induction hypothesis, $u'=\sum b_{i,\mathbf{c}}u_{\mathbf{c}}$ with $u_{\mathbf{c}}\in \mathcal{K}_{i}\cap U_{\A}^{-}(\mathfrak{g})$. If $j\neq i$, then $u=\sum b_{j}^{(n)}b_{i,\mathbf{c}}u_{\mathbf{c}}$ or $u=\sum b_{jk}b_{i,\mathbf{c}}u_{\mathbf{c}}$ is of the form $mb_{i,\mathbf{c}}m'$ since $b_{j}^{(n)}$ or $b_{jk}$ belongs to $\mathcal{K}_{i}\cap U_{\A}^{-}(\mathfrak{g})$. If $j=i$, then $u=\sum b_{i,(k,\mathbf{c})}u_{\mathbf{c}}$ is already in the form we wanted. 
    
\end{proof}
\vskip 2mm
Set $L_{\A}(\infty)=L(\infty)\cap U_{\A}^{-}(\mathfrak{g})$, then $L_{\A}(\infty)$ is stable under Kashiwara operators $\widetilde{e}_{il}, \widetilde{f}_{il}$. Therefore, $B(\infty)\subset L_{\A}(\infty)\slash  qL_{\A}(\infty) \subset L(\infty)\slash  qL(\infty)$.
Define 
\begin{equation*}
\begin{aligned}
& (\mathtt b_i^nU_q^-(\g))^\A:=\mathtt b_i^nU_q^-(\g)\cap U^-_{\A}(\g) \quad \text{for} \ i\in I^{\text{re}} \ \text{and}\ n\geq 1,\\
& (\mathtt b_{i,\mathbf c}U_q^-(\g))^\A:=\mathtt b_{i,\mathbf c}U_q^-(\g)\cap U^-_{\A}(\g) \quad \text{for} \ i\in I^{\text{im}} \ \text{and}\ \mathbf c \in \mathscr C_i\backslash \{0\}.
\end{aligned}
\end{equation*}
\vskip 2mm
\subsection{$\A$-form of $V(\lambda)$)}
Let $V(\lambda)=\U v_{\lambda}$ be the irreducible highest weight module with highest weight $\lambda\in P^{+}$. We define $V(\lambda)^{\A}=U_{\A}^{-}(\mathfrak{g})c_{\lambda}$. let $L_{\A}(\lambda)=V(\lambda)^{\A}\cap L(\lambda)$, we have $L_{\A}(\lambda)\subset \pi_{\lambda}(L_{\A}(\infty))$ and $B(\lambda)\subset L_{\A}(\lambda)\slash  qL_{\A}(\lambda) \subset L(\lambda)\slash  qL(\lambda)$

\vskip 2mm
By the same argument in \cite[Lemma 7.1]{CHW2} and \cite[Proposition 7.3]{CHW2}, let $\pi=-1$, we have the following proposition.
\begin{proposition}
Let $i\in I^{\text{re}}$ and $\mu\in P$.
\begin{itemize}
\item [(i)] For $u \in V(\lambda)_\mu$ with $n=-\mu(h_i)\geq 1$, we have
$$u=(-1)^{p(i)\begin{pmatrix}
    n\\2
\end{pmatrix}}\sum_{k\geq n}(-1)^{k-n}\begin{bmatrix}k-1\\k-n\end{bmatrix}_i\mathtt b_i^{(k)}A_i^{(k)}u.$$
\item [(ii)] If $n=-\mu(h_i)\geq 0$, then we have
$$V(\lambda)_{\mu}^{\mathbb A}=\sum_{k\geq n}\mathtt b_i^{(k)}V(\lambda)^{\mathbb A}_{\mu+k\alpha_i}=(\mathtt  b_i^{ n} V(\lambda))_\mu^{\mathbb A}.$$
\end{itemize}
\end{proposition}

\subsection{Existence of canonical bases.}
\vskip 2mm

Let $V$ be a finite-dimensional vector space over $\Q(q)$, $M$ be  an $\A$-submodule of $V$,  and $L_0$ (resp. $L_\infty$) be a free $\A_0$-submodule (resp. free $\A_\infty$-submodule) of $V$ such that $V\cong \Q(q)\otimes_{\A_0}L_0 \cong \Q(q)\otimes_{\mathbb A_\infty}L_\infty$.

\vskip 3mm

\begin{lemma}\cite[Lemma 7.1.1]{Kas91}\label{K1} 
 Let $V, M, L_0, L_\infty$  be as above.
\begin{itemize}
\item [(i)] Assume that the canonical map $M\cap L_0 \cap L_\infty\rightarrow M\cap L_0/M\cap qL_0$ is an isomorphism.  Then
    $$M\cap L_0\cong \Q[q]\otimes_{\Q}(M\cap L_0 \cap L_\infty),$$
    $$M\cap L_\infty \cong \Q[q^{-1}]\otimes_{\Q}(M\cap L_0 \cap L_\infty),$$
    $$M\cong \A\otimes_{\Q}(M\cap L_0 \cap L_\infty),$$
    $$M\cap L_0 \cap L_\infty\cong (M\cap L_\infty)/(M\cap q^{-1}L_\infty),$$
    $$M\cap L_0 \cap L_\infty\simeq (\Q(q)\otimes_{\A}M)\cap L_0/(\Q(q)\otimes_{\A}M)\cap qL_0.$$
\item [(ii)] Let $E$ be  an  $\Q$-vector space and $\varphi:E\rightarrow M\cap L_0 \cap L_\infty$ a homomorphism. Assume that
$M=\A\varphi(E)$ and $E\rightarrow L_0/qL_0, E\rightarrow L_\infty/q^{-1}L_\infty$ are injective. Then $E\rightarrow M\cap L_0 \cap L_\infty \rightarrow M\cap L_0/M \cap qL_0$ are isomorphisms.
\end{itemize}
\end{lemma}

\vskip 3mm

\begin{lemma}\cite[Lemma 7.1.2]{Kas91}\label{K2} 
 Let $V, M, L_0, L_\infty$ be as above  and  let $N$ be an $\A$-submodule of $M$.  Assume that
\begin{itemize}
\item [(1)] $N\cap L_0\cap L_\infty\cong N\cap L_0/N\cap qL_0$.
\item [(2)] There exist an $\Q$-vector space $F$ and a homomorphism $\varphi:F\rightarrow M\cap(L_0+N)\cap (L_\infty+N)$ such that\\
    (a) $M=\A\varphi(F)+N$,\\
    (b) the induced homomorphisms  $\varphi_0:F\rightarrow (L_0+N) /  (qL_0+N)$ and $\varphi_\infty:F\rightarrow (L_\infty+N) / (q^{-1}L_\infty+N)$ are injective.
\end{itemize}

\vskip 2mm

Then the following statements hold.

\vskip 2mm

\begin{itemize}
\item [(i)] $M\cap L_0\cap L_\infty\rightarrow M\cap L_0/M\cap qL_0$ is an isomorphism.
\item [ii)] $M\cap L_0/M\cap qL_0 \cong F\oplus (N\cap L_0/N\cap qL_0)$.
\end{itemize}
\end{lemma}
 We will prove the following  inductive statements on $r\geq 0$.
\vskip 3mm
\begin{itemize}
\item [X($r$)]: For any $\alpha\in Q_+(r)$, we have the following canonical isomorphism
$$U^-_{\A}(\g)_{-\alpha}\cap L(\infty)\cap L(\infty)^- \xrightarrow{\sim} \frac{U^-_{\A}(\g)_{-\alpha}\cap  L(\infty)}{U^-_{\A}(\g)_{-\alpha}\cap qL(\infty)}\xrightarrow{\sim}  L(\infty)_{-\alpha}/q L(\infty)_{-\alpha}.$$
\end{itemize}
We denote by $G_\infty$ the inverse of this isomorphism.

\vskip 2mm

\begin{itemize}
\item [Y($r$)]: For any $\alpha\in Q_+(r)$ and $\lambda\in P^{+}$, we have the following canonical isomorphism
$$V(\lambda)^{\A}_{\lambda-\alpha}\cap  L(\lambda)\cap L(\lambda)^- \xrightarrow{\sim} \frac{V(\lambda)^{\A}_{\lambda-\alpha}\cap  L(\lambda)}{V(\lambda)^{\A}_{\lambda-\alpha}\cap q L(\lambda)}\xrightarrow{\sim}  L(\lambda)_{\lambda-\alpha}/q L(\lambda)_{\lambda-\alpha}.$$
\end{itemize}
We denote by $G_\lambda$ the inverse of this isomorphism.

\vskip 2mm

\begin{itemize}
\item [Z($r$)]: For $\alpha\in Q_+(r)$, $(i,l)\in I^{\infty}$, and $n\geq 0$, assume that $b\in \widetilde{f}_{il}^n \, \mathcal B(\infty)_{-\alpha+l n\, \alpha_i}$.  Then we have
    $G_\infty(b)\in b_{il}^n \, U_q^-(\g).$
\end{itemize}

\vskip 3mm

If $r=0$, our assertions are obvious. Now we assume that A($r-1$), B($r-1$) and C($r-1$) are true.  Then by Lemma \ref{K1}, we have  the following result.

\vskip 3mm

\begin{lemma}
For $\alpha\in Q_+(r-1)$, we have
$$\begin{aligned}
& U^-_{\A}(\g)_{-\alpha}\cap  L(\infty)=\bigoplus_{b\in  B(\infty)_{-\alpha}}\Q[q]G_\infty(b), \quad
U^-_{\A}(\g)_{-\alpha}=\bigoplus_{b\in  B(\infty)_{-\alpha}}\A G_\infty(b),\\
& V(\lambda)^{\A}_{\lambda-\alpha}\cap  L(\lambda)=\bigoplus_{b\in B(\lambda)_{\lambda-\alpha}}\Q[q]G_\lambda(b), \quad
V(\lambda)^{\A}_{\lambda-\alpha}=\bigoplus_{b\in B(\lambda)_{\lambda-\alpha}}\A G_\lambda(b),
\end{aligned}$$
and
$$G_\infty(b)v_\lambda=G_\lambda(\overline \pi_\lambda (b)).$$
\end{lemma}

\begin{lemma}
For $\alpha\in Q_+(r-1)$, $b\in \mathcal B(\infty)_{-\alpha}$ (resp. $b\in \mathcal B(\lambda)_{\lambda-\alpha}$), we have $\overline{G_\infty(b)}=G_\infty(b)$ (resp. $\overline{G_\lambda(b)}=G_\lambda(b)$).
\begin{proof}
Let $Q=(G_\infty(b)-\overline{G_\infty(b)})/(-q-q^{-1})$. Then we have $Q\in U^-_{\mathbb A}(\g)_{-\alpha}\cap q\mathcal L(\infty)\cap \mathcal L(\infty)^-$ since $1/(-q-q^{-1})\in q\A_0$. By $A(r-1)$, the canonical isomorphism maps $Q$ to $0$. Hence, $Q=0$ and $\overline{G_\infty(b)}=G_\infty(b)$. 
\end{proof}
\end{lemma}

Let $i\in I^{\text{iso}}$, $\lambda \in P^{+}$ and $\alpha \in Q_{+}$.
For each partition  $\mathbf {c}=(1^{l_1}2^{l_2}\cdots k^{l_k}\cdots)$,
we define
\begin{equation*}
\begin{aligned}
& ({ b}_{i, \mathbf{c}} * U_{q}^{-}(\g))^{\A}_{-\alpha}
:= \sum_{k \ge 1} ({b}_{i,k}^{l_k} U_{q}^{-}(\g))^{\A}_{-\alpha}
= \sum_{k \ge 1} {b}_{i,k}^{l_k} (U_{\A}^{-}(\g)_{-\alpha + k l_{k} \alpha_{i}}), \\
& ({b}_{i, \mathbf{c}} * V(\lambda))^{\A}_{\lambda - \alpha} :
= \sum_{k \ge 1} ({b}_{i, k}^{l_k} V(\lambda))^{\A}_{\lambda - \alpha}
=\sum_{k \ge 1} {\mathtt b}_{i, k}^{l_k}(V(\lambda)^{\A}_{\lambda - \alpha + k l_k \alpha_{i}}).
\end{aligned}
\end{equation*}

Here $( b_{i,k}^{l_k}U^-_q(\g))^{\A}_{-\alpha}
= ({b}_{i,k}^{l_k} V(\lambda))^{\A}_{\lambda - \alpha} = 0   $ if $l_k=0$ or $p(i)=1$ and $l_{k}>1$.

\vskip 3mm
By a  similar argument in \cite[Proposition 5.5]{FKKT22}, we have 
\begin{proposition}\label{M}{\rm
Let $\alpha\in Q_+(r)$ and $\lambda\in P^+$.
\begin{itemize}
\item [(i)] For $i\in I^{\text{re}}$ and $n\geq 1$, we have
$$( b_i^nV(\lambda))^{\A}_{\lambda-\alpha}\cap  L(\lambda)\cap  L(\lambda)^-\xrightarrow{\sim}\frac{( b_i^nV(\lambda))^{\A}_{\lambda-\alpha}\cap L(\lambda)}{( b_i^nV(\lambda))^{\A}_{\lambda-\alpha}\cap q L(\lambda)}\cong \bigoplus_{b\in  B(\lambda)_{\lambda-\alpha}\cap \widetilde f_i^n B(\lambda)}\Q b.$$
\item [(ii)] For $i\in I^{\text{im}}\backslash I^{\text{iso}}$ and any composition $\mathbf c$ with $|\mathbf c|\neq 0$, we have
$$( b_{i,\mathbf c}V(\lambda))^{\A}_{\lambda-\alpha}\cap  L(\lambda)\cap  L(\lambda)^-\xrightarrow{\sim}\frac{( b_{i,\mathbf c}V(\lambda))^{ \A}_{\lambda-\alpha}\cap  L(\lambda)}{(b_{i,\mathbf c}V(\lambda))^{\mathbb A}_{\lambda-\alpha}\cap qL(\lambda)}\cong \bigoplus_{b\in \ B(\lambda)_{\lambda-\alpha}\cap \widetilde f_{i,\mathbf c} B(\lambda)}\Q b.$$
\item [(iii)] For $i\in I^{\text{iso}}$ and any partition $\mathbf c=1^{l_1}2^{l_2}\cdots k^{l_k}\cdots$, we have
$$( b_{i,\mathbf c}*V(\lambda))^{\A}_{\lambda-\alpha}\cap L(\lambda)\cap  L(\lambda)^-\xrightarrow{\sim}\frac{(b_{i,\mathbf c}*V(\lambda))^{\A}_{\lambda-\alpha}\cap  L(\lambda)}{(b_{i,\mathbf c}*V(\lambda))^{\A}_{\lambda-\alpha}\cap q L(\lambda)}\cong \bigoplus_{b\in B(\lambda)_{\lambda-\alpha}\cap(\widetilde f_{i,\mathbf c}*B(\lambda))}\Q b,$$
where $\widetilde f_{i,\mathbf c}*\mathcal B(\lambda):=\bigcup_{k\geq 1}  {\widetilde f}_{i,k}^{l_k}\mathcal B(\lambda)$.
\end{itemize}

}\end{proposition}

\vskip 3mm
By a  similar argument in \cite[Corollary 5.6]{FKKT22}, we have 
\begin{corollary}\label{M1}
Let $\alpha\in Q_+(r)$.
\begin{itemize}
\item [(i)] For $i\in I^{\text{re}}$ and $n\geq 1$, we have
$$( b_i^nU^-_q(\g))^{\A}_{-\alpha}\cap L(\infty)\cap  L(\infty)^-\xrightarrow{\sim}\frac{(b_i^nU^-_q(\g))^{\A}_{-\alpha}\cap  L(\infty)}{( b_i^nU^-_q(\g))^{\A}_{-\alpha}\cap q L(\infty)}\cong \bigoplus_{b\in  B(\infty)_{-\alpha}\cap \widetilde f_i^n B(\infty)}\Q b.$$
\item [(ii)] For $i\in I^{\text{im}}\backslash I^{\text{iso}}$ and $|\mathbf c|\neq 0$, we have
$$( b_{i,\mathbf c}U^-_q(\g))^{\A}_{-\alpha}\cap  L(\infty)\cap  L(\infty)^-\xrightarrow{\sim}\frac{( b_{i,\mathbf c}U^-_q(\g))^{\A}_{-\alpha}\cap  L(\infty)}{( b_{i,\mathbf c}U^-_q(\g))^{\A}_{-\alpha}\cap q L(\infty)}\cong \bigoplus_{b\in  B(\infty)_{-\alpha}\cap \widetilde f_{i,\mathbf c} B(\infty)}\Q b.$$
\item [(iii)]  For $i\in I^{\text{iso}}$ and any partition  $\mathbf c=(1^{l_1}2^{l_2}\cdots k^{l_k}\cdots)$, we have
$$( b_{i,\mathbf c}*U^-_q(\g))^{\A}_{-\alpha}\cap  L(\infty)\cap  L(\infty)^-\xrightarrow{\sim}\frac{(b_{i,\mathbf c}*U^-_q(\g))^{\A}_{-\alpha}\cap  L(\infty)}{(b_{i,\mathbf c}*U^-_q(\g))^{\A}_{-\alpha}\cap q L(\infty)}\cong \bigoplus_{b\in  B(\infty)_{-\alpha}\cap(\widetilde f_{i,\mathbf c}*B(\infty))}\Q b,$$
where $\widetilde f_{i,\mathbf c}*\mathcal B(\infty)=\bigcup_{k\geq 1}\widetilde f_{ik}^{l_k}\mathcal B(\infty)$.
\end{itemize}
\end{corollary}

\vskip 3mm

For $\alpha\in Q_+(r)$ and $(i,l)\in I^{\infty}$, We denote by $G_{il}$ the inverse of the isomorphism
$$( b_{il}U^-_q(\g))^{\A}_{-\alpha}\cap  L(\infty)\cap  L(\infty)^-\xrightarrow{\sim} \bigoplus_{b\in  B(\infty)_{-\alpha}\cap \widetilde f_{il} B(\infty)}\Q b.$$
Then Corollary \ref{M1} and Lemma \ref{K1} (i) implies that $( b_{il}^nU^-_q(\g))^{ \A}_{-\alpha}=\bigoplus_{b\in  B(\infty)_{-\alpha}\cap \widetilde f_{il}^n B(\infty)}\A G_{il}(b)$ for any $n\geq 1$.

\vskip 3mm

\begin{lemma}\label{6.9}{
Let $(i,l),(j,s)\in I^{\infty}$, $\alpha\in Q_+(r)$ and $b\in \widetilde f_{il}\mathcal B(\infty)\cap\widetilde f_{js}\mathcal B(\infty)\cap \mathcal B(\infty)_{-\alpha}$. Then we have
$$G_{il}(b)=G_{js}(b).$$
\begin{proof}
 Let us write $b=\widetilde{f}_{il}\cdots\widetilde{f}_{i_{k}l_{k}}\cdot 1$, where $(k,h)\in I^{\infty}$. If $i_k\in I^{\text{re}}$, then our claim can be proved by the same argument in \cite{Kas91}. So we  assume that $i_k\in I^{\text{im}}$. Take $\lambda\in P^+$ with $\langle h_{i_k} ,\lambda \rangle =0$ and $\langle h_{j}, \lambda \rangle \gg 0$ for all
$j\in I\backslash\{k\}$. Then 
$$V(\lambda)_{\lambda-\alpha}\simeq U^-_q(\g)_{-\alpha} \bigg/\sum_{n\geq 1}U^-_q(\g)_{-\alpha+n\alpha_k}\mathtt b_{i_kn}.$$
The same argument in \cite[Lemma 7.5.1]{Kas91} shows that
$$Q=G_{il}(b)-G_{js}(b)\in (\sum_{n\geq 1}U^-_q(\g)_{-\alpha+n\alpha_k}b_{i_kn})\cap U^-_{\A}(\g)_{-\alpha}\cap qL(\infty)\cap L(\infty)^-.$$
 Proposition \ref{5.38} implies
$$\rho(Q)\in (\sum_{n\geq 1}b_{i_kn}U^-_q(\g)_{-\alpha+n\alpha_{i_k}})\cap U^-_{\A}(\g)_{-\alpha}\cap q L(\infty)\cap L(\infty)^-.$$

If $k\in I^{\text{im}}\backslash I^{\text{iso}}$, we assume that $\rho(Q)= b_{i_k 1}u_1+\cdots+ b_{i_ks}u_s$. Since $\overline{\rho(Q)}=\rho(\overline Q)=\rho(Q)$, we have $\rho(Q)= b_{k1}\overline u_1+\cdots+b_{ks}\overline u_s$. Note that for each $1\leq l\leq s$, $ b_{i_kl}u_l=\widetilde f_{i_kl}\widetilde e_{i_kl}\rho(Q)\in  b_{i_kl}U^-_q(\g)\cap U^-_{\A}(\g)_{-\alpha}\cap qL(\infty)$ and $u_l=\overline u_l=\widetilde e_{i_kl}\rho(Q)$. Hence
$b_{i_kl}u_l\in (b_{i_kl}U^-_q(\g))^{\A}_{-\alpha}\cap q L(\infty)\cap  L(\infty)^-$ and Corollary \ref{M1} (ii) implies $\rho(Q)=0$.

If $k\in I^{\text{iso}}$, since $\rho(Q)\in (\sum_{n\geq 1} b_{i_kn} U^-_q(\g)_{-\alpha+n\alpha_k})\cap U^-_{\A}(\g)_{-\alpha}$, 
$\rho(Q)$ can be expressed as the form $\rho(Q)= b_{i_k1}u_1+\cdots+ b_{i_ks}u_s$ with $$u_l=\sum_{\substack{\mathbf c\in \mathscr C_{i_k}\ \text{and}  \\ \mathbf c\ \text{contains no} \ l+1,\cdots, s}} b_{k,\mathbf c}u_{\mathbf c}.$$
 For every $1\leq l\leq s$, we have $$ b_{i_kl}u_{l}=\widetilde f_{i_kl}\widetilde e_{i_kl}(\rho(Q)-\sum_{l<p\leq s}b_{i_kp}u_p)\in ( b_{i_kl}U^-_q(\g))^{\A}_{-\alpha}.$$
Hence $\rho(Q)\in ( b_{k,(1,\dots, s)}*U^-_q(\g))^{\A}_{-\alpha}$, and Corollary \ref{M1} (iii) implies $\rho(Q)=0$.
\end{proof}
}\end{lemma}


\vskip 3mm
\begin{proposition}\label{pfAr}  $({X}(r))$
     For any $\alpha\in Q_+(r)$,  the following canonical map is a  isomorphism
$$U^-_{\A}(\g)_{-\alpha}\cap L(\infty)\cap L(\infty)^- \xrightarrow{\sim} \frac{U^-_{\A}(\g)_{-\alpha}\cap  L(\infty)}{U^-_{\A}(\g)_{-\alpha}\cap qL(\infty)}\xrightarrow{\sim}  L(\infty)_{-\alpha}/q L(\infty)_{-\alpha}.$$
\end{proposition}
\begin{proof}

By Lemma \ref{6.9},  we can give a well-defined map
$$G:\mathcal L(\infty)_{-\alpha}/q\mathcal L(\infty)_{-\alpha}\rightarrow U^-_{\mathbb A}(\g)_{-\alpha}\cap \mathcal L(\infty)\cap \mathcal L(\infty)^-$$
with 
$$G(b)= G_{il}(b) \ \ \text{for} \  b\in \widetilde f_{il} B(\infty)\cap  B(\infty)_{-\alpha},  (i,l) \in I^{\infty}.$$
Thus we have $b=G(b)+q\mathcal L(\infty)$ and for any $n\geq 1$,
\begin{equation}\label{G(b)}
(b_{il}^nU^-_q(\g))^{\A}_{-\alpha}=\bigoplus_{b\in B(\infty)_{-\alpha}\cap \widetilde f_{il}^n B(\infty)}\A G(b).
\end{equation}
Since $U^-_{\A}(\g)_{-\alpha}=\sum_{(i,l)\in I^{\infty}}(\mathtt b_{il}U^-_q(\g))^{\A}_{-\alpha}$, we obtain $$U^-_{\A}(\g)_{-\alpha}=\sum_{b\in  B(\infty)_{-\alpha}}\A G(b).$$
Let $E=L(\infty)_{-\alpha}/q L(\infty)_{-\alpha}$ and $M=U^-_{\A}(\g)_{-\alpha}$.   By Lemma \ref{K1} (ii), we have
$$ L(\infty)_{-\alpha}/q L(\infty)_{-\alpha}\xrightarrow{G} U^-_{\A}(\g)_{-\alpha}\cap  L(\infty)\cap  L(\infty)^-\rightarrow \frac{U^-_{\A}(\g)_{-\alpha}\cap  L(\infty)}{U^-_{\A}(\g)_{-\alpha}\cap q L(\infty)}$$
are isomorphisms, which proves $X(r)$.
\end{proof}
$Z(r)$ follows from (\ref{G(b)}). Finally, we will prove $B(r)$.

\vskip 3mm

\begin{lemma}\label{6.11}{
Let $\alpha\in Q_+(r)$, $b\in \mathcal B(\infty)_{-\alpha}$ and $\lambda\in P^+$. If $\overline{\pi}_\lambda (b)=0$, then $G(b)v_\lambda=0$.
\begin{proof}
Choose $(i,l)\in I^{\infty}$ such that $\widetilde e_{il}b\neq 0$. From $\overline{\pi}_\lambda (b)=0$, we know that $bv_{\lambda}\in qL(\lambda)$. Then $G(b)v_\lambda\in (b_{il}V(\lambda))^{\A}_{\lambda-\alpha}\cap q L(\lambda)\cap  L(\lambda)^-$. By Proposition \ref{M}, we have $G(b)v_\lambda=0$.
\end{proof}
}\end{lemma}

\vskip 3mm
\begin{proposition}\label{PfBr}$(Y(r))$
For any $\alpha\in Q_+(r)$ and $\lambda\in P^{+}$, we have the following canonical isomorphism
$$V(\lambda)^{\A}_{\lambda-\alpha}\cap  L(\lambda)\cap L(\lambda)^- \xrightarrow{\sim} \frac{V(\lambda)^{\A}_{\lambda-\alpha}\cap  L(\lambda)}{V(\lambda)^{\A}_{\lambda-\alpha}\cap q L(\lambda)}\xrightarrow{\sim}  L(\lambda)_{\lambda-\alpha}/q L(\lambda)_{\lambda-\alpha}.$$

\end{proposition}
\begin{proof}

By  Lemma \ref{6.11}, we have
$V(\lambda)^{\A}_{\lambda-\alpha}=\sum_{\substack{b\in  B(\infty)_{-\alpha}\\ \overline{\pi}_\lambda(b)\neq 0}}\A G(b)v_\lambda$.
Let $$E=\sum_{\substack{b\in  B(\infty)_{-\alpha}\\ \overline{\pi}_\lambda(b)\neq 0}}\Q G(b)v_\lambda\subseteq V(\lambda)^{\A}_{\lambda-\alpha}\cap  L(\lambda)\cap L(\lambda)^-.$$
Since $\{b\in B(\infty)_{-\alpha}\mid \overline{\pi}_\lambda(b)\neq 0\}\xrightarrow{\sim} B(\lambda)_{\lambda-\alpha}$ and (\ref{tg}), we have an isomorphism
$$E\xrightarrow{\sim}  L(\lambda)_{\lambda-\alpha}/q L(\lambda)_{\lambda-\alpha}$$ given by
$$G(b)v_\lambda\longmapsto G(b)v_\lambda+q L(\lambda)=\overline{\pi}_\lambda(b).$$
By Lemma \ref{K1} (ii), 
$$E\xrightarrow{\sim} V(\lambda)^{\A}_{\lambda-\alpha}\cap  L(\lambda)\cap L(\lambda)^- \xrightarrow{\sim} \frac{V(\lambda)^{\A}_{\lambda-\alpha}\cap  L(\lambda)}{V(\lambda)^{\A}_{\lambda-\alpha}\cap q L(\lambda)}\cong  L(\lambda)_{\lambda-\alpha}/qL(\lambda)_{\lambda-\alpha}$$
are isomorphisms, which proves $Y(r)$.
\end{proof}
\vskip 3mm

To summarize, we  obtain the main goal of this paper.

\vskip 3mm

\begin{theorem}\label{A}{
There exist canonical isomorphisms
\begin{equation*}
\begin{aligned}
U^-_{\A}(\g)\cap  L(\infty)\cap L(\infty)^-
& \xrightarrow{\sim}  L(\infty)/q L(\infty),\\
V(\lambda)^{\A}\cap  L(\lambda)\cap L(\lambda)^-
& \xrightarrow{\sim}  L(\lambda)/q L(\lambda) \ \ \ (\lambda \in P^{+}).
\end{aligned}
\end{equation*}
}
\end{theorem}

\vskip 3mm

\begin{definition} \hfill

\vskip 2mm

{\rm

(a) ${\mathbf B} := \{ G_\infty(b) \mid b \in B(\infty) \}$ is called the {\it global basis} of $U_{\A}^{-}(\g)$
corresponding to $B(\infty)$.

\vskip 2mm

(b) ${\mathbf B}^{\lambda} := \{ G_{\lambda}(b) \mid b \in B(\lambda) \}$ is called the {\it global basis} of $V(\lambda)^{\A}$ corresponding to $B(\lambda)$.

}
\end{definition}

\vskip 3mm

\begin{remark}

The global bases ${\mathbf B}$ and ${\mathbf B}^{\lambda}$ are unique because they are stable under the bar involution. our global bases are   pseudo-bases in the sense of Definition \ref{def:crystal basis Blambda} and Definition \ref{4.15}, but not a genuine
basis in general.

\end{remark}

\end{document}